\documentclass[english,3p,preprint,review]{elsarticle}
\usepackage[T1]{fontenc}
\usepackage[latin9]{inputenc}
\usepackage{verbatim}
\usepackage{prettyref}
\usepackage{float}
\usepackage{units}
\usepackage{amstext}
\usepackage{amssymb}
\usepackage{graphicx}
\usepackage{esint}

\makeatletter

\providecommand{\tabularnewline}{\\}

\usepackage{lineno}
\usepackage{hyperref}

\makeatother

\usepackage{babel}
\begin{document}

\title{Spatially Heterogeneous Biofilm Simulations using an Immersed Boundary
Method with Lagrangian Nodes Defined by Bacterial Locations}

\author[rvt]{Jason F. Hammond}

\ead{jason.hammond@kirtland.af.mil}

\author[focal2]{Elizabeth J. Stewart}

\ead{ejstewar@umich.edu}

\author[focal1]{John G. Younger}

\ead{jyounger@med.umich.edu}

\author[focal2]{Michael J. Solomon}

\ead{mjsolo@umich.edu}

\author[cu]{David M. Bortz\fnref{fn1}\corref{cor1}}

\fntext[fn1]{Math Biology Group, Applied Mathematics, 526 UCB, University of Colorado,
Boulder, CO 80309-0526, Phone: (303) 492-7569, Fax: (303) 492-4066
}

\ead{dmbortz@colorado.edu}

\ead[url]{http://mathbio.colorado.edu}

\cortext[cor1]{Corresponding author}

\address[rvt]{AFRL, High Power Microwave Division}

\address[cu]{Department of Applied Mathematics, University of Colorado Boulder}

\address[focal1]{Department of Emergency Medicine, University of Michigan Ann Arbor }

\address[focal2]{Department of Chemical Engineering, University of Michigan Ann Arbor}
\begin{abstract}
In this work we consider how surface-adherent bacterial biofilm communities
respond in flowing systems. We simulate the fluid-structure interaction
and separation process using the immersed boundary method. In these
simulations we model and simulate different density and viscosity
values of the biofilm than that of the surrounding fluid. The simulation
also includes breakable springs connecting the bacteria in the biofilm.
This allows the inclusion of erosion and detachment into the simulation.
We use the incompressible Navier-Stokes (N-S) equations to describe
the motion of the flowing fluid. We discretize the fluid equations
using finite differences and use a geometric multigrid method to solve
the resulting equations at each time step. The use of multigrid is
necessary because of the dramatically different densities and viscosities
between the biofilm and the surrounding fluid. We investigate and
simulate the model in both two and three dimensions. 

Our method differs from previous attempts of using IBM for modeling
biofilm/flow interactions in the following ways: the density and viscosity
of the biofilm can differ from the surrounding fluid, and the Lagrangian
node locations correspond to experimentally measured bacterial cell
locations from 3D images taken of \emph{Staphylococcus epidermidis}
in a biofilm. 
\end{abstract}
\maketitle

\section{\label{sec:Introduction-to-Biofilm}Introduction}

In this paper we investigate the response and fragmentation of a biofilm
attached to the interior of a tube and subjected to a flowing fluid.
Specifically, we study the mechanisms of biofilm fluid response and
detachment in terms of varying biofilm density, elasticity, and viscosity.
In the simulations, we model biofilms attached to the walls of 3-dimensional
square tubes using an extension of the immersed boundary method (IBM)
(originally developed by Peskin \citep{Peskin1977}). Our approach
differs from the traditional IBM in several ways. We use experimentally
measured biofilm bacterial cell locations as initial positions for
our Lagrangian nodes whereas traditional IBM refines the Lagrangian
mesh along with the Eulerian mesh. As a result we also have to adapt
the Dirac delta approximation to scale with the radius of the bacteria
rather than with the mesh width. 

In this introduction, we first provide a brief background on the biology
and biomechanics of bacterial biofilms in \prettyref{sub:Bacterial-Biofilms}.
In \prettyref{sec:Mathematical-Models}, we discuss some alternative
mathematical models that have been used to model biofilms (along with
advantages and disadvantages). In \prettyref{sec:Immersed-Boundary-Method},
we introduce the immersed boundary method, and in \prettyref{sub:Variable-Viscosity},
we discuss the significance of including variable viscosity.

\subsection{\label{sub:Bacterial-Biofilms}Bacterial Biofilms}

Biofilms are a phenotype of bacteria that are found in health, industrial
and natural settings. In the medical field, biofilms occur on devices
such as contact lenses, catheters, and mechanical heart valves. In
industrial settings, they occur in and on water pipes, storage tanks,
ship hulls, filters, food preparation facilities, etc. In natural
settings, they can be found as slime on rocks in bodies of water or
as dental plaque on teeth. 

\begin{figure}
\begin{centering}
\includegraphics[scale=0.3]{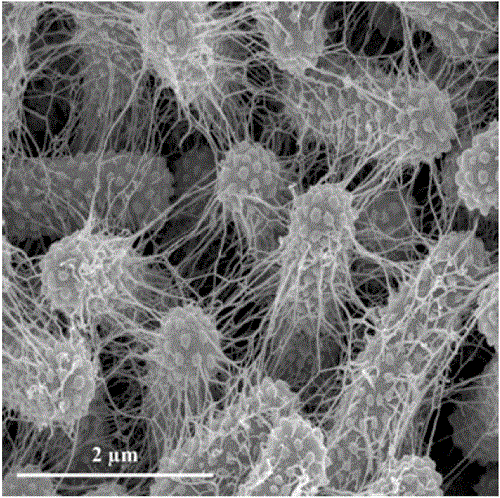}
\par\end{centering}

\caption{\label{fig:BiofilmPicture}Scanning electron microscopy of sessile
K. pneumoniae LM21 performed in mature biofilm formed on Thermanox
slides in microfermentor system after 48 hours of development at 20000
times magnification. This image is from Balestrino et al. \citep{Balestrino2008}
(used with permission).}
\end{figure}

Physically, biofilms are immobile and consist of a community of bacterial
cells embedded in a dense surface-adherent extracellular matrix (ECM)
of polysaccharides. Biofilms are mechanically strong structures that
tend to deform and fragment rather than completely dislodge when subjected
to flows. \prettyref{fig:BiofilmPicture} contains an electron micrograph
of a biofilm of \textit{Klebsiella pneumoniae}, clearly showing the
ECM interconnecting the bacterial cells. The physical properties of
the ECM are central to the growth, attachment, and detachment of biofilms.
The focus of this work is on the biomechanical response of the ECM
to fluid flow resulting in deformation and separation.  

An important feature of biofilms is that they are known to behave
like \textit{viscoelastic fluids} \citep{Klapper2002}. In other words,
they exhibit both viscous and elastic responses upon deformation.
Describing the exact viscoelastic behavior has been the subject of
much experimental, theoretical, and computational research \citep{Klapper2002,Rupp2005,Aravas2008,Pavlovsky2013,Lau2009,Kreft2001,Picioreanu2004,Picioreanu1999,PicioreanuLoodsdrecht2000DisContMod,Picioreanu2001}.
For example, Klapper et al.~ and Pavlosky et al.~use a linear Jeffrey\textquoteright{}s
constitutive law \citep{Klapper2002,Pavlovsky2013} while Lau et al.~use
a Voigt standard linear solid model for viscoelastic materials \citep{Lau2009}.
Our model includes elasticity of the biofilm by using simple linear
springs (as done in \citep{Alpkvist2007}) to connect the bacterial
cells, and includes the viscosity of the biofilm with a modification
of the constitutive equations for stress.

\subsection{Mathematical Models of Biofilms\label{sec:Mathematical-Models}}

Much research into the mathematical modeling of biofilm growth and
fluid/structure interactions has been conducted in the last three
decades \citep{Picioreanu1999,PicioreanuLoodsdrecht2000DisContMod,Picioreanu2001,ZhangCoganWang2008One,ZhangCoganWang2008Two,Kreft1998,Kreft2001}.
Below, we summarize several modeling and simulation strategies. This
is not intended to be exhaustive and we direct the interested reader
to \citep{WangZhang2010Review,Klapper2010} for more in-depth reviews. 

The first attempts at mathematical modeling of biofilms were conducted
in the early 1980s \citep{Rittmann1980,Rittmann1982,Kissel1984}.
Picioreanu and others \citep{Kreft2001,Picioreanu2004,Kreft1998}
advocated for an individual based (IB) approach, which models the
behavior of each bacteria, encompassing ideas such as cell division,
cell motility, metabolism, and death to simulate the growth and formation
of colonies. Hybrid discrete-continuum models were the first methods
to couple the flow with the the biofilm computationally in 2D and
3D simulations. Picioreanu, van Loodsdrecht, and Heijnen developed
and used these hybrid discrete-continuum models to incorporate the
flow over the irregular biofilm's surfaces, convective and diffusive
mass transfer of substrate, bacterial growth, and biomass spreading
\citep{Picioreanu1999,PicioreanuLoodsdrecht2000DisContMod,Picioreanu2001}. 

The most sophisticated (purely) continuum models developed are the
phase-field models, which use a one-fluid/two-component formulation
in which the ECM and the bacteria are modeled as one fluid component,
while the collective ensemble of nutrient substrates and the surrounding
fluid are the other \citep{ZhangCoganWang2008One,ZhangCoganWang2008Two}.
Two-dimensional simulations of both biofilm growth and biofilm-flow
interaction are presented in \citep{ZhangCoganWang2008Two}, in which
shear induced deformation and detachment are illustrated. 

We note that our model differs from both these approaches in the way
that we model the biofilm. We treat the bacteria in the biofilm as
discrete points, where the nodal locations in our simulations correspond
to the locations of the bacterial cells within the biofilm. This contrasts
from the continuum phase-field models that only include averaged biomechanical
properties of the biofilm. With our mathematical formulation, just
as in the individual based models, we can obtain the cumulative local
stresses as well as attribute different local properties to the biofilm.
Our model can be thought of as an extension of the individual based
models, where we accurately account for the interactions with the
fluid as well as include the possibility of fragmentation. We also
assume that on the time scale of our simulations there is no biofilm
growth; so we ignore such factors as nutrient concentrations and growth
rates.

\subsection{\label{sec:Immersed-Boundary-Method}Immersed Boundary Method}

The overall goal is to simulate the response of a biofilm attached
to the walls of both 2- and 3-dimensional square tubes. We do this
using an extension of the immersed boundary method. In this section,
we introduce the immersed boundary method and provide the framework
necessary for us to later extend the IBM for our use. 

The immersed boundary method (IBM) was originally developed by Peskin
to study blood flow in the heart \citep{Peskin1977}. The IBM has
been used previously to model and simulate biofilm/fluid interactions
by Dillon, Fauci, et al.~in \citep{Dillon1996}, and by Alpkvist
and Klapper in \citep{Alpkvist2007}. The authors successfully coupled
the fluid to the biofilm; however, they make the assumption that the
biofilm has the same density and viscosity as the surrounding fluid.
This choice substantially simplifies the task of solving the N-S equations
but does not account for the fact that biofilms typically have $500\times$
larger viscosity and $12\%$ larger density than water \citep{Klapper2002,Ro1991}.
They also use a random distribution of points within a biofilm-shaped
shell to represent the biofilm, which does not account for the true
spatial distribution of bacterial cells within a biofilm. 

The IBM has been used more recently in the modeling of immersed elastic
structures in viscous flows in \citep{Huang2009IBM,Luo2008,Zhuo2011,Strychalski2012},
in which the authors use constitutive viscoelastic models including
Maxwell, Voigt, and Jeffrey's models to incorporate forces in the
immersed structures into the IBM. Similarly, our ultimate goal is
to establish an appropriate constitutive model for the forces in the
biofilm with the help of experimental collaborators and to include
this in our IBM formulation.

\subsection{\label{sub:Variable-Viscosity}Variable Viscosity}

It is agreed upon in the biofilm research community that biofilms
behave like a viscoelastic fluid. There has been a few efforts to
match the behavior of biofilms with conventional mechanical viscoelastic
models \citep{Klapper2002,Aravas2008,Pavlovsky2013,Lau2009}. However,
these efforts have not produced a consensus on the model to use for
viscoelasticity in biofilms. This is, in part, due to the fact that
the viscoelastic properties in biofilms is highly variable with different
growth conditions \citep{Chen2005effectsAdhesiveStrength} and even
in the same growth conditions \citep{Aggarwal2010}. 

The incorporation of spatially variable viscosity in the immersed
boundary method is an area that has yet to be well developed. Luo
et al.~couple the immersed viscoelastic structure to the fluid flow
in an immersed boundary type formulation, but they solve the fluid
equations separately from the equations governing the motion of the
immersed viscoelastic solid and then couple the solutions at their
physical interface \citep{Luo2008}. This formulation will not work
in our case because we couple the biofilm to the fluid in the entire
domain so that it will behave as a viscoelastic fluid. Another approach
to including viscosity into the immersed boundary method is by replacing
the simple elastic springs with viscoelastic links, which will change
the value of the external force, $\mathbf{f}$, in the Navier-Stokes
equations (Equation \prettyref{eq:N-S1} below). This type of strategy
was used first by Bottino in \citep{Bottino1998} to model general
viscoelastic connections in actin cytoskeleton of ameboid cells and
also by Dillon and Zhuo in \citep{Zhuo2011} to model sperm motility. 

There are two natural ways to add viscosity to the biofilm. The first
is through the use of local damping forces in addition to the springs
that we use for the elastic component. In this way, we can define
the forces between any two connected bacterial cells using a typical
mechanical model for a viscoelastic material. The second way is to
treat the entire domain as a continuous Newtonian viscous fluid with
a spatially varying viscous coefficient. The use of dashpot damping
with our current mathematical formulation is quicker to implement
but has serious stability restrictions in the simulation, and thus
we omit this method here.%
\footnote{We will pursue this approach in future work when change our numerical
scheme to a semi-implicit or implicit method.%
} 

The approach in this work is to treat the fluid in the entire domain
as a Newtonian viscous fluid with a spatially varying viscous coefficient.
The core idea involves replacing the viscous term in the Navier-Stokes
equation with one that can apply to nonuniform viscosity in the fluid.
We note that to the best of our knowledge this approach has not yet
been attempted in the single fluid immersed boundary method. There
have been attempts using two materials (fluid-fluid or fluid-solid),
coupling them at their interface, in which the stress is adapted in
the viscoelastic fluid or solid to account for a different viscosity
\citep{ZhangCoganWang2008Two,Luo2008}. However, in our approach,
we couple the biofilm to the fluid within the entire biofilm region,
not just at the interface. Thus, we must adapt the forces in the Navier-Stokes
equations to account for the viscous and elastic stresses on the surface
and within the biofilm.

We now describe the organization of this paper. In \prettyref{sec:Mathematical-Formulation}
we provide our mathematical formulation, which is a variation of the
immersed boundary method literature. In \prettyref{sec:Numerical-Method},
we describe our numerical method, based on a multigrid approach. In
\prettyref{sec:Validation}, we provide numerical validation of our
method. In \prettyref{sec:Simulations-Results}, we provide simulation
results in both two and three dimensions, running our simulations
for a variety of experimentally obtained biofilms with varying parameters
such as spring constants, densities, viscosities. Finally, we provide
conclusions in \prettyref{sec:Biofilm-Simulation-Conclusions} and
a discussion of possibilities for future work in \prettyref{sec:FW}.

\section{Mathematical Formulation\label{sec:Mathematical-Formulation}}

In this section, we provide the mathematical formulation for our simulations.
We use an Eulerian mesh to describe the system as a whole and solve
the dimensionless N-S equations at each time step on this mesh. The
Lagrangian nodes are used only to compute information about the biofilm
(location, velocity, local density, force) and then transfer the information
back onto the Eulerian mesh using the Dirac delta function.%
\footnote{The Dirac Delta function is approximated in the actual implementation
(see Equation \prettyref{eq:Dirac Delta Approx}).%
}\negthinspace{}\negthinspace{} For convenience, we provide a list
in Appendix A of the variables and parameters used in this work. 

We now introduce the mathematical equations used in our model. The
dependent Eulerian variables are velocity $\mathbf{u}(\mathbf{x},\, t)$,
pressure $p(\mathbf{x},\, t)$, density $\rho(\mathbf{x},\, t)$,
and Eulerian force density $\mathbf{f}(\mathbf{x},\, t)$, where $\mathbf{x}$
is the independent Eulerian variable and $t$ is time. The dependent
Lagrangian variables are position of the nodes $\mathbf{X}(\mathbf{q},\, t)$,
velocity of the nodes $\mathbf{U}(\mathbf{q},\, t)$, and the Lagrangian
force density $\mathbf{F}(\mathbf{q},\, t)$, where $\mathbf{q}=(q,r,s)$
is the independent Lagrangian variable. The equations of motion for
the biofilm-fluid interaction are 
\begin{eqnarray}
\rho(\mathbf{x},\, t)\left(\frac{\partial\mathbf{u}}{\partial t}+\mathbf{u}\cdot\triangledown\mathbf{u}\right) & = & -\triangledown p+\triangledown\cdot\left(\mu(\mathbf{x},\, t)\left(\triangledown\mathbf{u}+\left(\triangledown\mathbf{u}\right)^{T}\right)\right)+\mathbf{f}(\mathbf{x},\, t)\label{eq:N-S1VarVisc}\\
\triangledown\cdot\mathbf{u} & = & 0,\label{eq:Incompressible}\\
\frac{\partial\mathbf{X}}{\partial t}(\mathbf{q},\, t) & = & \mathbf{U}(\mathbf{X}(\mathbf{q},\, t),\, t),\label{eq:LagVel}\\
\mathbf{f}(\mathbf{x},\, t) & = & \int_{\Omega_{b}}\mathbf{F}(\mathbf{q},\, t)\delta(\mathbf{x}-\mathbf{X}(\mathbf{q},\, t))\textrm{d}\mathbf{q},\label{eq:ForceCouple}\\
\mathbf{\rho}(\mathbf{x},\, t) & = & \rho_{0}+\int_{\Omega_{b}}\rho_{b}\delta(\mathbf{x}-\mathbf{X}(\mathbf{q},\, t))\textrm{d}\mathbf{q},\label{eq:DensCouple}\\
\mathbf{U}(\mathbf{X}(\mathbf{q},\, t),\, t) & = & \int_{\Omega}\mathbf{u}(\mathbf{x},\, t)\delta(\mathbf{x}-\mathbf{X}(\mathbf{q},\, t))\textrm{d}\mathbf{x},\label{eq:VelCouple}
\end{eqnarray}
where $\mu$ is the dynamic viscosity, $\rho_{0}$ is the mass density
of the fluid, $\rho_{b}$ is the additional mass density of the biofilm
from that of the surrounding fluid, $\Omega$ is the flow domain,
$\Omega_{b}\subset\Omega$ is the space occupied by only the biofilm,
and $\delta(\mathbf{x})$ is the Dirac delta function. Equations \prettyref{eq:N-S1VarVisc}
and \prettyref{eq:Incompressible} are the incompressible Navier-Stokes
(N-S) equations with spatially varying viscosity and a forcing term
that represents the forces applied by the biofilm on the fluid. Equation
\prettyref{eq:LagVel} is the equation of motion of the biofilm, where
$\mathbf{U}(\mathbf{q},\, t)$ is the velocity of the biofilm. The
systems of PDE's given by \prettyref{eq:N-S1VarVisc}-\prettyref{eq:Incompressible}
is coupled to \prettyref{eq:LagVel} by the integrals given in \prettyref{eq:ForceCouple}-\prettyref{eq:VelCouple}.

To avoid numerical inaccuracies due to roundoff errors, we non-dimensionalize
these equations using the non-dimensional variables defined as 

\begin{tabular}{cccc}
$t^{*}=\frac{t}{T}$,  & $\mathbf{x}^{*}=\frac{\mathbf{x}}{L}$,  & $\mathbf{u}^{*}=\frac{\mathbf{u}}{u_{0}}$, & $p^{*}=\frac{p-p_{L_{tube}}}{p_{0}-p_{L_{tube}}}$,\tabularnewline
$\mathbf{\triangledown}^{*}=L\mathbf{\triangledown}$, & $\rho^{*}=\frac{\rho}{\rho_{0}}$, & $\mathbf{f}^{*}=\frac{\mathbf{f}}{f_{0}}$, & $\mu^{*}=\frac{\mu}{\mu_{0}}$,\tabularnewline
\end{tabular} \linebreak{}
where $p_{0}$ is the pressure at the upstream end of the tube, $p_{L_{tube}}$
is the pressure at the downstream end of the tube, $T$ is the characteristic
time scale, $f_{0}$ is the characteristic force density, and $L$
is the characteristic length. We use the scaling parameters defined
in \prettyref{tab: 2D sim params-1}. Dropping the stars from the
dimensionless variables, equations \prettyref{eq:Incompressible}
and \prettyref{eq:ForceCouple}-\prettyref{eq:VelCouple} remain the
same as in the case with dimensions, while equations \prettyref{eq:N-S1VarVisc}
and \prettyref{eq:LagVel} become
\begin{eqnarray}
\sigma\rho(\mathbf{x},\, t)\frac{\partial\mathbf{u}}{\partial t}+\rho(\mathbf{x},\, t)\mathbf{u}\cdot\triangledown\mathbf{u} & = & -\varepsilon\triangledown p+Re^{-1}\triangledown\cdot\left(\mu(\mathbf{x},\, t)\left(\triangledown\mathbf{u}+\left(\triangledown\mathbf{u}\right)^{T}\right)\right)+\frac{Lf_{0}}{\rho_{0}u_{0}^{2}}\mathbf{f}(\mathbf{x},\, t),\label{eq:N-SnondimFull}\\
\sigma\frac{\partial\mathbf{X}}{\partial t}(\mathbf{q},\, t) & = & \mathbf{U}(\mathbf{q},\, t),
\end{eqnarray}
where $\sigma=\frac{L}{Tu_{0}}$ is the Strouhal number, $\varepsilon=\frac{p_{0}-p_{L_{tube}}}{\rho_{0}u_{0}^{2}}$
is the Euler number, and $Re=\frac{\rho_{0}Lu_{0}}{\mu}$ is the Reynolds
number of the fluid. 

The initial velocity profile is the exact solution to the incompressible
Navier-Stokes equations in a square or circular tube with rigid walls
and no-slip conditions at the walls. The velocity profile for a circular
cylinder can be found in many textbooks in fluid dynamics (such as
\citep{Zamir2000}), and a series solution for the laminar flow velocity
profile for a square tube was derived by Spiga and Morini in \citep{Spiga1994}.

\section{Numerical Method\label{sec:Numerical-Method}}

In this section, we describe the numerical formulation for our simulations.
Our numerical task is to solve the system defined by Equations \prettyref{eq:N-S1VarVisc}-\prettyref{eq:VelCouple},
and we now provide the details of our numerical approach. 

The incompressible flow Navier-Stokes equations, \prettyref{eq:N-S1VarVisc}-\prettyref{eq:Incompressible},
are discretized on a fixed uniform Eulerian lattice, while the biofilm
equations are discretized on a moving Lagrangian array of points that
do not necessarily coincide with the fixed Eulerian mesh points of
the fluid computation. We represent the interaction equations \prettyref{eq:ForceCouple}-\prettyref{eq:VelCouple}
with a smoothed approximation $\tilde{\delta}$ to the Dirac delta
function (see \prettyref{sub:Dirac-Delta-Approximation}). Our numerical
approach was inspired by the solving technique used by Zhu and Peskin
in \citep{Zhu2002} to simulate a flapping filament in a soap film. 

The discretized equations corresponding to \prettyref{eq:ForceCouple}-\prettyref{eq:VelCouple}
are given by
\begin{eqnarray}
\mathbf{f}^{n}(\mathbf{x}) & = & \sum_{s=1}^{\eta}\mathbf{F}^{n}(s)\tilde{\delta}(\mathbf{x}-\mathbf{X}^{n}(s),\,\omega),\label{eq:discrForceCouple}\\
\mathbf{\rho}^{n}(\mathbf{x}) & = & \rho_{0}+\sum_{s=1}^{\eta}\rho_{b}\tilde{\delta}(\mathbf{x}-\mathbf{X}^{n}(s),\,\omega)d_{0}^{3},\label{eq:discrDensCouple}\\
\mathbf{U}^{n+1}(s) & = & \sum_{\mathbf{x}}\mathbf{u}^{n+1}(\mathbf{x})\tilde{\delta}(\mathbf{x}-\mathbf{X}^{n}(s),\, h)h^{3},\label{eq:discrVelCouple}
\end{eqnarray}
where the superscript $n$ denotes numerical approximations at a particular
time step $n$, $\eta$ is the total number of Lagrangian discretization
points, the sum in \prettyref{eq:discrVelCouple} is over all the
discrete points of the form $\mathbf{x}=(ih,\, jh,\, kh)$ with $i,\, j,$
and $k$ are integers, $h$ is the Eulerian mesh width, and $d_{0}^{3}$
is the average volume element of the Lagrangian nodes (computed by
dividing the total volume of the biofilm by the total number of Lagrangian
nodes distributed within it). Following convention, we replace $(q,\, r,\, s)$
from the mathematical formulation with only $s$, which we use as
an indexed label with a unique number assigned to each Lagrangian
point \citep{Zhu2002}. In \prettyref{eq:discrForceCouple}, $\mathbf{F}(s)$
is now the total elastic force on the Lagrangian node associated with
marker $s$, as opposed to an elastic force density. This is because
we calculate the force explicitly depending on which other nodes it
is connected to.

\subsection{\label{sub:Dirac-Delta-Approximation}Dirac Delta Approximation}

In \citep{Peskin2002}, Peskin defines $\delta_{h}(\mathbf{x})$ as
\begin{eqnarray}
\delta_{h}(\mathbf{x})=h^{-3}\phi\left(\frac{x}{h}\right)\phi\left(\frac{y}{h}\right)\phi\left(\frac{z}{h}\right) &  & ,\label{eq:Dirac Delta Post 1}
\end{eqnarray}
where $\phi(r)$ is 
\begin{equation}
\phi(r)=\phi_{1}=\left\{ \begin{array}{cc}
\frac{1}{8}\left(3-2|r|+\sqrt{1+4|r|-4r^{2}}\right); & \textrm{if }|r|\le1\,,\\
\frac{1}{8}\left(5-2|r|-\sqrt{-7+12|r|-4r^{2}}\right); & \textrm{if }1\le|r|\le2\,,.\\
0; & \textrm{if }|r|\ge2\,.
\end{array}\right.\label{eq:PhiPeskin}
\end{equation}

We replace this $\delta_{h}$, that is used in standard IBM implementations,
with one that scales with $\omega$ instead of $h$ as
\begin{equation}
\tilde{\delta}(\mathbf{x},\,\omega)=\omega^{-3}\phi\left(\frac{x}{\omega}\right)\phi\left(\frac{y}{\omega}\right)\phi\left(\frac{z}{\omega}\right).\label{eq:Dirac Delta Approx}
\end{equation}
We deviate from the standard scaling of the Dirac Delta approximation
for two reasons. The first is that we wish to give a presence to the
bacterial cells that is representative of the true volume of the cells.
Thus, in the simulations, we make $\omega$ in \prettyref{eq:discrForceCouple}
and \prettyref{eq:discrDensCouple} equal to the radius of a bacterial
cell that we are modelling. Equation \prettyref{eq:discrForceCouple}
then spreads the force over a volume that is slightly larger than
the cell, ensuring that the entire space occupied by the cell in the
fluid is influenced by the force. The second reason we use this scaling
is because, during the mesh refinement analysis described in \prettyref{sub:Eulerian-Grid-Refinement},
we discovered that the implementation with the scaling by $h$ restricts
us to less than first-order convergence of the velocity, $\mathbf{u}$.
Using a scaling that is independent of the mesh-width fixes this issue
and leads to greater than first order convergence. 

\begin{figure}
\begin{centering}
\includegraphics[width=4in]{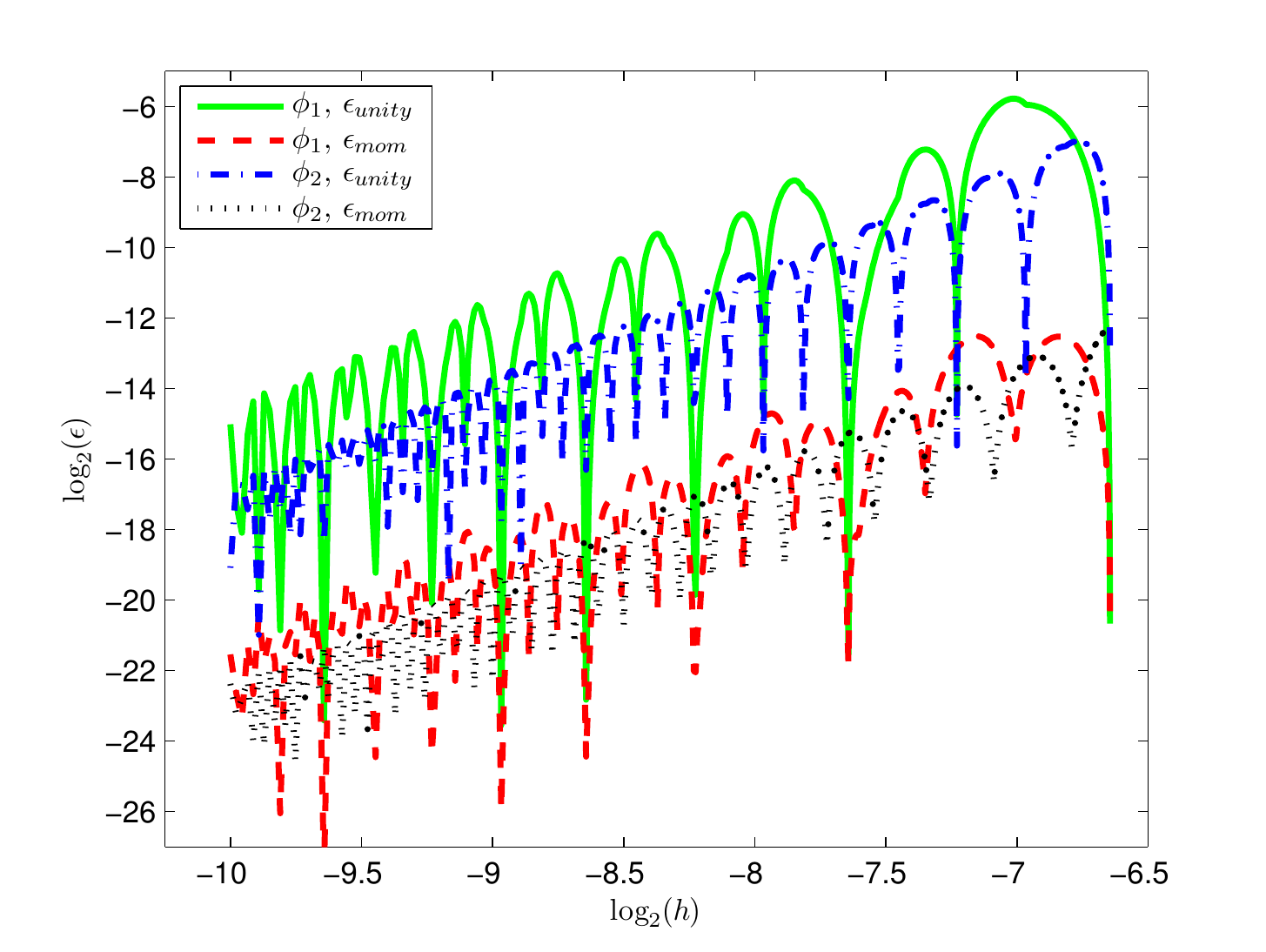}
\par\end{centering}

\caption{\label{fig:ComarePhi}We show $\epsilon_{unity}(\nicefrac{1}{100},h)$
and $\epsilon_{mom}(\nicefrac{1}{100},h)$ for $h\in\left(\frac{1}{1024},\frac{1}{100}\right)$.
$\phi_{1}$ is the $\phi$ given in \prettyref{eq:PhiPeskin} and
$\phi_{2}$ is the $\phi$ given in \prettyref{eq:PhiCosAlternative}.
We show $\log_{2}$ in the $x$ and $y$ axes so that the convergence
rate appears as the slope of the line segments. }
\end{figure}

However, this modification does have the negative consequence of losing
two desirable conditions that were previously satisfied by the Dirac
delta approximation, $\delta_{h}$. Specifically, with $\delta_{h}$
as defined in \prettyref{eq:PhiPeskin}, the \textit{unity condition,}
\begin{equation}
\sum_{\mathbf{x}\in g_{h}}\delta_{h}(\mathbf{x}-\mathbf{X})h^{3}=1,\:\forall\mathbf{X},\label{eq:unityConditionDirac}
\end{equation}
and the \textit{first-moment condition,} 
\begin{equation}
\sum_{\mathbf{x}\in g_{h}}(\mathbf{x}-\mathbf{X})\delta_{h}(\mathbf{x}-\mathbf{X})h^{3}=0,\:\forall\mathbf{X}\,,\label{eq:1stMomentConditionDirac}
\end{equation}
are both satisfied. However, using $\tilde{\delta}$ in place of $\delta_{h}$,
these conditions fail to hold true for all $\mathbf{X}$ when $\omega\ne h$.
In practice, this is not a major concern as many IBM formulations
use a Dirac delta approximation that satisfies the unity condition,
\prettyref{eq:unityConditionDirac}, but does not satisfy the first-moment
condition, \prettyref{eq:1stMomentConditionDirac}. For example, in
\citep{Zhu2002}, Peskin and Zhu replace $\phi(r)$ in $\delta_{h}$
with
\begin{equation}
\phi(r)=\phi_{2}=\left\{ \begin{array}{cc}
\frac{1}{4}\left(1+\cos\left(\frac{\pi r}{2}\right)\right); & \textrm{if }|r|\le2,\\
0; & \textrm{if }|r|>2.
\end{array}\right.\label{eq:PhiCosAlternative}
\end{equation}
We do have to choose whether we want to use $\phi(r)$ as defined
by \prettyref{eq:PhiPeskin} or by \prettyref{eq:PhiCosAlternative}.
For either choice of $\phi$, it is true that both 
\[
\lim_{h\rightarrow0}\int\tilde{\delta}(\mathbf{x}-\mathbf{X})d\mathbf{x}=1
\]
 and 
\[
\lim_{h\rightarrow0}\int(\mathbf{x}-\mathbf{X})\tilde{\delta}(\mathbf{x}-\mathbf{X})d\mathbf{x}=0\,.
\]
Note that, in the limit as $h\rightarrow0$, we see greater than $O(h^{2})$
convergence to \prettyref{eq:unityConditionDirac} and \prettyref{eq:1stMomentConditionDirac}
(see \prettyref{fig:ComarePhi}), which is consistent with the theoretical
convergence rate for a Riemann sum. Therefore, we choose to use the
$\phi$ for which the summations in \prettyref{eq:unityConditionDirac}
and \prettyref{eq:1stMomentConditionDirac} are closest to $1$ and
$0$, respectively, for the values of $\omega$ and $h$ used in our
simulations so that we can have the most accurate discrete approximation
of the Dirac delta function.

We now define two error metrics to determine how well $\tilde{\delta}$
(using either $\phi_{1}$ or $\phi_{2}$) satisfies the unity and
first-moment conditions. Analogous metrics and comparisons could be
conducted in higher dimensions, but for simplicity we provide a one-dimensional
comparison. Using $\omega=\frac{1}{100}$ and the one-dimensional
version of \prettyref{eq:Dirac Delta Approx}, we define 
\begin{equation}
\epsilon_{unity}(\omega,h)=\max_{\mathbf{X}\in[0,h]}\left|\left(\sum_{\mathbf{x}\in g_{h}}\tilde{\delta}(\mathbf{x}-\mathbf{X},\,\omega)h\right)-1\right|\label{eq:PhiErrorUnity}
\end{equation}
 and 
\begin{equation}
\epsilon_{mom}(\omega,h)=\max_{\mathbf{X}\in[0,h]}\left|\sum_{\mathbf{x}\in g_{h}}\left(\mathbf{x}-\mathbf{X}\right)\tilde{\delta}(\mathbf{x}-\mathbf{X},\,\omega)h\right|.\label{eq:PhiErrorFirstMomentum}
\end{equation}
We only have to find the maximum over $\mathbf{X}\in[0,h]$ since
the summations are periodic with period $h$, because $\mathbf{X}=0$
and $\mathbf{X}=h$ both correspond to a Lagrangian point being at
the same location as an Eulerian point. We find $\epsilon_{unity}(\nicefrac{1}{100},h)$
and $\epsilon_{mom}(\nicefrac{1}{100},h)$ for values of $h\in\left(\frac{1}{1024},\frac{1}{100}\right)$
and using both $\phi_{1}$ and $\phi_{2}$. These values are compared
in \prettyref{fig:ComarePhi}, and show that using $\phi_{2}$ provides
a better approximation of the Dirac delta function in terms of matching
the values of these summations for most values of $h$. Notice too,
that the values of $\epsilon_{unity}$ and $\epsilon_{mom}$ are exactly
zero for certain relationships between $h$ and $\omega$. For example,
if $h=\frac{\omega}{z}$ with $z\in\mathbb{N}$, then $\epsilon_{unity}$
and $\epsilon_{mom}$ are exactly zero when using $\phi_{1}$. Thus
if our problem allows for $h=\frac{\omega}{z}$, then we would use
$\phi_{1}$. However, with the choices we have made for $\omega$
and $h$ in our numerical scheme, it is a better choice to use $\phi_{2}$
in our simulations. 

In the transfer of information from the Lagrangian grid to the Eulerian
grid, we scale $\phi$ by $\omega$ (see \prettyref{eq:discrForceCouple}
and \prettyref{eq:discrDensCouple}). However, in the transfer of
velocity from the Eulerian to Lagrangian grid (Eq. \prettyref{eq:discrVelCouple}),
we scale $\phi$ with $h$ instead of $\omega$ in order to capture
the velocity only at the center of mass of the bacterial cells. 

\subsection{Elastic Forces}

In \citep{Alpkvist2007}, Alpkvist and Klapper use Hooke's Law to
describe the elastic force between the connected Lagrangian nodes.
We also use this method as our first attempt to model the interconnecting
force in the biofilm. Thus, the elastic force on each Lagrangian point
using Hooke's Law is
\begin{equation}
\mathbf{F}^{n}(s)=\sum_{k=1}^{\eta}I_{s,k}\frac{\mathbf{d}_{s,k}}{d_{s,k}}T_{s,k}\,,\label{eq:springF}
\end{equation}
where $T$ is the tension between nodes $s$ and $k$, $I$ is the
connectivity matrix defined as 
\begin{eqnarray*}
I_{s,k} & = & \left\{ \begin{array}{cc}
1 & \textrm{ bacteria }s\textrm{ connected to bacteria }k\\
0 & \textrm{otherwise}\:,
\end{array}\right.
\end{eqnarray*}
and $\mathbf{d}_{s,k}$ is the vector pointing from Lagrangian node
$s$ to $k$ with magnitude $d_{s,k}$. The tension from the spring
connecting node $s$ and $k$ is formulated as 
\[
T_{s,k}=K_{s,k}(d_{s,k}-r_{s,k}),
\]
where $r_{s,k}$ is the rest length of the spring connecting nodes
$s$ and $k$, and $K_{s,k}$ is its Hookean spring coefficient. We
choose to define each spring coefficient as
\begin{equation}
K_{s,k}=\frac{F_{max}}{r_{s,k}},\label{eq:defineK}
\end{equation}
where $F_{max}$ is the force required to break the spring. We define
the spring coefficients in this way to ensure that all of the springs,
regardless of initial length, break with a force of $F_{max}$ when
they are stretched to a length of $2r_{s,k}$. In our simulations,
we vary $F_{max}$ to attain specific results (such as detachment;
see \prettyref{sub:Discussion-of-Elastic} more details). As is done
in \citep{Alpkvist2007}, we model the failure of the ECM by breaking
the connections between the Lagrangian nodes as the springs used to
connect them exceed twice their resting length. We note that this
condition, however, is not based on experimental evidence, and in
future work we will adapt this breaking criteria according to experimental
results. In \prettyref{sec:FW}, we discuss future adaptations to
the breaking criteria in terms of the yield stress of polymers.

\begin{figure}
\begin{centering}
(a)\includegraphics[width=5in]{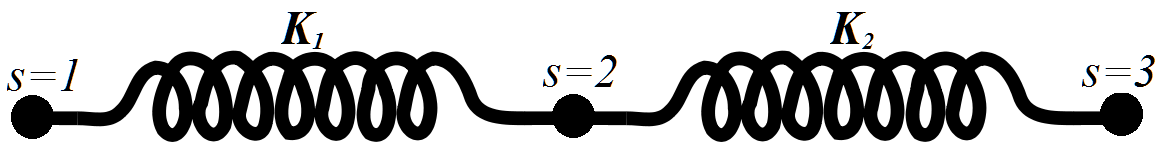}\\
(b)\includegraphics[width=3in]{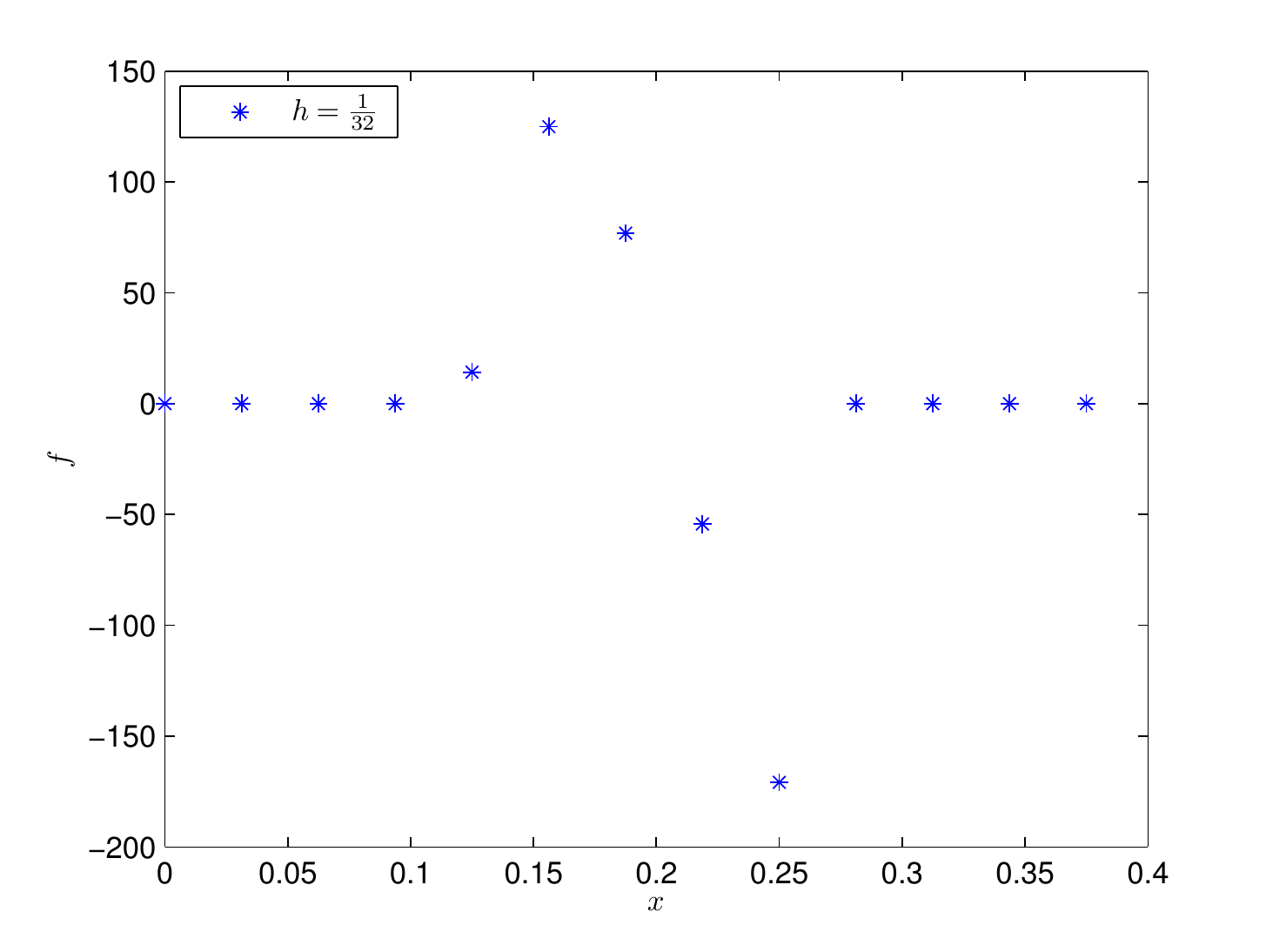}(c)\includegraphics[width=3in]{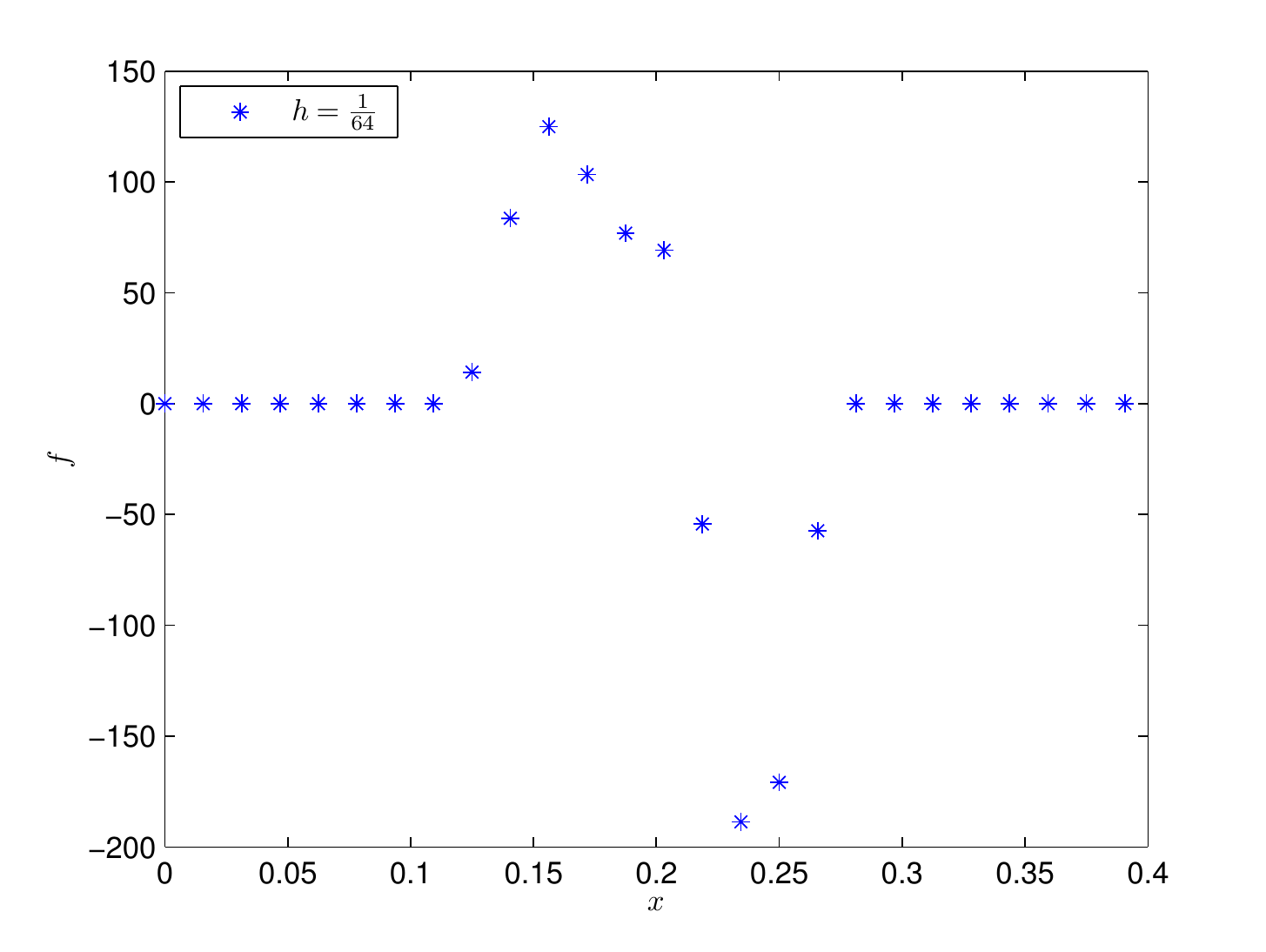}\\
(d)\includegraphics[width=3in]{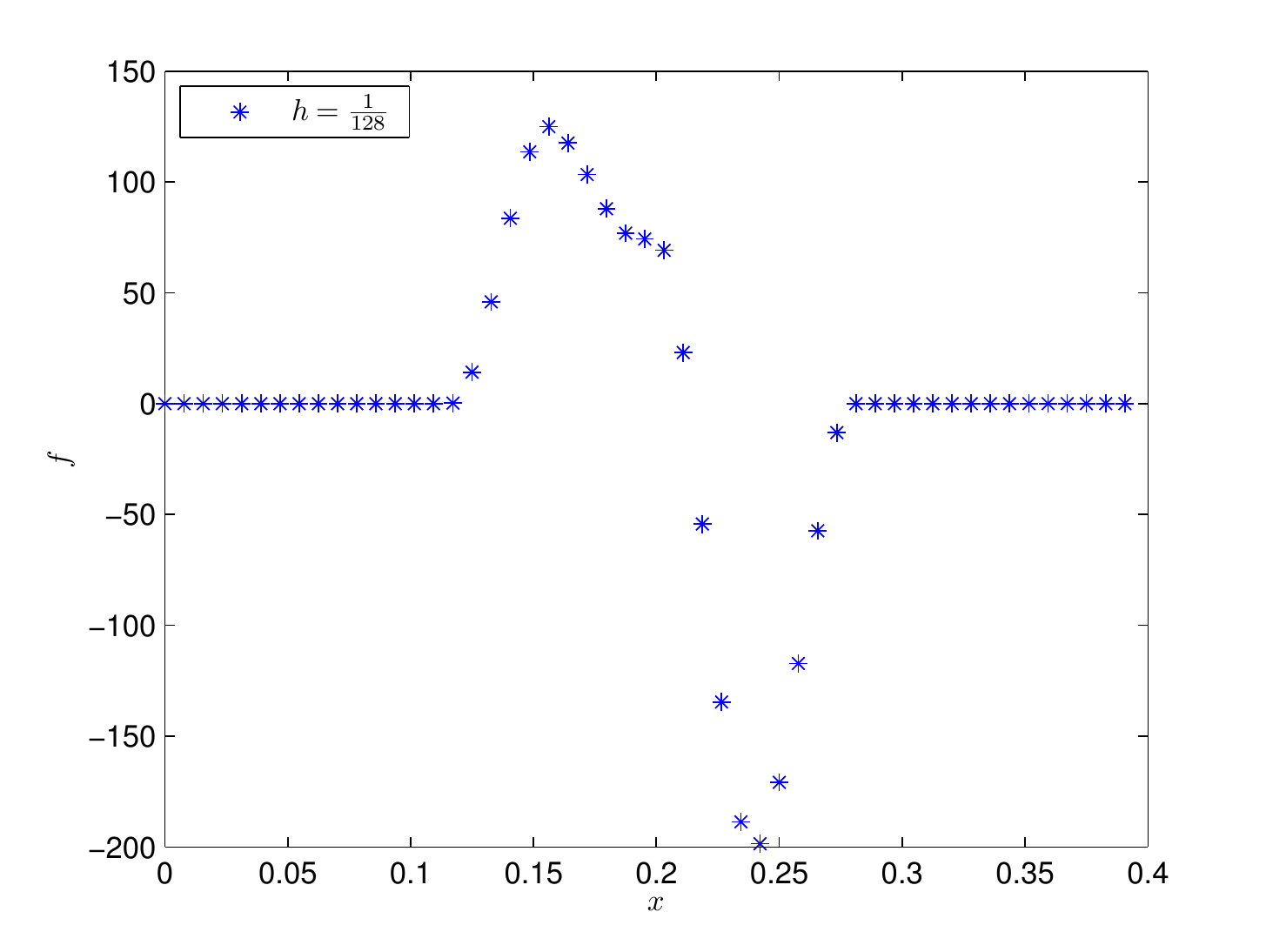}(e)\includegraphics[width=3in]{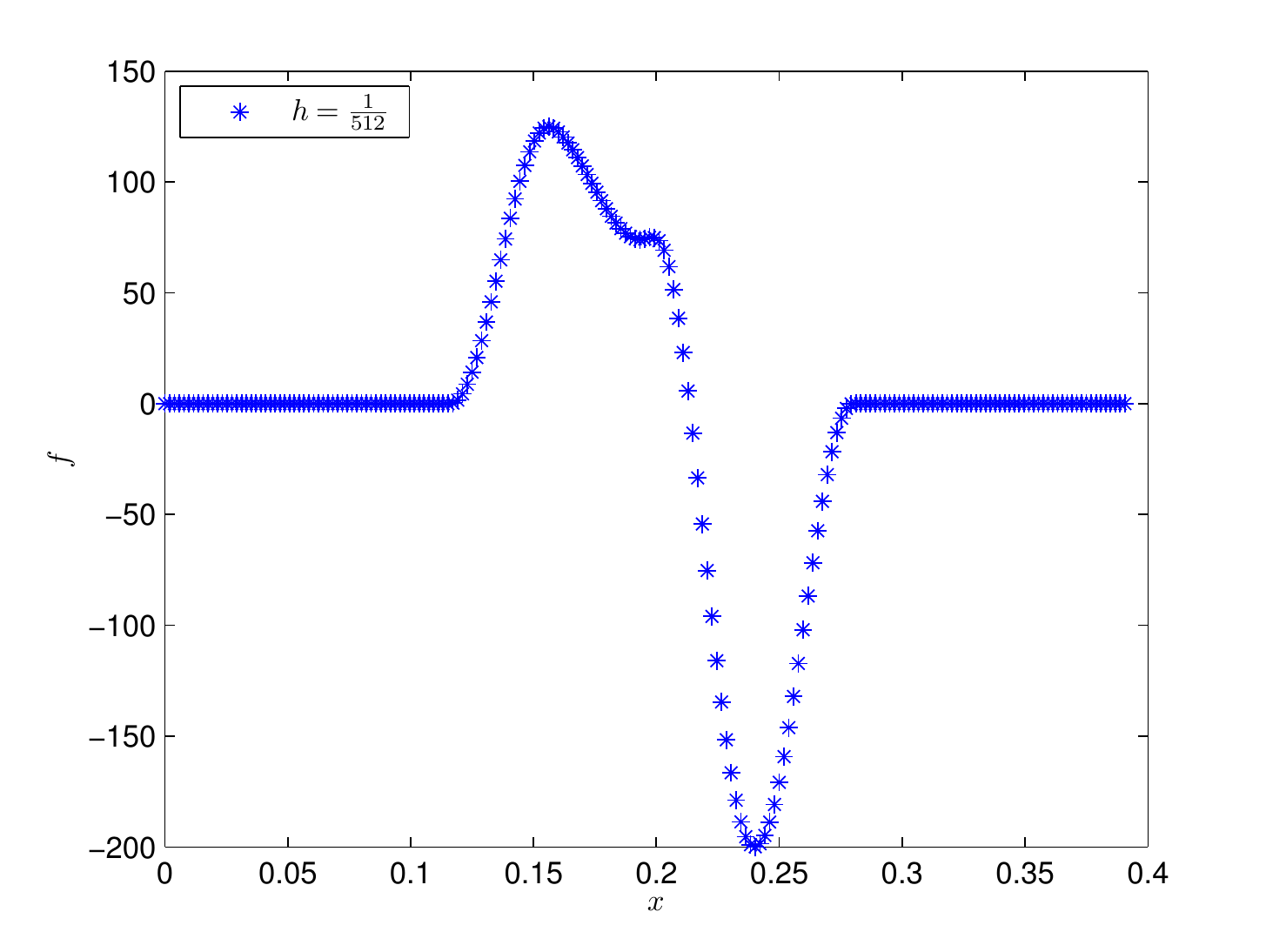}
\par\end{centering}

\caption{\label{fig:ForceEx}Examples of 1D $f(x)$ with three Lagrangian nodes
(a), one at $X(s=1)=\nicefrac{5}{32}$ on an Eulerian node and the
others at $X(s=2)=\frac{5}{32}+\frac{2\pi}{150}$ and $X(s=3)=\frac{5}{32}+\frac{4\pi}{150}$
(this position is chosen so that the second and third nodes are close
enough to the first one to show the interaction of the forces), with
$\omega=\frac{1}{50}$. These positions were chosen so that the cells
would be close enough to each other to have an overlapping region
on the Eulerian grid after the transfer of the forces from the Lagrangian
grid. In this demonstration, $F(s=1)=5,\, F(s=3)=-8,$ and $F(s=2)=-(F(1)+F(3))=3$.
These plots illustrate the effect of using different spatial steps:
(b) $h=\frac{1}{32}$, (c) $h=\frac{1}{64}$, (d) $h=\frac{1}{128}$,
(e) $h=\frac{1}{512}$. }
\end{figure}

We conclude this section with a 1D illustration to show how 3 linearly
connected cells transfer their elastic forces to the Eulerian grid
via Equation \prettyref{eq:discrForceCouple}, and how this is effected
by varying $h$ (see \prettyref{fig:ForceEx}). In this example, the
three nodes are close enough to each other that their forces add in
the overlapping regions. Sub-figures (b)-(e) illustrate how, with
finer discretizations (smaller $h$), the elastic forces on the bacteria
are more accurately represented in the Eulerian grid. 

\subsection{Variable Viscosity\label{sub:varVisc}}

\begin{figure}
\begin{centering}
(a)\includegraphics[width=3in]{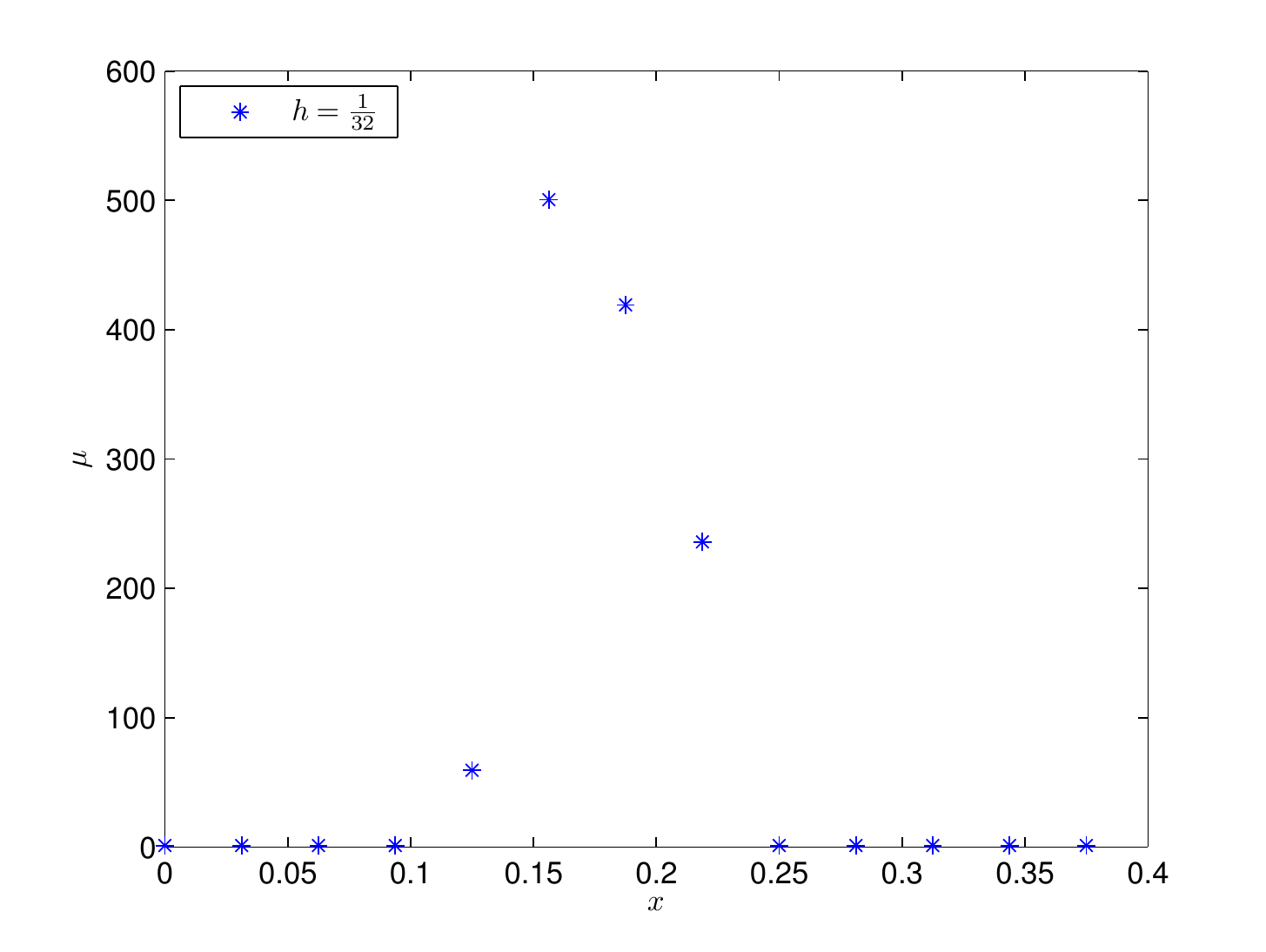}(b)\includegraphics[width=3in]{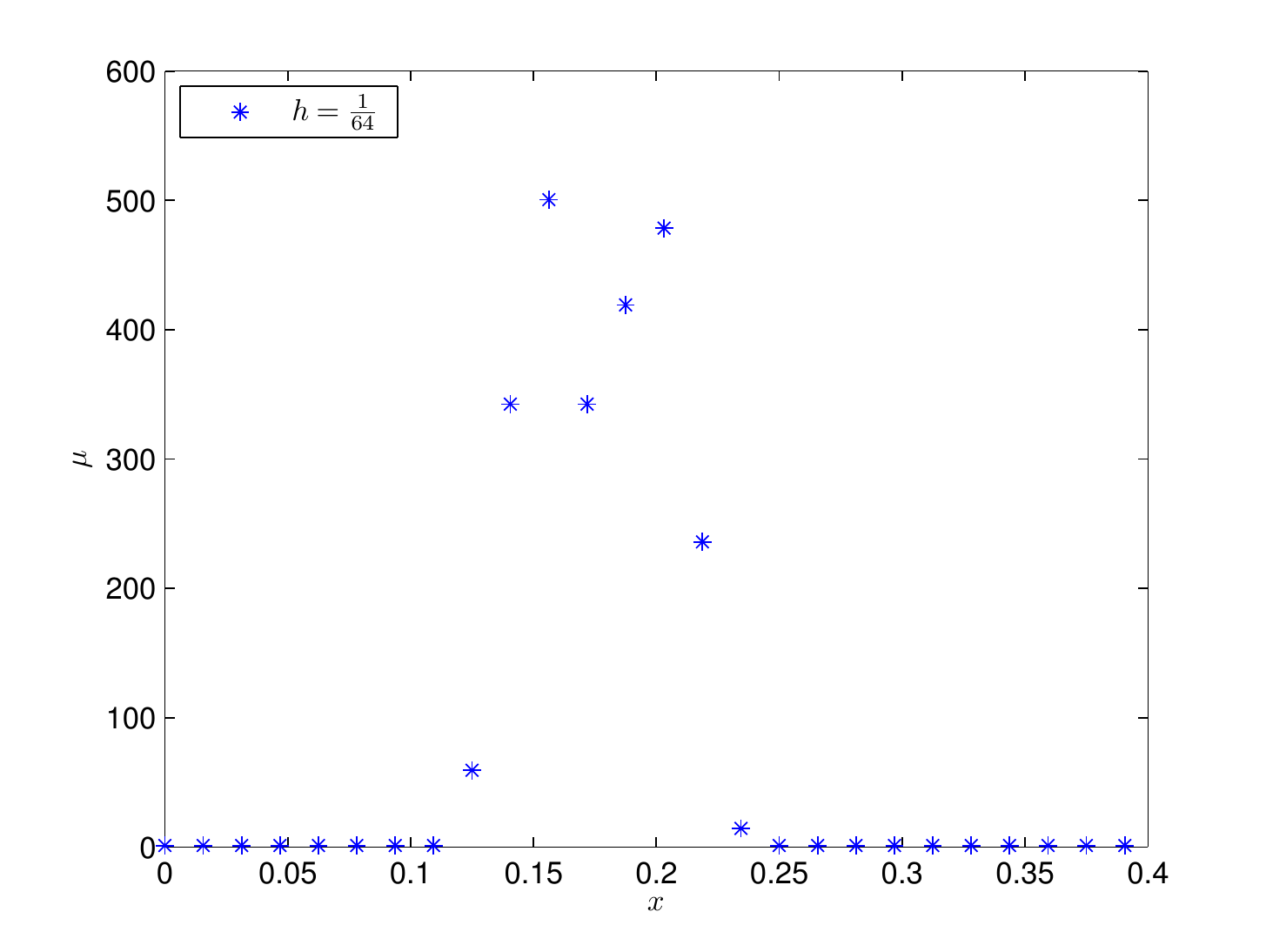}\\
(c)\includegraphics[width=3in]{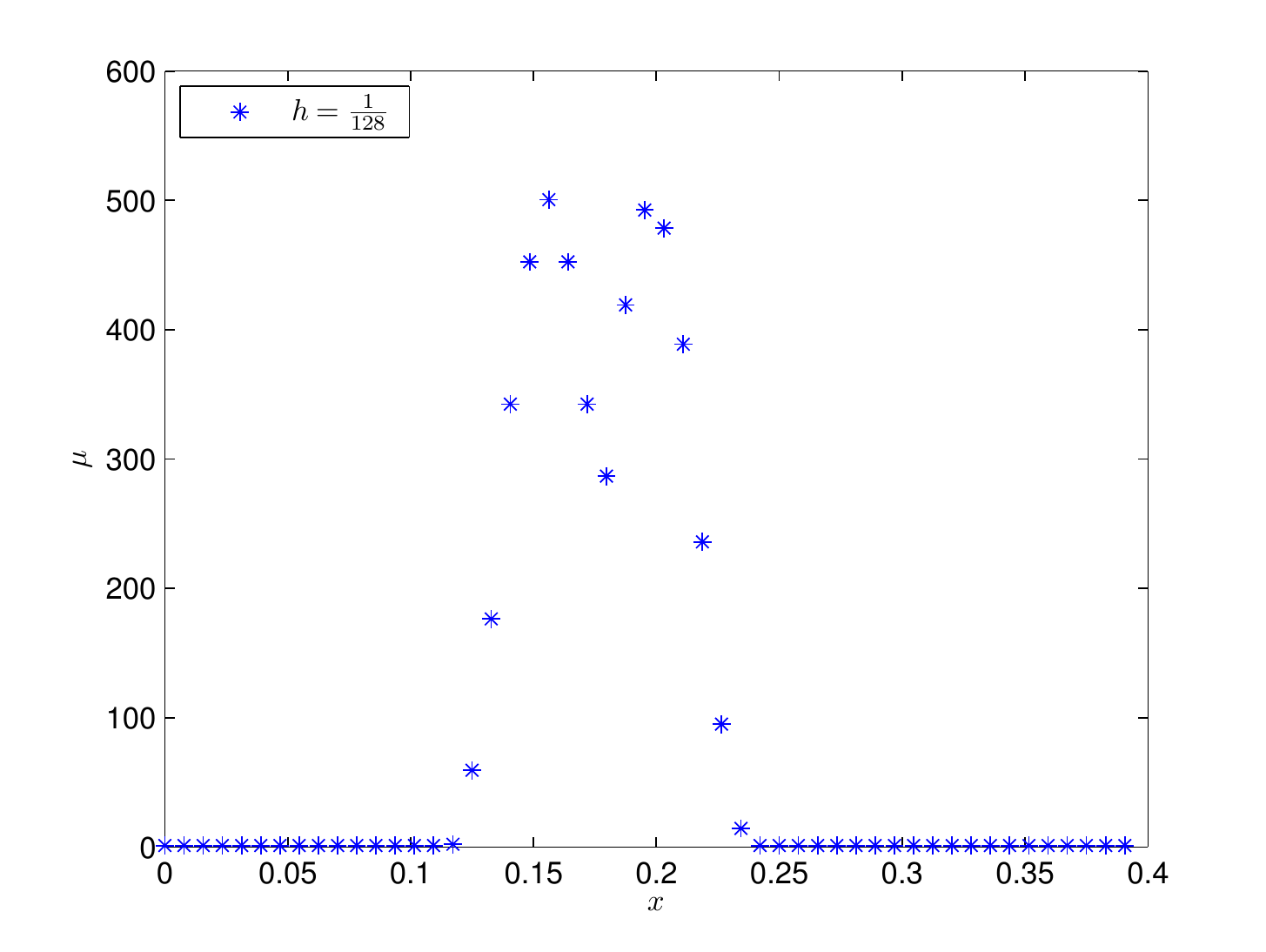}(d)\includegraphics[width=3in]{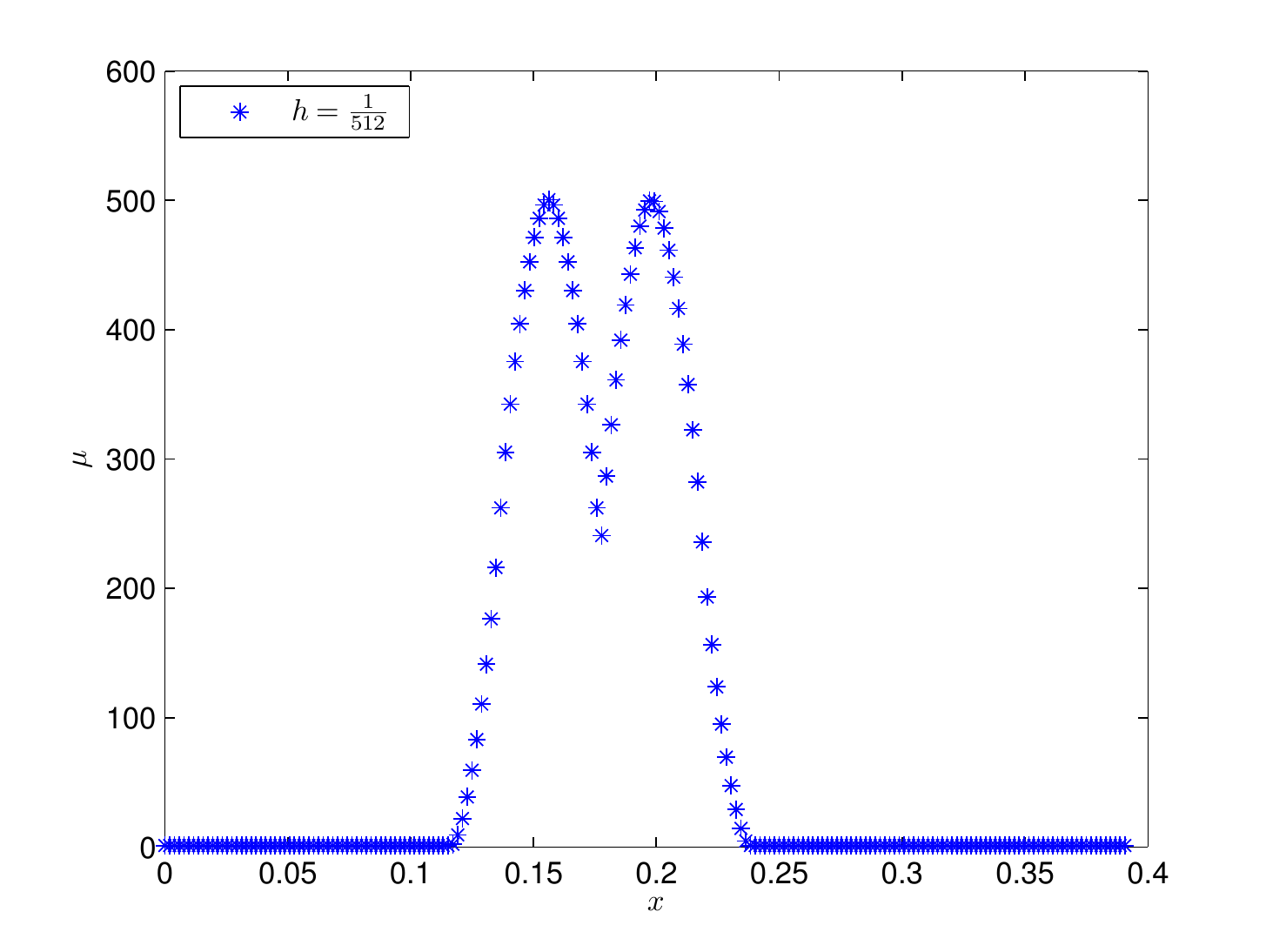}
\par\end{centering}

\caption{\label{fig:ViscEx}Examples of 1D $\mu(x)$ with two Lagrangian nodes,
one at $x=\nicefrac{5}{32}$ on an Eulerian node and the other at
$x=\frac{5}{32}+\frac{2\pi}{150}$, with $\omega=\frac{1}{50}$. These
Lagrangian nodes (cells) are close enough that there is a region of
interaction. These plots illustrate the effect of using different
spatial steps: (a) $h=\frac{1}{32}$, (b) $h=\frac{1}{64}$ (c) $h=\frac{1}{128}$
(d) $h=\frac{1}{512}$}
\end{figure}

It is known that ECM density decreases with distance from an individual
cell. To account for this, the exact form of $\mu(\mathbf{x})$ used
in our simulations is
\begin{equation}
\mu(\mathbf{x})=\max_{1\le s\le\eta,\, s\in\mathbb{N}}\left[\left(2\omega\right)^{D}\left(\mu_{max}-\mu_{out}\right)\tilde{\delta}(\mathbf{x}-\mathbf{X}(s),\,\omega)+\mu_{out}\right],\label{eq:ViscExp}
\end{equation}
where $\mu_{max}$ is the viscosity at a bacterial node, $\mu_{out}$
is the viscosity of the surrounding fluid, $D$ is the spatial dimension,
and $\omega$ is a parameter we can use to stretch the influence of
the additional viscosity. We made this choice for $\mu(\mathbf{x})$
because we wanted a viscosity that would decrease at the same rate
as the elastic force with the distance from the bacterial cell. See
\prettyref{fig:ViscEx} for a 1D example of the effect of $\omega$
and $h$ on the viscosity distribution, $\mu(\mathbf{x})$, from two
interacting cells. In the future, we can change this function to suit
the specific viscous properties of the particular biofilm. 

\subsection{Solution Strategy}

We employ a \textit{projection method} (\citep{BrownCortezMinion2001})
to solve the incompressible Navier-Stokes equations numerically, building
on the method used by Zhu and Peskin in \citep{Zhu2002}. This method
introduces a velocity field (at an intermediate time), $\mathbf{\tilde{u}}(\mathbf{x},\, t)$,
which is the solution to the difference equation
\begin{equation}
\rho^{n}\left(\sigma\frac{\tilde{\mathbf{u}}_{k}^{n+1}-\mathbf{u}_{k}^{n}}{\triangle t}+\frac{1}{2}(\mathbf{u}\cdot\mathbf{D}^{0}\mathbf{u}_{k}+\mathbf{D}^{0}\cdot(\mathbf{u}\mathbf{u}_{k}))^{n}\right)=\frac{1}{Re}\mathbf{D}^{0}\cdot\left(\mu^{n}\left(\mathbf{D}^{0}\tilde{\mathbf{u}}_{k}^{n+1}+D_{h,k}^{0}\tilde{\mathbf{u}}^{n+1}\right)\right)+\frac{Lf_{0}}{\rho_{0}u_{0}^{2}}\mathbf{f}_{k}^{n}\,,\label{eq:NSdis1-1}
\end{equation}
where $k=1,\,2,\,3$ and the subscripts $k$ denote the $k^{\textrm{th}}$
component of that vector. The finite difference operators in \prettyref{eq:NSdis1-1}
are originally defined by 
\begin{eqnarray*}
L_{h}(\phi(\mathbf{x})) & =\sum_{i=1}^{3} & \frac{\phi(\mathbf{x}+h\mathbf{e}_{i})+\phi(\mathbf{x}-h\mathbf{e}_{i})-2\phi(\mathbf{x})}{h^{2}}\,,
\end{eqnarray*}
\begin{eqnarray*}
D_{h,i}^{0}(\phi) & = & \frac{\phi(\mathbf{x}+h\mathbf{e}_{i})-\phi(\mathbf{x}-h\mathbf{e}_{i})}{2h}\,,\\
\mathbf{D}^{0} & = & (D_{h,1}^{0},D_{h,2}^{0},D_{h,3}^{0}),
\end{eqnarray*}
where $\mathbf{e}_{i}$ is the unit vector in the $i^{th}$ direction.
Additionally, $\mathbf{D}^{0}\cdot\left(a\mathbf{D}^{0}\phi\right)$
is defined for scalar functions, $a(\mathbf{x})$ and $\phi(\mathbf{x})$,
using the midpoint values of $a$ as 
\begin{equation}
\mathbf{D}^{0}\cdot\left(a\mathbf{D}^{0}\phi\right)=\sum_{i=1}^{3}\frac{a\left(\mathbf{x}+\frac{h}{2}\mathbf{e}_{i}\right)\frac{\phi(\mathbf{x}+h\mathbf{e}_{i})-\phi(\mathbf{x})}{h}-a\left(\mathbf{x}-\frac{h}{2}\mathbf{e}_{i}\right)\frac{\phi(\mathbf{x})-\phi(\mathbf{x}-h\mathbf{e}_{i})}{h}}{h}\,,\label{eq:divergenceGrad}
\end{equation}
 and $\mathbf{D}^{0}\cdot\left(aD_{h,k}^{0}\mathbf{u}\right)$ is
defined as
\begin{eqnarray}
\mathbf{D}^{0}\cdot\left(aD_{h,k}^{0}\mathbf{u}\right) & = & \frac{a\left(\mathbf{x}+\frac{h}{2}\mathbf{e}_{k}\right)\frac{\mathbf{u}_{k}(\mathbf{x}+h\mathbf{e}_{k})-\mathbf{u}_{k}(\mathbf{x})}{h}-a\left(\mathbf{x}-\frac{h}{2}\mathbf{e}_{k}\right)\frac{\mathbf{u}_{k}(\mathbf{x})-\mathbf{u}_{k}(\mathbf{x}-h\mathbf{e}_{k})}{h}}{h}+\nonumber \\
 &  & +\sum_{i\ne k}^{3}\left(\frac{a\left(\mathbf{x}+h\mathbf{e}_{i}\right)\frac{\mathbf{u}_{i}(\mathbf{x}+h\mathbf{e}_{i}+h\mathbf{e}_{k})-\mathbf{u}_{i}(\mathbf{x}+h\mathbf{e}_{i}-h\mathbf{e}_{k})}{2h}}{2h}\right.\label{eq:divergencePartial}\\
 &  & -\left.\frac{a\left(\mathbf{x}-h\mathbf{e}_{i}\right)\frac{\mathbf{u}_{i}(\mathbf{x}-h\mathbf{e}_{i}+h\mathbf{e}_{k})-\mathbf{u}_{i}(\mathbf{x}-h\mathbf{e}_{i}-h\mathbf{e}_{k})}{2h}}{2h}\right)\,.\nonumber 
\end{eqnarray}

Note that in the case of constant viscosity, using the incompressibility
constraint, \prettyref{eq:N-S1VarVisc} is reduced to the standard
N-S equation used in IBM,
\begin{eqnarray}
\rho(\mathbf{x},\, t)\left(\frac{\partial\mathbf{u}}{\partial t}+\mathbf{u}\cdot\triangledown\mathbf{u}\right) & = & -\triangledown p+\mu\triangle\mathbf{u}+\mathbf{f}(\mathbf{x},\, t).\label{eq:N-S1}
\end{eqnarray}
Therefore, in the case of constant viscosity, equation \prettyref{eq:NSdis1-1}
is replaced with
\begin{equation}
\rho^{n}\left(\sigma\frac{\tilde{\mathbf{u}}_{k}^{n+1}-\mathbf{u}_{k}^{n}}{\triangle t}+\frac{1}{2}(\mathbf{u}\cdot\mathbf{D}^{0}\mathbf{u}_{k}+\mathbf{D}^{0}\cdot(\mathbf{u}\mathbf{u}_{k}))^{n}\right)=\frac{1}{Re}L_{h}(\tilde{\mathbf{u}}_{k}^{n+1})+\frac{Lf_{0}}{\rho_{0}u_{0}^{2}}\mathbf{f}_{k}^{n}\,.\label{eq:NSdis1}
\end{equation}
We solve \prettyref{eq:NSdis1} as opposed to \prettyref{eq:NSdis1-1}
in the constant viscosity case, since each component of $\tilde{\mathbf{u}}$,
$\tilde{\mathbf{u}}_{k}$, is independent of each other in \prettyref{eq:NSdis1},
and can therefore be solved for separately and simultaneously.

To complete the discretized incompressible Navier-Stokes system, we
have the following two equations:
\begin{eqnarray}
\sigma\rho^{n}\left(\frac{\mathbf{u}^{n+1}-\mathbf{\tilde{u}}^{n+1}}{\triangle t}\right) & = & -\varepsilon\mathbf{D}^{0}p^{n+1}\,,\label{eq:NSdis2}\\
\mathbf{D}^{0}\cdot\mathbf{u}^{n+1} & = & 0\,.\label{eq:NSdis3}
\end{eqnarray}
We point out here that summing equations \prettyref{eq:NSdis1-1}
and \prettyref{eq:NSdis2} leads to the discretized version of \prettyref{eq:N-S1},
with the exception that the evaluation of the viscous term is at the
intermediate value of the velocity. We solve for pressure by applying
$\mathbf{D}^{0}$ to both sides of \prettyref{eq:NSdis2} and using
\prettyref{eq:NSdis3} to obtain
\begin{equation}
\mathbf{D}^{0}\cdot\left(\frac{1}{\rho^{n}}\mathbf{D}^{0}p^{n+1}\right)=\frac{\sigma}{\varepsilon}\frac{\mathbf{D}^{0}\cdot\tilde{\mathbf{u}}^{n+1}}{\Delta t}.\label{eq:NSdis4}
\end{equation}

To solve equations \prettyref{eq:NSdis1-1} and \prettyref{eq:NSdis4},
we use Gauss-Seidel as a smoother in a multigrid solver. At each time
step, we solve \prettyref{eq:NSdis1} for $\mathbf{\tilde{u}}^{n+1}$,
substitute it into \prettyref{eq:NSdis4}, solve for $p^{n+1}$, and
finally solve for $u^{n+1}$ using \prettyref{eq:NSdis2}. Then the
velocity is transferred from the Eulerian points to the Lagrangian
points using \prettyref{eq:discrVelCouple}. With $\mathbf{U}^{n+1}$
computed, the new Lagrangian node locations are computed using Euler's
method as
\[
\mathbf{X}^{n+1}(s)=\frac{\Delta t}{\sigma}\mathbf{U}^{n+1}(s)+\mathbf{X}^{n}(s).
\]
The forces between the Lagrangian points are then recalculated and
transferred to the Eulerian points using \prettyref{eq:discrForceCouple}.
Finally, the values of $\rho^{n}$ and $\mu^{n}$ are evaluated using
the new Lagrangian locations.

\subsection{Multigrid}

In this section, we discuss the elements of multigrid that we use
in our solution strategy. For more details on multigrid, see \citep{Briggs2000}. 

In our solver, we use the conventional Gauss-Seidel iterative method
with red-black ordering. In the multigrid scheme, we use full-weighting
restriction to go from fine to coarse grids, and we use linear interpolation
to go from coarse back to fine grids. The finest grid is the grid
with step size $h$ and the grids become coarser by increasing the
step size by a factor of 2. This halves the number of nodes in each
dimension, allowing for significantly faster computations on the coarser
grids. The number of levels in the multigrid solver depends on both
the dimensions of the computational domain as well as $h$. In our
simulations, we iterate using multigrid V-cycles until we reach a
sufficiently low value for the norm of the residual, 
\[
\left\Vert f^{h}-A^{h}\tilde{v}^{h}\right\Vert ,
\]
at each time step. Here, $A^{h}v^{h}=f^{h}$ is the linear discretization
of a PDE, and the residual is $r^{h}=f^{h}-A^{h}\tilde{v}^{h}$, where
$\tilde{v}$ is an approximation to $v$. The residual provides a
bound on the true error in the solution of the linear system since
we have this relationship between the error and the relative residual
error:
\begin{equation}
\frac{\left\Vert e^{h}\right\Vert }{\left\Vert v^{h}\right\Vert }\le cond(A^{h})\frac{\left\Vert r^{h}\right\Vert }{\left\Vert f^{h}\right\Vert }\,,\label{eq:resCondrelationship}
\end{equation}
where $e^{h}=v^{h}-\tilde{v}^{h}$, $r^{h}=f^{h}-A^{h}\tilde{v}^{h}$,
and $cond(A^{h})$ is the condition number of $A^{h}$ \citep{Atkinson1989}.
In \prettyref{sec:Validation}, we give approximations for the condition
numbers for our matrices.

\subsection{Boundary Conditions}

\begin{figure}
\begin{centering}
(a)\includegraphics[width=2in]{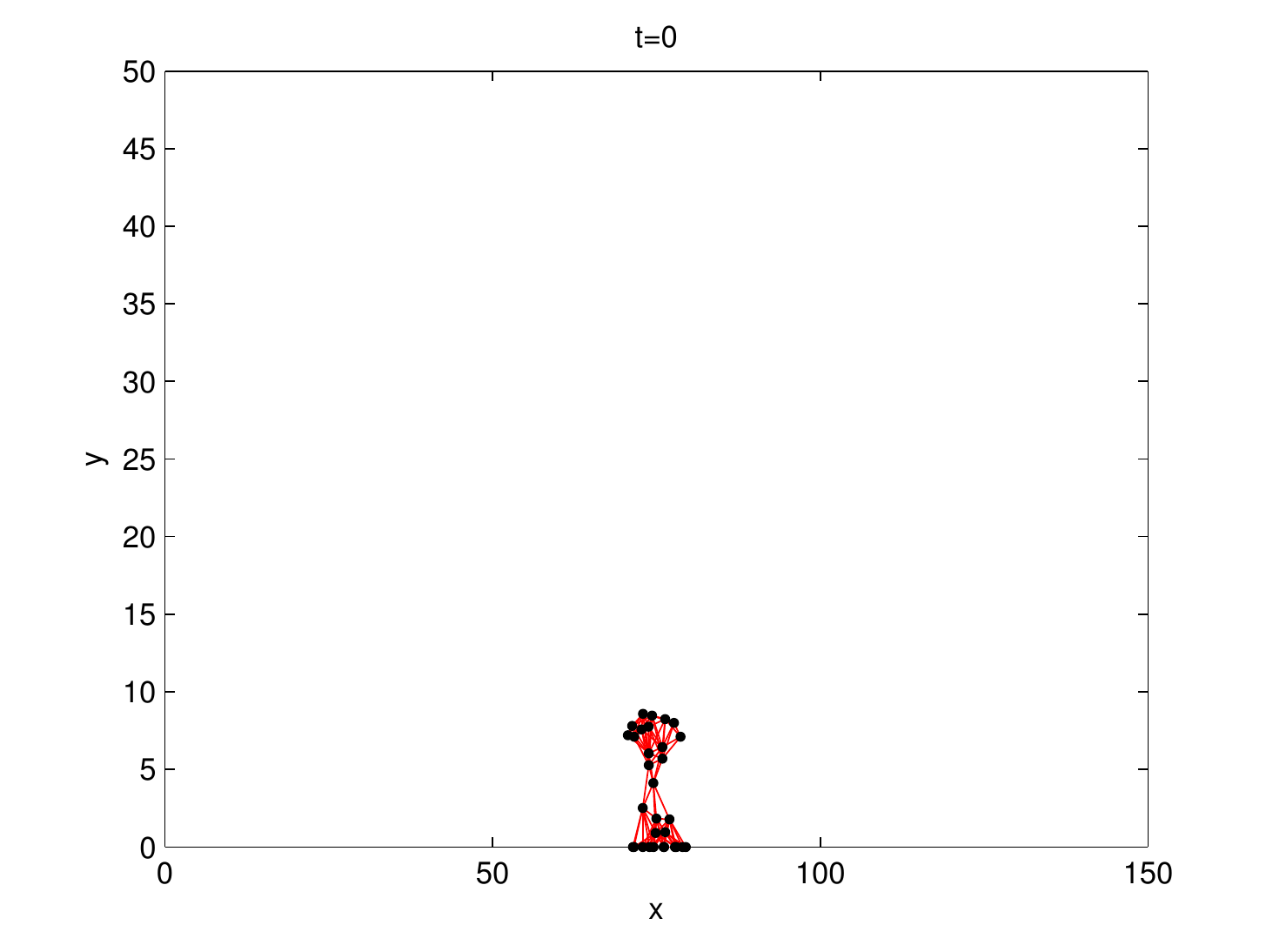}(b)\includegraphics[width=2in]{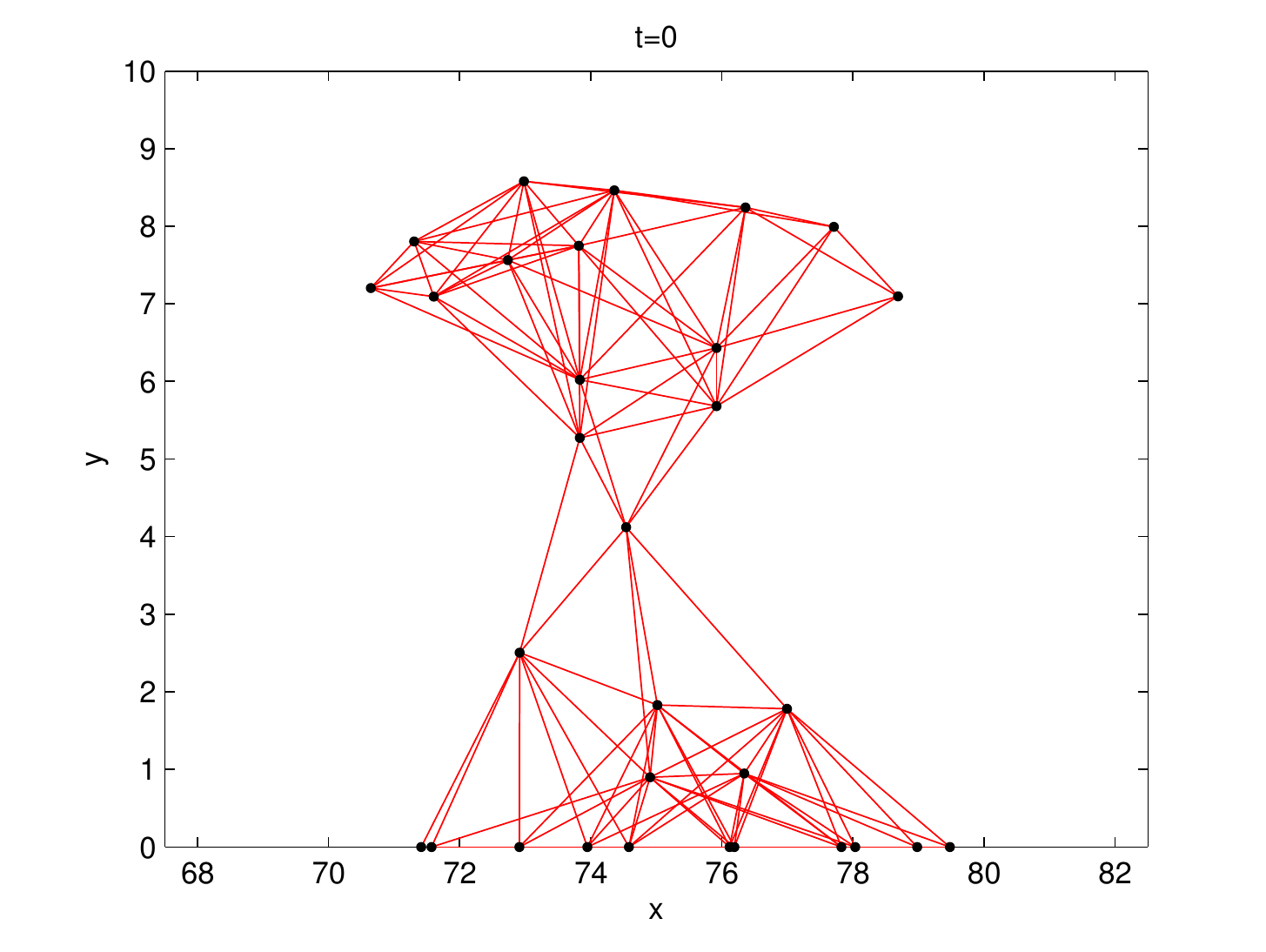}
\par\end{centering}

\begin{centering}
(c)\includegraphics[width=2in]{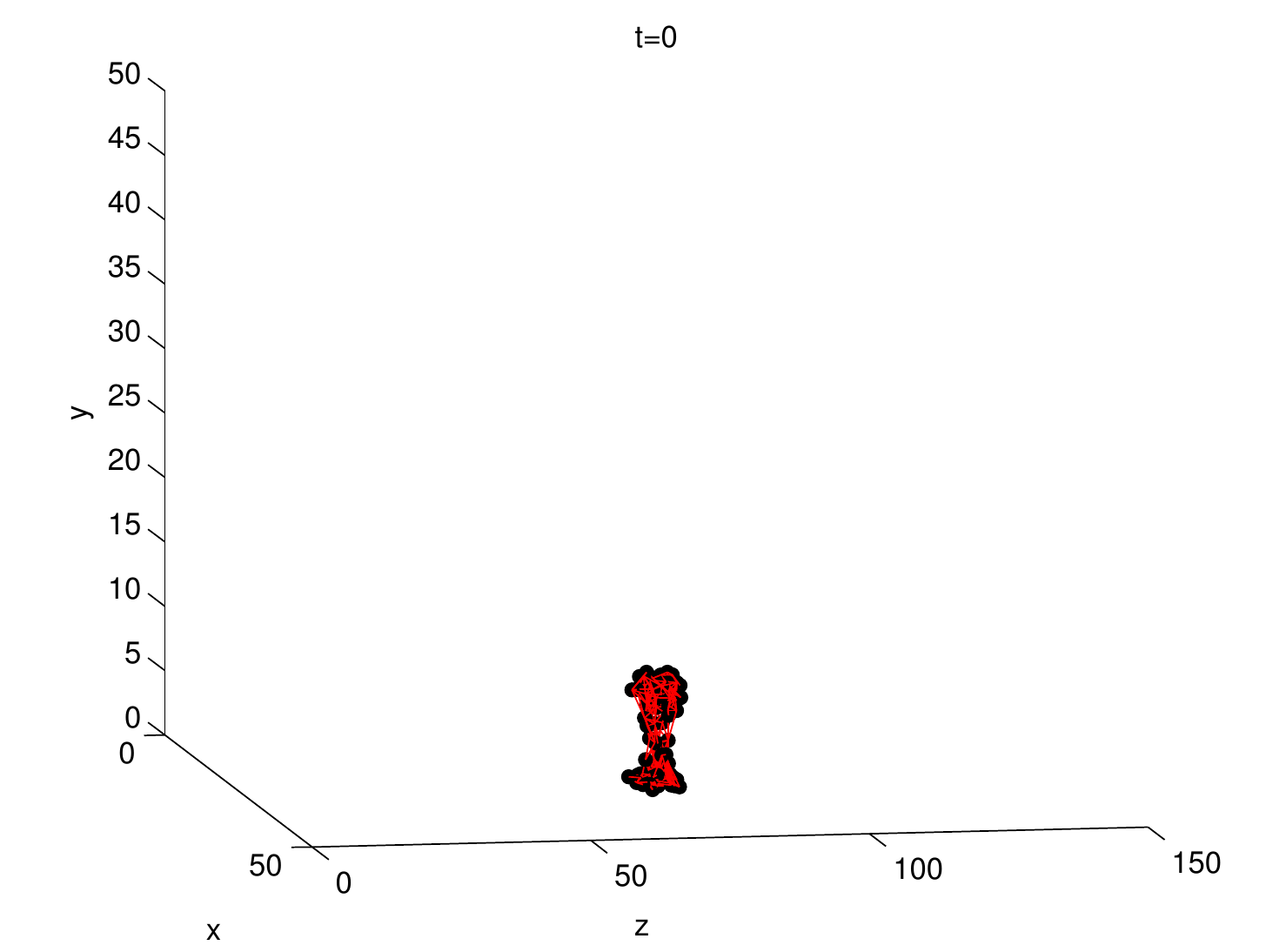}(d)\includegraphics[width=2in]{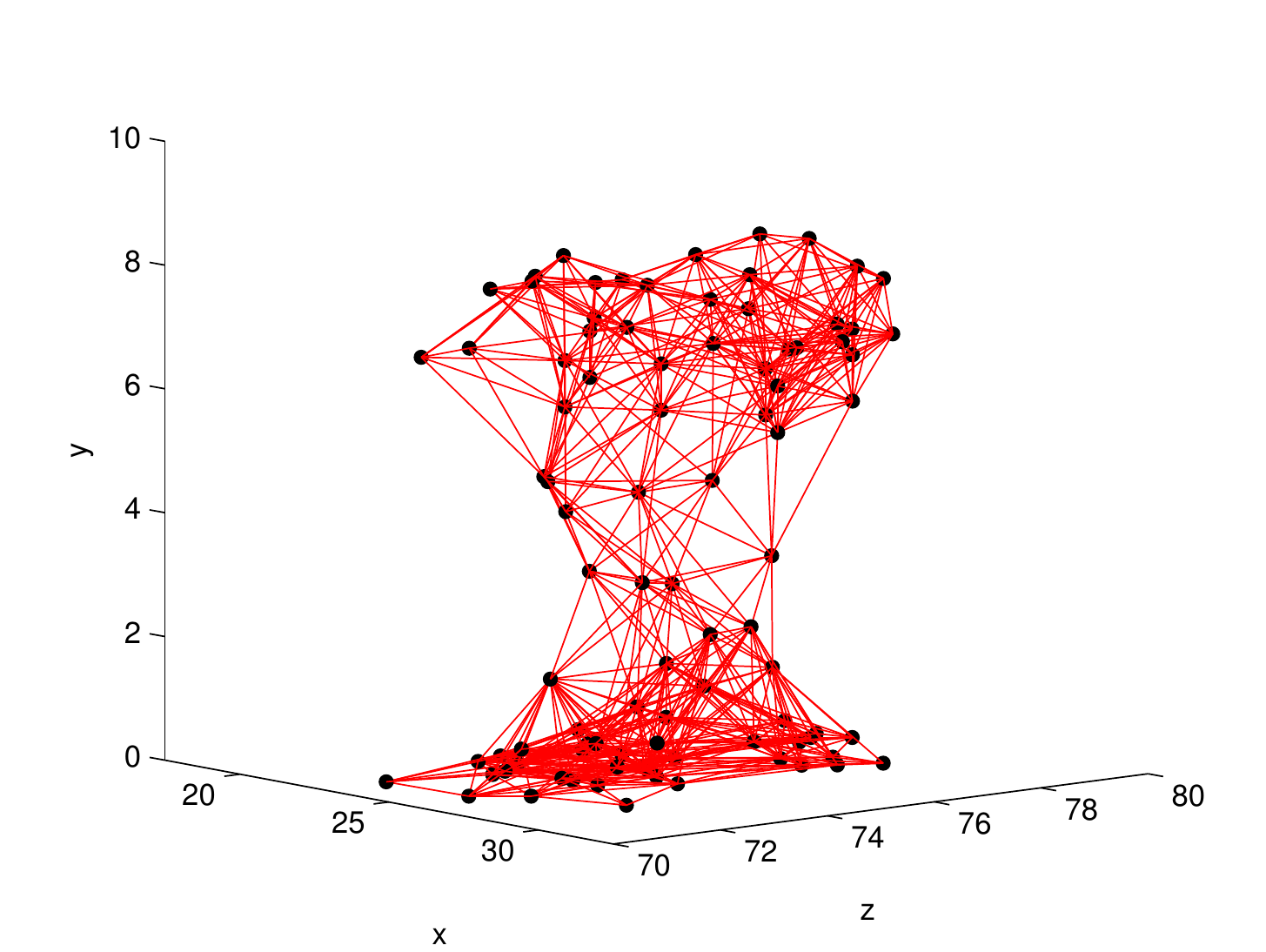}
\par\end{centering}

\caption{\label{fig:Comp Domain Examples}Shows examples of the computational
domains with a sample biofilm. Axes units are microns. (a) is 2 dimensional
and (c) is 3 dimensional. (b) and (d) show enlarged images of the
biofilms. The points shown in these plots show the initial position
of all of the bacterial cells in the biofilm.}
\end{figure}

The computational domain used in our example simulations is a section
of a tube with the biofilm centered in the direction along the axis
of the tube (see \prettyref{fig:Comp Domain Examples}). In the 2D
case, flow is along the $x$-axis and, in 3D, it is along the $z$-axis.
The boundary conditions we used in these simulations were derived
from exact solutions for the velocity and pressure in the case of
laminar flow. We now provide the boundary conditions in both the 2D
and 3D cases.

\subsubsection{2D Boundaries\label{sub:2D-Boundaries}}

The no-slip boundary condition exists at the walls of tube and requires
that the velocity be zero there, so we use that as the boundary condition
at the walls. The velocity at the upstream boundary ($x=0$) is held
at the laminar flow velocity (shown in \prettyref{fig:LaminarFlowLeftBound}(a))
given by 
\begin{equation}
u_{1}(y)=\frac{\kappa}{2\mu}(y^{2}-2ay),\label{eq:2DlaminarSoln}
\end{equation}
where $a$ is the radius of the tube, $y$ is the displacement from
the bottom edge of the tube, $\kappa$ is the linear rate at which
the pressure decreases through the tube, and $u_{1}$ is the x-component
of the velocity (i.e., $\mathbf{u}=\left(u_{1},\, u_{2}\right)$).
At the downstream boundary, a Neumann condition is applied to the
velocity by enforcing that 
\[
\left(\frac{\partial}{\partial x}\mathbf{u}(x,y)\right)_{x=x_{down}}=0\qquad\forall y,
\]
where $x_{down}$ represents the $x$ value at the downstream boundary.

The boundary conditions for pressure come from the laminar flow equation
for pressure given by
\begin{equation}
p(x)=\kappa x+p(0).\label{eq:pressLam}
\end{equation}
In our simulations, we hold the pressure at the upstream boundary
at $p(0)$ and at $p(x_{down})$ at the downstream. At the top boundary,
we hold the pressure at the values given by \prettyref{eq:pressLam}
and, at the bottom boundary (the boundary on which the biofilm is
attached), we use a Neumann boundary
\[
\left(\frac{\partial}{\partial y}p(x,y)\right)_{y=0}=0,\qquad\forall x.
\]

\subsubsection{3D Boundaries}

\begin{figure}[h]
\begin{centering}
(a)\includegraphics[width=2in]{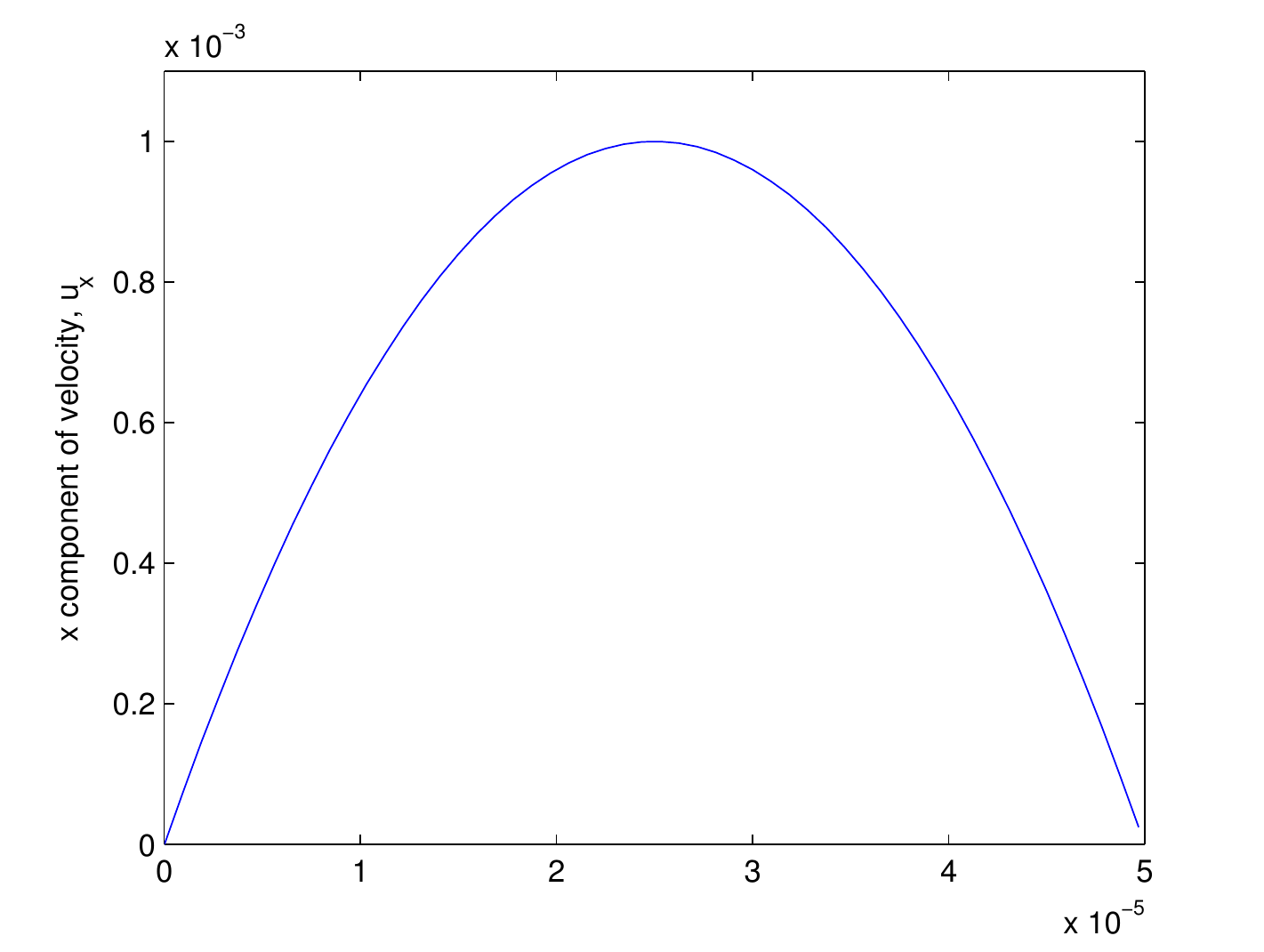}(b)\includegraphics[width=2in]{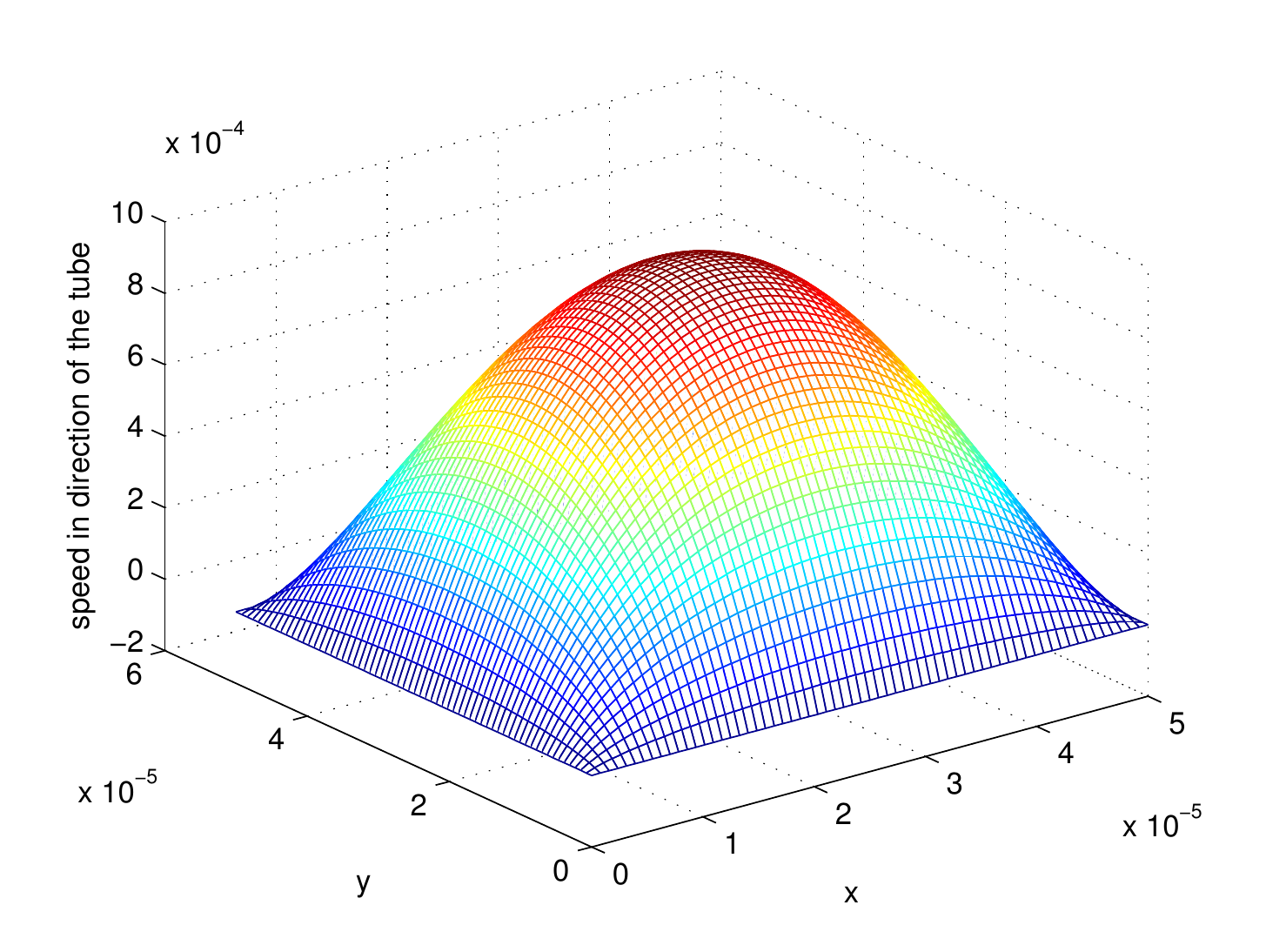}
\par\end{centering}

\caption{\label{fig:LaminarFlowLeftBound}We illustrate upstream boundary conditions:
(a) 2D laminar flow velocity profile, (b) 3D laminar flow velocity
profile. Units are in meters for (a) and (b).}
\end{figure}

In the 3D simulations, we orient the square tube along the $z$-axis
(see \prettyref{fig:Comp Domain Examples}(c)). The no-slip boundary
condition exists at the walls of tube and requires that the velocity
be zero there so we use that as the boundary condition at the walls.
Derived by Spiga in \citep{Spiga1994}, the velocity at the upstream
boundary is held at the laminar flow velocity (shown in \prettyref{fig:LaminarFlowLeftBound}(b))
given by 
\begin{equation}
u_{3}(x,y)=-\frac{16\kappa a^{2}}{\mu\pi^{4}}\sum_{n,m>0,\,\textrm{odd}}\frac{\sin\left(n\pi\nicefrac{x}{a}\right)\sin\left(m\pi\nicefrac{y}{a}\right)}{nm\left(n^{2}+m^{2}\right)},\label{eq:3DlaminarFlow}
\end{equation}
where $a$ is the width of the tube and $u_{3}$ is the z-component
of the velocity (i.e., $\mathbf{u}=\left(u_{1},\, u_{2},\, u_{3}\right)$).
At the downstream boundary, a Neumann condition is applied to the
velocity by enforcing that 
\[
\left(\frac{\partial}{\partial z}\mathbf{u}(x,y,z)\right)_{z=z_{down}}=0\qquad\forall x,\, y.
\]

The boundary conditions for pressure come from the laminar flow equation
for pressure given by
\begin{equation}
p(z)=\kappa z+p(0),\label{eq:pressLam3D}
\end{equation}
where $z=0$ is the upstream boundary. In our simulations, we hold
the pressure at the upstream boundary at $p(0)$ and at $p(z_{down})$
at the downstream. At the top and side boundaries, we hold the pressure
at the values given by Equation \prettyref{eq:pressLam3D} and, at
the bottom boundary (side with attached biofilm), we use the Neumann
condition given by
\[
\left(\frac{\partial}{\partial y}p(x,y,z)\right)_{y=0}=0\qquad\forall x,\, z.
\]

\section{\label{sec:Validation}Validation}

It is crucial to validate our numerical method using known results
and mesh refinement convergence analysis. Thus, below we provide numerical
evidence that both our 2D and 3D simulations and numerical methods
are working as they should.

For the purposes of the 2D and 3D convergence analysis conducted in
later sections, we require the following notation. We present the
following notation in 3D (the 2D versions are analogous but without
the $z$ elements). Define the \textit{Eulerian grid function $p$-norm}
for an arbitrary 3D vector field, $\mathbf{w}(\mathbf{x})=(w_{1}(\mathbf{x}),w_{2}(\mathbf{x}),w_{3}(\mathbf{x}))$,
by
\begin{equation}
\left\Vert \mathbf{w}\right\Vert _{p}=\left(\sum_{i,j,k}\left|\mathbf{w}(x_{i},y_{j},z_{k})\right|{}^{p}h^{D}\right)^{\nicefrac{1}{p}},\label{eq:EulPnorm}
\end{equation}
where $D$ is the spatial dimension, $1\le p<\infty$, and 
\[
\left|\mathbf{w}(x_{i},y_{j},z_{k})\right|=\sqrt{w_{1}(x_{i},y_{j},z_{k})^{2}+w_{2}(x_{i},y_{j},z_{k})^{2}+w_{3}(x_{i},y_{j},z_{k})^{2}}.
\]
Then 
\[
\left\Vert \mathbf{w}\right\Vert _{\infty}=\max_{i,j,k}\left|\mathbf{w}(x_{i},y_{j},z_{k})\right|.
\]
Additionally, on the Lagrangian grid define the \textit{Lagrangian
grid function $p$-norm} for a vector field, $\mathbf{X}=(X_{1}(s),X_{2}(s),X_{3}(s))$,
as 
\[
\left\Vert \mathbf{X}\right\Vert _{p}=\left(\sum_{s=1}^{\eta}\left|\left(X_{1}(s),\, X_{2}(s),\, X_{3}(s)\right)\right|^{p}d_{0}^{D}\right)^{\nicefrac{1}{p}},
\]
where $1\le p<\infty$ and $d_{0}^{D}$ is the average volume element
of the Lagrangian nodes. Then 
\[
\left\Vert \mathbf{X}\right\Vert _{\infty}=\max_{1\le s\le\eta}\left|\left(X_{1}(s),\, X_{2}(s),\, X_{3}(s)\right)\right|.
\]
Note that both of these grid function norms are derived from using
discretizations of the integrals used in a typical function p-norm
(see Appendix A of \citep{LeVeque2007} for more details). 

There are three parts to our simulation validation process: 1) we
illustrate that in the absence of the biofilm our numerical simulation
converges to the analytical solution; 2) we verify that our multigrid
technique is correctly accelerating the convergence of our chosen
relaxation scheme; and 3) we determine the convergence rate of the
simulations with a biofilm using a mesh refinement convergence analysis. 

Before discussing the results of our validation process, we provide
a brief description of the two primary sources of error present in
our simulations, \prettyref{sub:Discussion-of-Errors}. We also setup
our simulations with a detailed description of initial Lagrangian
node positions in \prettyref{sub:SimSetup}. Then we provide 2D simulation
validation in \prettyref{sub:Two-Dimensional-Validation} and 3D validation
in \prettyref{sub:Three-Dimensional-Validation}.

\subsection{\label{sub:Discussion-of-Errors}Discussion of Errors}

In our numerical scheme, we have two sources of error: 1) \textit{discretization
error} is introduced by discretizing the Navier-Stokes equations in
space and time; and 2) \textit{algebraic error} is introduced when
we attempt to solve the resultant systems of linearized equations. 

As it is impossible to compute the true algebraic error, we use the
norm of the residual to deduce an upper bound on the algebraic error
using Equation \prettyref{eq:resCondrelationship}. Recall that Equation
\prettyref{eq:resCondrelationship} indicates that the relative algebraic
error at each timestep is no larger than the condition number of the
matrix times the relative residual norm (recall that $\frac{\left\Vert r\right\Vert }{\left\Vert f\right\Vert }$
is the relative residual norm). We do not construct these matrices
during the actual simulations because we do not need them to solve
the systems. However, we did construct them to find their condition
numbers and found that the condition numbers for the matrices used
in the computations for $\tilde{\mathbf{u}}^{n+1}$ and $p^{n+1}$
are $O\left(h^{-2}\right)$ (this is true for both 2D and 3D simulations).
The simulations that resulted in the plots given in \prettyref{sub:Two-Dimensional-Simulations}
and \prettyref{sub:Three-Dimensional-Simulations} were run using
$h=\frac{1}{128}$, and thus the matrix condition numbers were approximately
$10^{4}$.

Our goal in the simulations is to ensure that the algebraic error
falls well below the discretization error at each time step, so the
total error will be dominated by the discretization error. In theory,
the discretization error is at best $O(h^{2})\approx C\left(\nicefrac{1}{128}\right)^{2}$$\approx C\times6\times10^{-5}$,
for $C>0$, with our discretization. Using a stopping criteria of
$10^{-9}$ for the relative residual at each timestep should suffice
(i.e., from \prettyref{eq:resCondrelationship}). We continue to the
next time step only when the computed relative residual, $\frac{\left\Vert r^{h}\right\Vert }{\left\Vert f^{h}\right\Vert }\le10^{-9}$,
because this implies $\left\Vert e^{h}\right\Vert \le cond(A)*10^{-9}\approx10^{-5}$. 

Another factor influencing the capability of our simulations is that
after, extensive simulation, we discover that our linear solver is
limited to converging to a relative residual norm of about $10^{-11}$
(possibly from machine precision issues). With $h=\frac{1}{512}$,
the condition number is $O(10^{5})$ and the discretization error
is $O(10^{-6})$, so the algebraic error is at best bounded by about
$10^{-11}\times10^{5}=10^{-6}$ (see Equation \prettyref{eq:resCondrelationship}),
and we can no longer be certain that the algebraic error falls below
the discretization error at each timestep. For this reason, we restrict
$h$ to be larger than $\frac{1}{512}$ in all of the simulations
and convergence analysis.

\subsection{\label{sub:SimSetup}Simulation Setup}

\begin{figure}
\begin{centering}
(a)\includegraphics[width=3in]{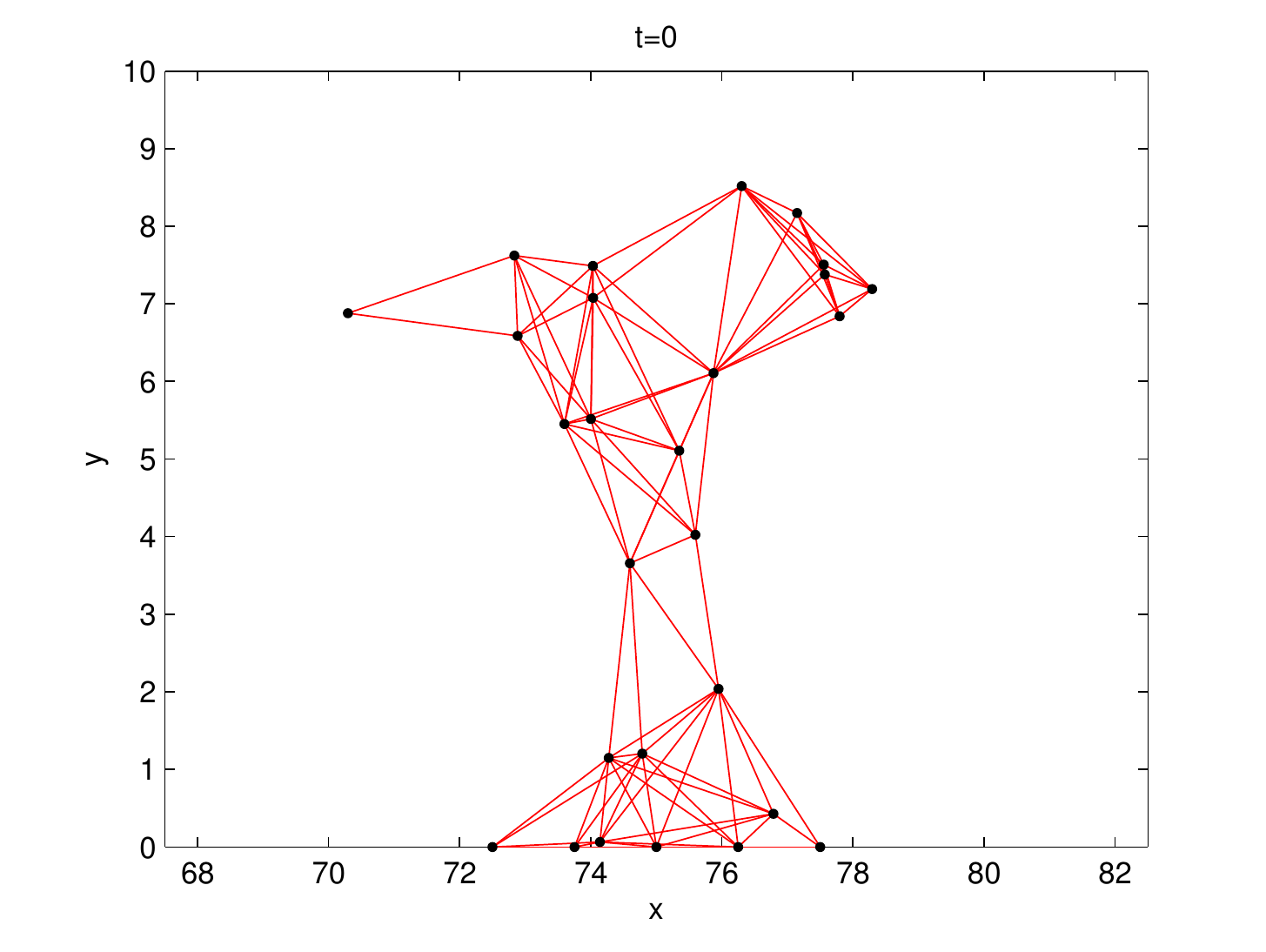}(b)\includegraphics[width=3in]{sample3Dbiofilm1}
\par\end{centering}

\caption{\label{fig:Mushroom-shaped-biofilm2Dsims}Mushroom-shaped biofilm
at $t=0$ in the middle of the computational domain attached to the
bottom $(y=0)$ of the tube. (a) This is the shape used in the 2D
simulations, and (b) 3D mushroom shaped biofilm at $t=0$.}
\end{figure}

In these convergence simulations, we constructed an experimentally
motivated mushroom shaped biofilm (shown in \prettyref{fig:Mushroom-shaped-biofilm2Dsims}(a)).
We carved this shape from a $1.6\,\mu\textrm{m}$ slice cut from data
points generated in the Younger and Solomon labs at the University
of Michigan. These data points are 3D bacterial cell locations from
3D Leica SP2 confocal laser scanning microscopy images taken of \textit{Staphylococcus
epidermidis} RP62A (ATCC 35984) grown in a Stovall 3 channel flow
cell for 24 hours at 37 C under a wall shear stress of $0.01\,\textrm{Pa}$.
For further details of how the coordinates were computed, see Stewart
et al.~\citep{Stewart}.

From this data, the average Lagrangian volume element, $d_{0}^{3}$,
is calculated to be approximately $4.036\,\mu\textrm{m}^{3}$, and
thus we use $d_{0}=1.59\,\mu\textrm{m}$ in both the 2D and 3D simulations.
We connect the initial distribution of cells with a distance based
connection criteria. Our inspiration for the connection distance criteria
came from the closeup images of biofilms such as the one shown in
\prettyref{fig:BiofilmPicture}. We observed that each bacterial cell
is connected to neighboring cells that are within about $2d_{0}$.
Thus, we varied the connection criteria in our algorithm between $1.5$-$2.5\times d_{0}$
in an effort to find one that resulted in a biofilm that was sufficiently
connected but not overcrowded. This resulted in the choice of a connection
criteria of $d_{c}=2.8\,\mu\textrm{m}$. In other words, we placed
spring connections between Lagrangian nodes at the beginning of the
simulation with every node connected to every other node less than
$2.8\,\mu\textrm{m}$ away. Admittedly, this value of $d_{c}$ is
arbitrary, and future work will include deriving a method to determine
this connection criteria through image analysis of closeup images
of biofilms similar to \prettyref{fig:BiofilmPicture}. The mushroom
shaped biofilm has a height of about $8.5\,\mu\textrm{m}$ and width
of about $8\,\mu\textrm{m}$ (see \prettyref{fig:Mushroom-shaped-biofilm2Dsims}).
In the convergence simulations, the maximum spring force, $F_{max}$,
is set to $5\times10^{-6}\,\textrm{N}$. The fluid parameters for
these convergence simulations are provided later in \prettyref{tab: 2D sim params}.

Note that $\tilde{\delta}$ is a function of $\omega$, a scaling
parameter we must choose that determines the volume/area of influence
when the forces and density are transferred to the Eulerian grid.
We must also point out here that a more accurate representation of
the Dirac delta function occurs when $\omega\ge h$. Thus, for the
purpose of the convergence simulations, we use $\omega=1.0\,\mu\textrm{m}$
in the transfer equations, \prettyref{eq:discrForceCouple} and \prettyref{eq:discrDensCouple}.
However, in our simulations, we use $\omega=0.5\,\mu\textrm{m}$ since
the actual radius of \textit{Staphylococcus epidermidis} is known
to be about $0.5\,\mu\textrm{m}$ \citep{Todar2012}. Using a characteristic
length of $L=50\,\mu\textrm{m}$, we have the non-dimensionalized
$\omega^{*}=\frac{\omega}{L}=\frac{1}{50}$. Dropping the star from
the dimensionless variable, we use $\omega=\frac{1}{50}$ in the convergence
simulations and $\omega=\frac{1}{100}$ in the results simulations.
We desire that $\omega\ge h$, so that the Lagrangian forces are spread
at least two Eulerian mesh widths in every direction (as is done in
the traditional IBM \citep{Peskin2002}). Using $\omega=\frac{1}{50}$
in the convergence simulations allows $h=\frac{1}{64},\,\frac{1}{128},\,\frac{1}{256}$
to obey these criteria. We again note that one of the reasons for
using $\omega$ in the scaling of \prettyref{eq:Dirac Delta Approx}
as opposed to $h$ is that better spatial convergence rates are achieved
since the scaling is independent of $h$.

\subsection{\label{sub:Two-Dimensional-Validation}Two-Dimensional Validation}

\begin{figure}
\begin{centering}
\includegraphics[width=4in]{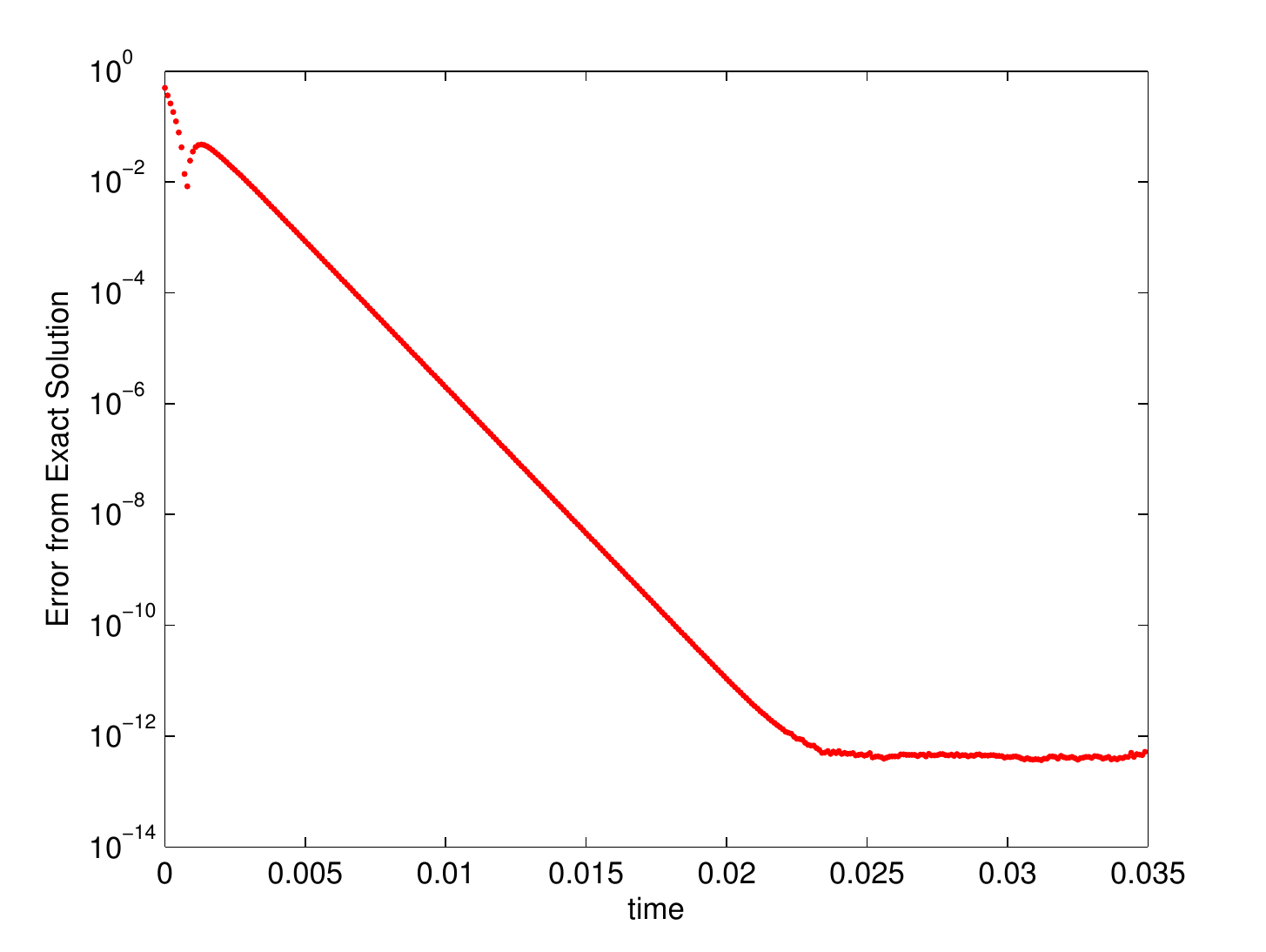}
\par\end{centering}

\caption{\label{fig:2DLaminarConvergence}Exact error in biofilm-free 2D simulation
with $h=\frac{1}{128}$. The solution is compared to the exact solution,
\prettyref{eq:2DlaminarSoln}, using a maximum norm. }
\end{figure}

In the absence of a biofilm, we expect that using a centered finite
difference approximation for the second derivatives allows \textit{exact}
convergence to the second-order polynomial solution (Equation \prettyref{eq:2DlaminarSoln}).
That is to say, we expect the numerical solution to converge to the
analytical within machine precision. Reassuringly, we find that the
biofilm-free simulations converge exactly to the steady state laminar
flow.%
\footnote{In the 3D simulations, we do not see exact convergence since the laminar
solution is not a second-order polynomial. See \prettyref{sub:Three-Dimensional-Simulations}
for details on the convergence properties of the 3D laminar flow case. %
}\negthinspace{}\negthinspace{} To illustrate, we started with an
initial velocity profile that is one-half of that of the laminar flow
velocity profile. The error in the simulation converged (within machine
precision) in less than 300 time iterations for all spatial resolutions
(see \prettyref{fig:2DLaminarConvergence} for example with $h=\nicefrac{1}{128}$
and $dt=0.0001$). 

\subsubsection{\label{sub:Multigrid-performance}Multigrid performance}

\begin{figure}[h]
\begin{centering}
(a)\includegraphics[width=3in]{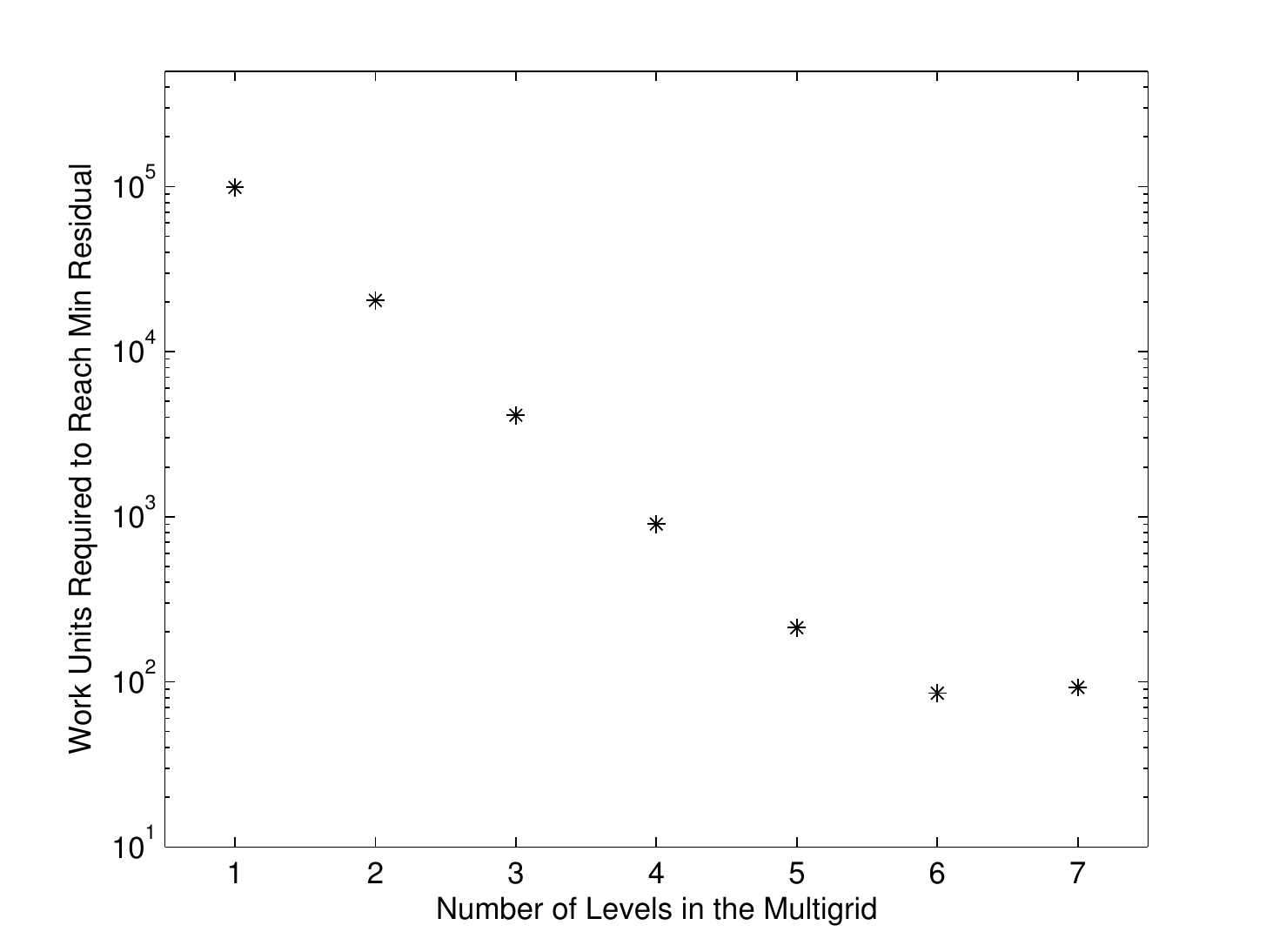}(b)\includegraphics[width=3in]{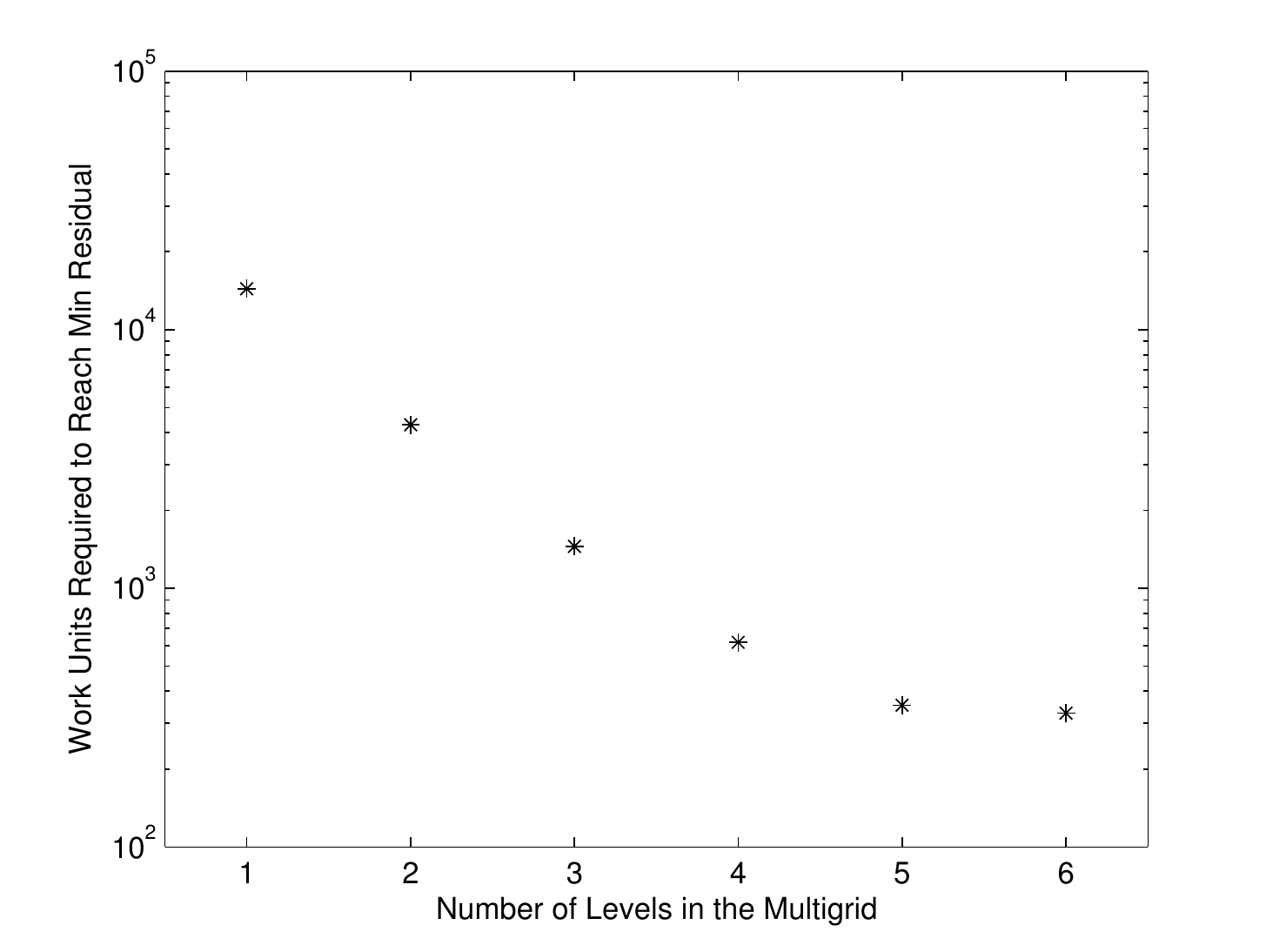}
\par\end{centering}

\caption{\label{fig:Work Units Decrease 2D}Decrease in work units required
to reach the minimum residual error as the number of multigrid levels
is increased in the (a) 2D simulations, (b) and 3D simulations. }
\end{figure}

Next, we provide numerical evidence that the multigrid technique convergences
optimally to the solutions of \prettyref{eq:NSdis1} and \prettyref{eq:NSdis4}.
Define a \textit{work unit} as the cost of performing one relaxation
on the finest grid (see \citep{Briggs2000}). In \prettyref{fig:Work Units Decrease 2D}(a),
we depict (for the pressure computation) the work units required to
reach the minimum residual error as a function of allowed levels in
the multigrid. This result shows that the number of work units required
decreases significantly with each added multigrid level. This means
that the multigrid method correctly accelerates the convergence of
our iterative method by doing computations on the coarser grids. For
example, with just one allowed level of multigrid, the relaxation
uses only the finest resolution grid and requires about $10^{5}$
work units, whereas with 6 multigrid levels we only require about
$10^{2}$ work units to achieve the same error. Note that there is
no reduction in the number of required work units with the addition
of a $7^{\textrm{th}}$ level in the multigrid, so we use at most
6 levels in our 2D solvers. The data in this plot was obtained using
our 2D simulation with a mushroom shaped biofilm similar to those
shown in \prettyref{sub:Two-Dimensional-Simulations}.

\subsubsection{Empirical Estimate of Convergence Rate in Time}

Similar to the development by Mori and Peskin in \citep{Mori2008},
we define a measure of error by 
\begin{equation}
E_{p}(q(T);\Delta t)=\left\Vert q^{\Delta t}(T)-q^{\Delta t/2}(T)\right\Vert _{p},\label{eq:ErrorTemporal}
\end{equation}
which is the error difference at time $t=T$ in a computed quantity,
$q$, using a temporal refinement of a half timestep. Then, an \textit{empirical
estimate} for the convergence rate is calculated using 
\begin{equation}
r_{p}(q(T);\Delta t)=\log_{2}\left(\frac{E_{p}(q(T);\Delta t)}{E_{p}(q(T);\frac{\Delta t}{2})}\right).\label{eq:temporal convergence rate}
\end{equation}
 We compute the approximate convergence rate in time using the $E_{2}$
and $E_{\infty}$ errors in the Eulerian variable, $\mathbf{u}$,
and in the Lagrangian variable, $\mathbf{X}$. We simulate until $t=T=0.01\,\textrm{s}$
using temporal step sizes that ranged from $\Delta t=\nicefrac{1}{5000}$
to $\Delta t=\nicefrac{1}{80000}$, decreasing by a factor of $2$
at each level. The Eulerian grid is discretized with a step size of
$h=\nicefrac{1}{256}$.

The empirical convergence rates from our temporal refinement are provided
in \prettyref{tab:ConvergenceTimeRefinement}. The immersed boundary
method, as we have implemented it, is formally second-order in space
and first order in time, but, for problems with sharp interfaces that
do not have smooth solutions, it is limited to first-order accuracy
in space and time. Thus for our problem we expect only first order
accuracy. The convergence rates in time given in \prettyref{tab:ConvergenceTimeRefinement}
show first-order convergence in time as is expected. In \prettyref{fig:ErrRef}(a),
we depict the exact values of $E_{p}(q(T);\Delta t)$ for $q=\mathbf{X}$
and $q=\mathbf{u}$. We show $\log_{2}$ in the $x$ and $y$ axes
so that the empirical convergence rates from \prettyref{tab:ConvergenceTimeRefinement}
appear as the slope of the line segments.

\subsubsection{\label{sub:Eulerian-Grid-Refinement}Empirical Estimate of Convergence
Rate in Space}

For this refinement study, we define a measure of error by 
\begin{equation}
E_{p}(q(T);h)=\left\Vert q^{h}(T)-I_{h}^{2h}\left(q^{h/2}(T)\right)\right\Vert _{p},\label{eq:ErrorSpatial}
\end{equation}
which is the error difference at time $t=T$ in a computed quantity,
$q$, using a spatial refinement of a half. In this definition, $I_{h}^{2h}$
is the restriction operator from a fine to a coarse grid. Then, an
empirical estimate for the convergence rate is calculated using 
\begin{equation}
r_{p}(q(T);h)=\log_{2}\left(\frac{E_{p}(q(T);h)}{E_{p}(q(T);\frac{h}{2})}\right).\label{eq:convergenceRate}
\end{equation}
We note that the estimates for convergence rates given by \prettyref{eq:temporal convergence rate}
and \prettyref{eq:convergenceRate} have a fairly simple derivation
using a Taylor series expansion (see \citep{Ferziger2002} or \citep{LeVeque2007}). 

\begin{figure}[H]
\begin{centering}
(a)\includegraphics[width=3in]{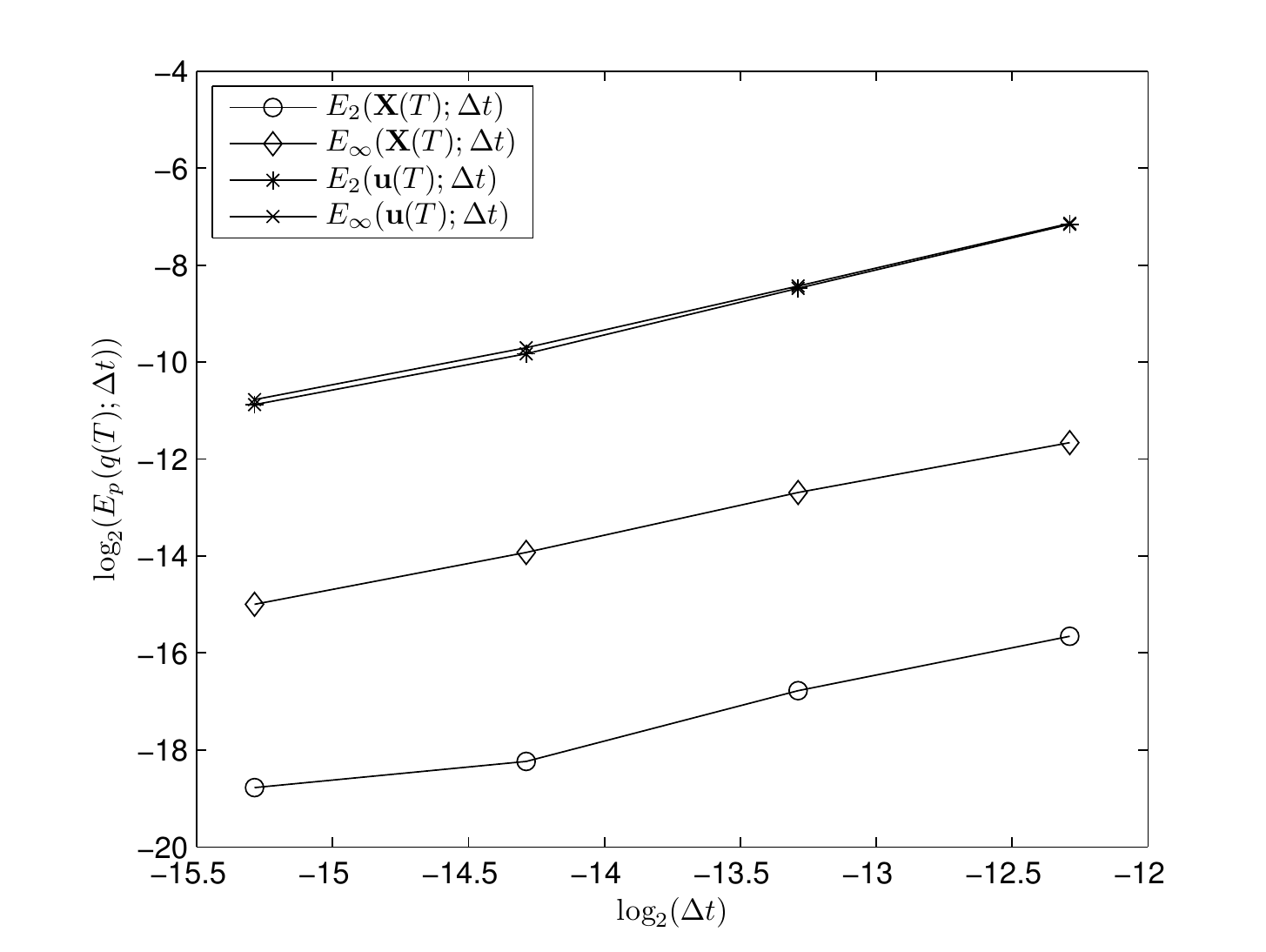}(b)\includegraphics[width=3in]{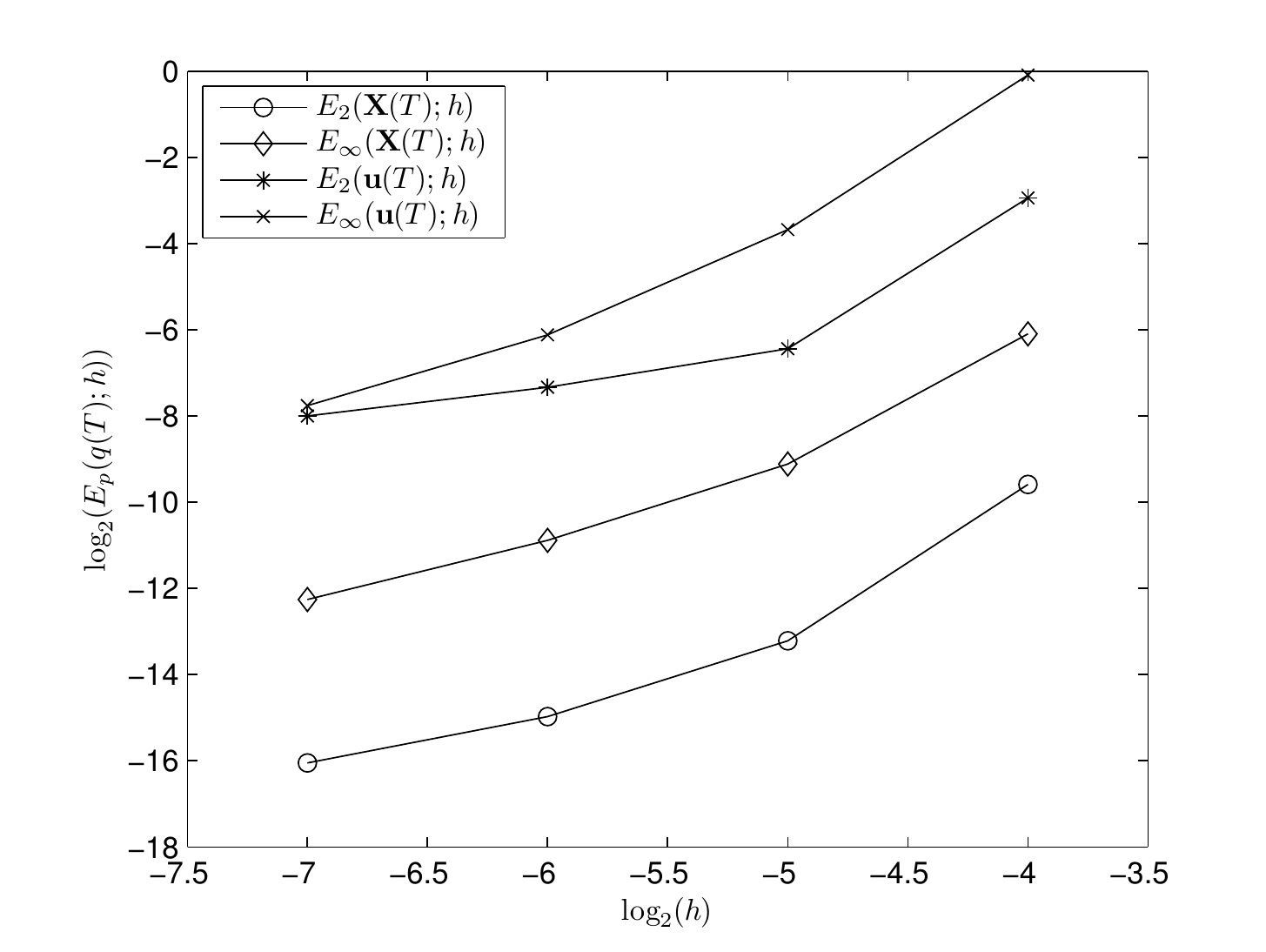}\\
(c)\includegraphics[width=3in]{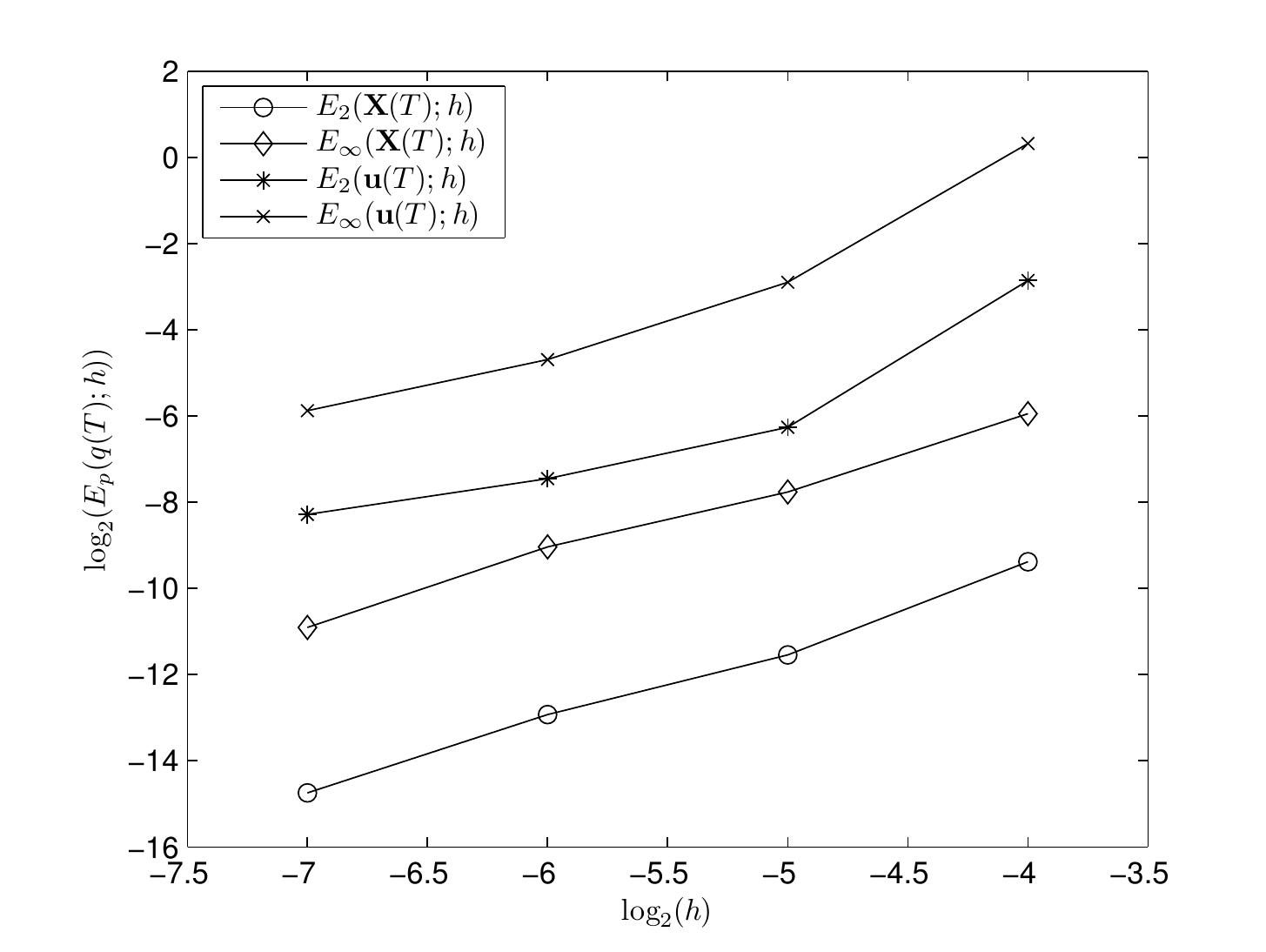}(d)\includegraphics[width=3in]{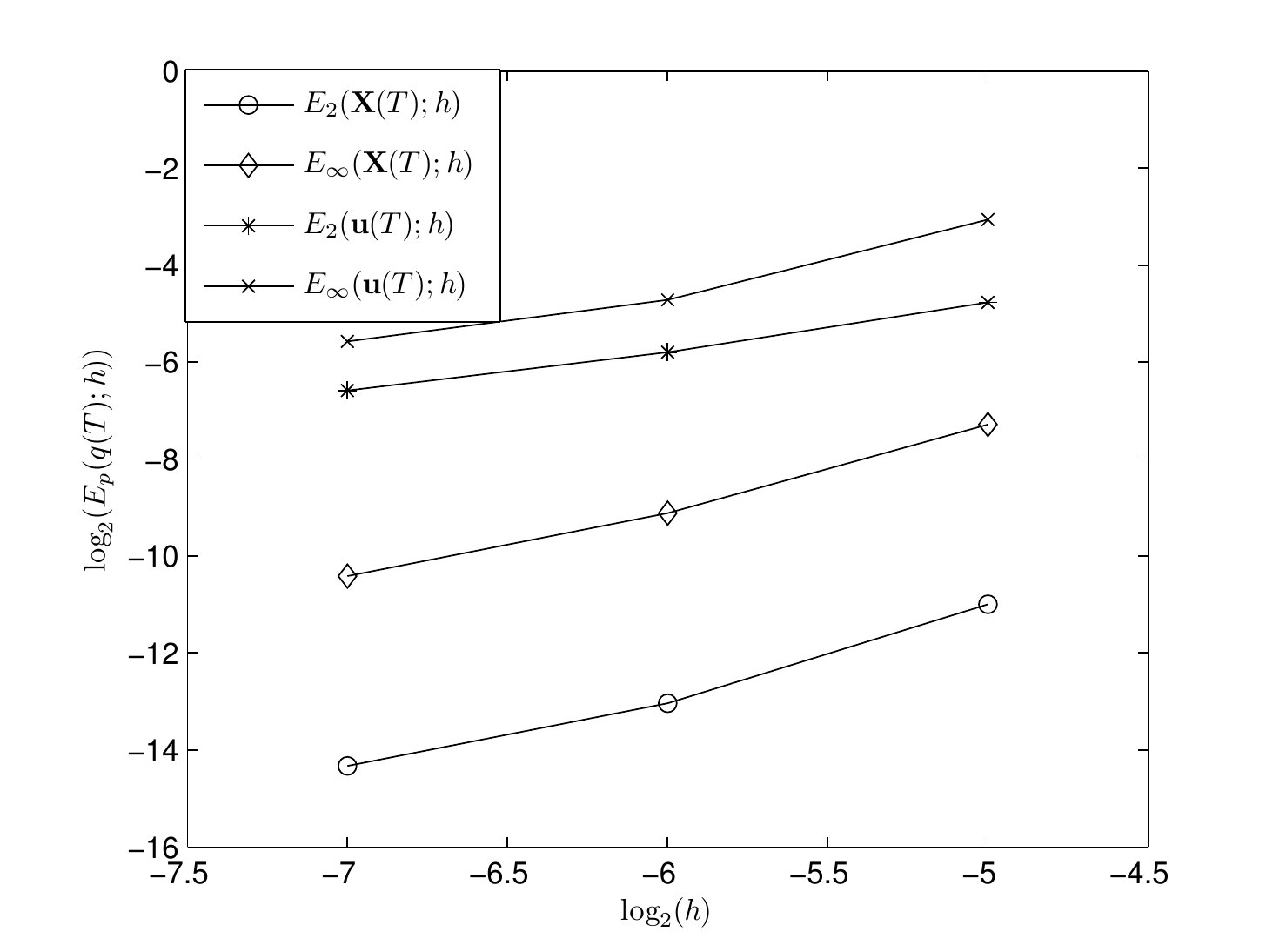}
\par\end{centering}

\caption{\label{fig:ErrRef}Empirical error estimates with (a) temporal refinement:
$E_{p}(q(T);\Delta t)$ is the $p$-norm of the error as defined by
\prettyref{eq:ErrorTemporal}, (b) spatial refinement with constant
density and viscosity: $E_{p}(q(T);h)$ is the $p$-norm of the error
as defined by \prettyref{eq:ErrorSpatial}, (c) spatial refinement
and increased biofilm density, (d) spatial refinement and increased
biofilm viscosity. In all plots, we show $\log_{2}$ in the $x$ and
$y$ axes so that the empirical convergence rate appears as the slope
of the line segments. }
\end{figure}

In the spatial refinement analysis, we did not refine the Lagrangian
grid with the Eulerian grid, so the same number of Lagrangian points
were present in all of the simulations. In addition, full weighting
restriction is used in the definition of the error, \prettyref{eq:ErrorSpatial},
for the error in $\mathbf{u}$. We also used a fixed timestep of $\Delta t=10^{-4}$
until $t=T=0.01\,\textrm{s}$ for all of these simulations. The computed
convergence rates from this refinement are provided in \prettyref{tab:ConvergenceEulerianRefinement}.
The $\infty$-norm convergence rates given in \prettyref{tab:ConvergenceEulerianRefinement}
show greater than first-order convergence in space for the error in
the Lagrangian variable $\mathbf{X}$ and in the Eulerian variable
$\mathbf{u}$. The seemingly large convergence rates for the lower
resolution grids ($h=\frac{1}{16},\,\frac{1}{32},\,\frac{1}{64}$
) can be explained by the fact that using $\omega=\frac{1}{50}$ in
the Dirac delta approximations does not allow the Lagrangian forces
to be adequately represented in the Eulerian grid. This leads to larger
errors in the coarse-grid simulations. Therefore, the best estimates
for the convergence rates are the ones using the three resolutions
all obeying $\omega>h$ given in the 4\textsuperscript{th} and 7\textsuperscript{th}
columns of \prettyref{tab:ConvergenceEulerianRefinement}. In \prettyref{fig:ErrRef}(b),
we depict the exact values of $E_{p}(q(T);h)$ for $q=\mathbf{X}$
and $q=\mathbf{u}$. We show $\log_{2}$ in the $x$ and $y$ axes
so that the empirical convergence rates from \prettyref{tab:ConvergenceEulerianRefinement}
appear as the slope of the line segments. We discuss possibilities
for improvement in the convergence rates later in the conclusion sections. 

We also used a grid refinement analysis to find the empirical convergence
rate with spatial refinement when the density of the biofilm is two
times that of the surrounding fluid. This analysis was done to show
that the first order convergence rate is maintained with the increased
density in the biofilms. The results of this convergence analysis
are shown in \prettyref{tab:ConvergenceEulerianRefinementDensity}
and \prettyref{fig:ErrRef}(c).

Finally, we compute the empirical convergence rates for our 2D simulation
with variable viscosity. In this convergence study, (\prettyref{tab:ConvergenceEulerianRefinementViscSpr}
and \prettyref{fig:ErrRef}(d)) we use a non-dimensionalized value
of biofilm viscosity of $\mu_{max}=500$, which means that the viscosity
at the location of a Lagrangian node is 500 times that of the surrounding
fluid. First-order convergence in space is maintained, even with this
very large biofilm viscosity. 

\begin{table}[H]
\begin{centering}
\begin{tabular}{|l|l|l|l|l|l|l|}
\hline 
q & $r_{2}(q(T);.0002)$ & $r_{2}(q(T);10^{-4})$ & $r_{2}(q(T);\frac{10^{-4}}{2})$ & $r_{\infty}(q(T);.0002)$ & $r_{\infty}(q(T);10^{-4})$ & $r_{\infty}(q(T);\frac{10^{-4}}{2})$\tabularnewline
\hline 
$\mathbf{u}$ & 1.32 & 1.31 & 1.08 & 1.29 & 1.28 & 1.07\tabularnewline
\hline 
$\mathbf{X}$ & 1.12 & 1.45 & 0.54 & 1.04 & 1.23 & 1.07\tabularnewline
\hline 
\end{tabular}
\par\end{centering}

\caption{\label{tab:ConvergenceTimeRefinement}Empirical convergence rates
with temporal refinement. $r_{p}(q(T);dt)$ is the convergence rate
in the variable, $q$, at $t=T$ using the $p$-norm and the three
time steps $dt,\,\nicefrac{dt}{2},\,\nicefrac{dt}{4}$.}
\end{table}

\begin{table}[H]
\begin{centering}
\begin{tabular}{|l|l|l|l|l|l|l|}
\hline 
q & $r_{2}(q(T);\frac{1}{16})$ & $r_{2}(q(T);\frac{1}{32})$ & $r_{2}(q(T);\frac{1}{64})$ & $r_{\infty}(q(T);\frac{1}{16})$ & $r_{\infty}(q(T);\frac{1}{32})$ & $r_{\infty}(q(T);\frac{1}{64})$\tabularnewline
\hline 
$\mathbf{u}$ & 3.51 & 0.89 & 0.65 & 3.59 & 2.44 & 1.65\tabularnewline
\hline 
$\mathbf{X}$ & 3.62 & 1.80 & 1.07 & 3.02 & 1.77 & 1.37\tabularnewline
\hline 
\end{tabular}
\par\end{centering}

\caption{\label{tab:ConvergenceEulerianRefinement}Empirical convergence rates
with spatial refinement. $r_{p}(q(T);h)$ is the convergence rate
in the variable, $q$, at $t=T$ using the $p$-norm and the three
Eulerian step sizes $h,\,\nicefrac{h}{2},\,\nicefrac{h}{4}$.}
\end{table}

\begin{table}[H]
\begin{centering}
\begin{tabular}{|l|l|l|l|l|l|l|}
\hline 
q & $r_{2}(q(T);\frac{1}{16})$ & $r_{2}(q(T);\frac{1}{32})$ & $r_{2}(q(T);\frac{1}{64})$ & $r_{\infty}(q(T);\frac{1}{16})$ & $r_{\infty}(q(T);\frac{1}{32})$ & $r_{\infty}(q(T);\frac{1}{64})$\tabularnewline
\hline 
$\mathbf{u}$ & 3.42 & 1.19 & 0.84 & 3.23 & 1.79 & 1.20\tabularnewline
\hline 
$\mathbf{X}$ & 2.20 & 1.38 & 1.82 & 1.83 & 1.26 & 1.88\tabularnewline
\hline 
\end{tabular}
\par\end{centering}

\caption{\label{tab:ConvergenceEulerianRefinementDensity}Empirical convergence
rates with spatial refinement and increased biofilm density. $r_{p}(q(T);h)$
is the convergence rate in the variable, $q$, at $t=T$ using the
$p$-norm and the three Eulerian step sizes $h,\,\nicefrac{h}{2},\,\nicefrac{h}{4}$.
In this experiment, the density of the biofilm is double that of the
surrounding fluid.}
\end{table}

\begin{table}[H]
\begin{centering}
\begin{tabular}{|l|l|l|l|l|}
\hline 
q & $r_{2}(q(T);\frac{1}{32})$ & $r_{2}(q(T);\frac{1}{64})$ & $r_{\infty}(q(T);\frac{1}{32})$ & $r_{\infty}(q(T);\frac{1}{64})$\tabularnewline
\hline 
$\mathbf{u}$ & 1.04 & 0.78 & 1.67 & 0.86\tabularnewline
\hline 
$\mathbf{X}$ & 2.04 & 1.29 & 1.82 & 1.31\tabularnewline
\hline 
\end{tabular}
\par\end{centering}

\caption{\label{tab:ConvergenceEulerianRefinementViscSpr}Empirical convergence
rates with spatial refinement and increased biofilm viscosity. $r_{p}(q(T);h)$
is the convergence rate in the variable, $q$, at $t=T$ using the
$p$-norm and the three Eulerian step sizes $h,\,\nicefrac{h}{2},\,\nicefrac{h}{4}$.
In this experiment the viscosity of the biofilm is 500 times that
of the surrounding fluid.}
\end{table}

\subsubsection{\label{sub:Time-Step-Stability-Restrictions}Time-Step Stability
Restrictions}

Finally, we investigated the stability of the method computationally
as it depends on the spatial and temporal refinement and the stiffness
of the springs. Analytically, stability applies to a numerical scheme
and not to a computational run, but here we follow Mori and Peskin
in \citep{Mori2008} and give a simple definition of the stability
for each computational run. Using the square of the $2$-norm defined
by \prettyref{eq:EulPnorm} on $\mathbf{u}$ (i.e.~$\left\Vert \mathbf{u}\right\Vert _{p}^{2}$)
gives a value which is proportional to the kinetic energy in the system.
We call the simulation \textit{stable} if magnitude of the total velocity
(as measured by the total kinetic energy) does not have a time of
extreme growth during the simulation. Moreover, this kinetic energy
should remain relatively close to the value of the total kinetic energy
in the case of no biofilm. Using this definition of stability, we
found, through experimentation with many combinations of $h$, $\Delta t$,
and $F_{max}$, that we have timestep restrictions that scale with
the mesh-width, $h$, and with the maximum Lagrangian force, $F_{max}$.
The restrictions are approximately given by 
\[
\triangle t\le C_{1}h\,
\]
and 
\[
\triangle t\le\frac{C_{2}}{F_{max}}\,,
\]
where $C_{1}$ and $C_{2}$ are positive proportionality constants.
Specific values of $C_{1}$ and $C_{2}$ change depending on the parameters
of the simulation. In future simulations, we hope to avoid these timestep
restrictions by using an implicit or semi-implicit method as is done
in \citep{Mori2008} and \citep{Newren2007}. All of the simulations
shown in this work and used in the convergence testing used time-steps
satisfying these two restrictions.

\subsection{\label{sub:Three-Dimensional-Validation}Three-Dimensional Validation}

In this subsection, we provide some numerical evidence validating
the 3D simulations. We first validate in the absence of a biofilm
using the exact laminar flow solution. Then, we validate the multigrid
method in the presence of a biofilm and, finally, we provide the empirical
convergence rates for the simulation in the presence of a biofilm. 

\begin{figure}[h]
\begin{centering}
\includegraphics[width=4in]{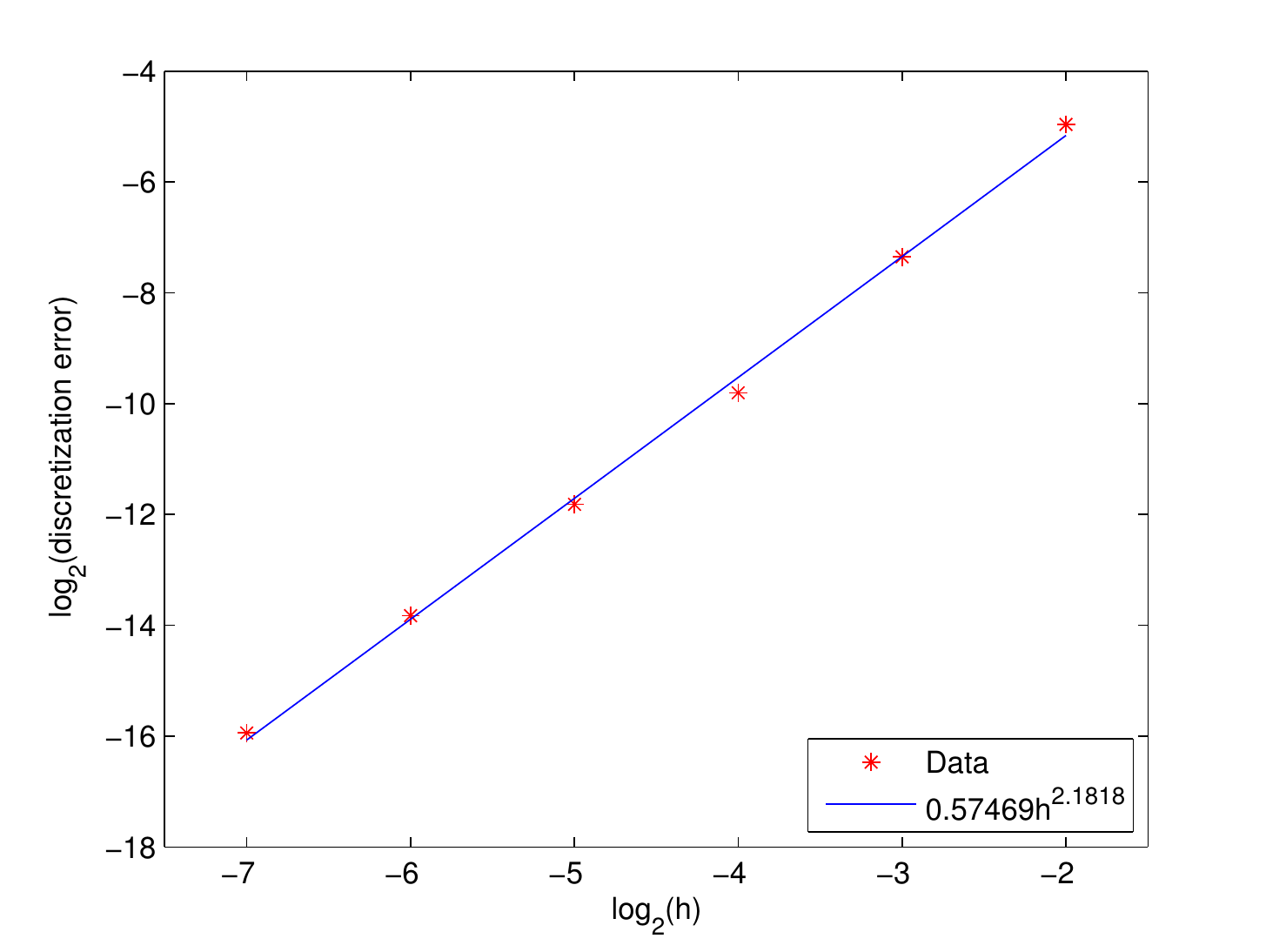}
\par\end{centering}

\caption{\label{fig:Laminar convergence rate}We show the convergence rate
in the laminar flow case. These points are discretization errors for
each of the spatial step sizes ($\frac{1}{4},\,\frac{1}{8},\,\frac{1}{16},\,\frac{1}{32},\,\frac{1}{64},\,\frac{1}{128}$). }
\end{figure}

We first tested the rate of convergence of our method on the laminar
flow case without the interference of a biofilm. To illustrate the
convergence rate in the absence of a biofilm, we started with an initial
velocity profile that is one-half of that of the laminar flow velocity
profile, given by \prettyref{eq:3DlaminarFlow}. We ran the simulation
enough timesteps until the approximate solution converged, with only
discretization error remaining, to the exact solution for six spatial
step sizes, $h=\left\{ \frac{1}{4},\,\frac{1}{8},\,\frac{1}{16},\,\frac{1}{32},\,\frac{1}{64},\,\frac{1}{128}\right\} $.
We computed the discretization error (using the exact laminar solution,
\prettyref{eq:3DlaminarFlow}, for computations) for each of the step
sizes and found that the error is $O(h^{2})$. This can be seen in
\prettyref{fig:Laminar convergence rate}, where on the vertical axis
we have the $\log_{2}$ of the error so that the convergence rate
appears as the slope in the plot. 

Next, in \prettyref{fig:Work Units Decrease 2D}(b), we depict (for
the pressure computation) the work units required to reach the minimum
residual error as a function of allowed levels in the multigrid approach.
This again implies that the multigrid method correctly accelerates
the convergence of our iterative method for the 3D simulations with
a biofilm. Note that there is only a slight reduction in the number
of required work units with the addition of a $6^{\textrm{th}}$ level
in the multigrid, and we saw no reduction with 7 levels, so we use
at most 6 levels in our 3D solvers.

\begin{figure}[h]
\begin{centering}
(a)\includegraphics[width=3in]{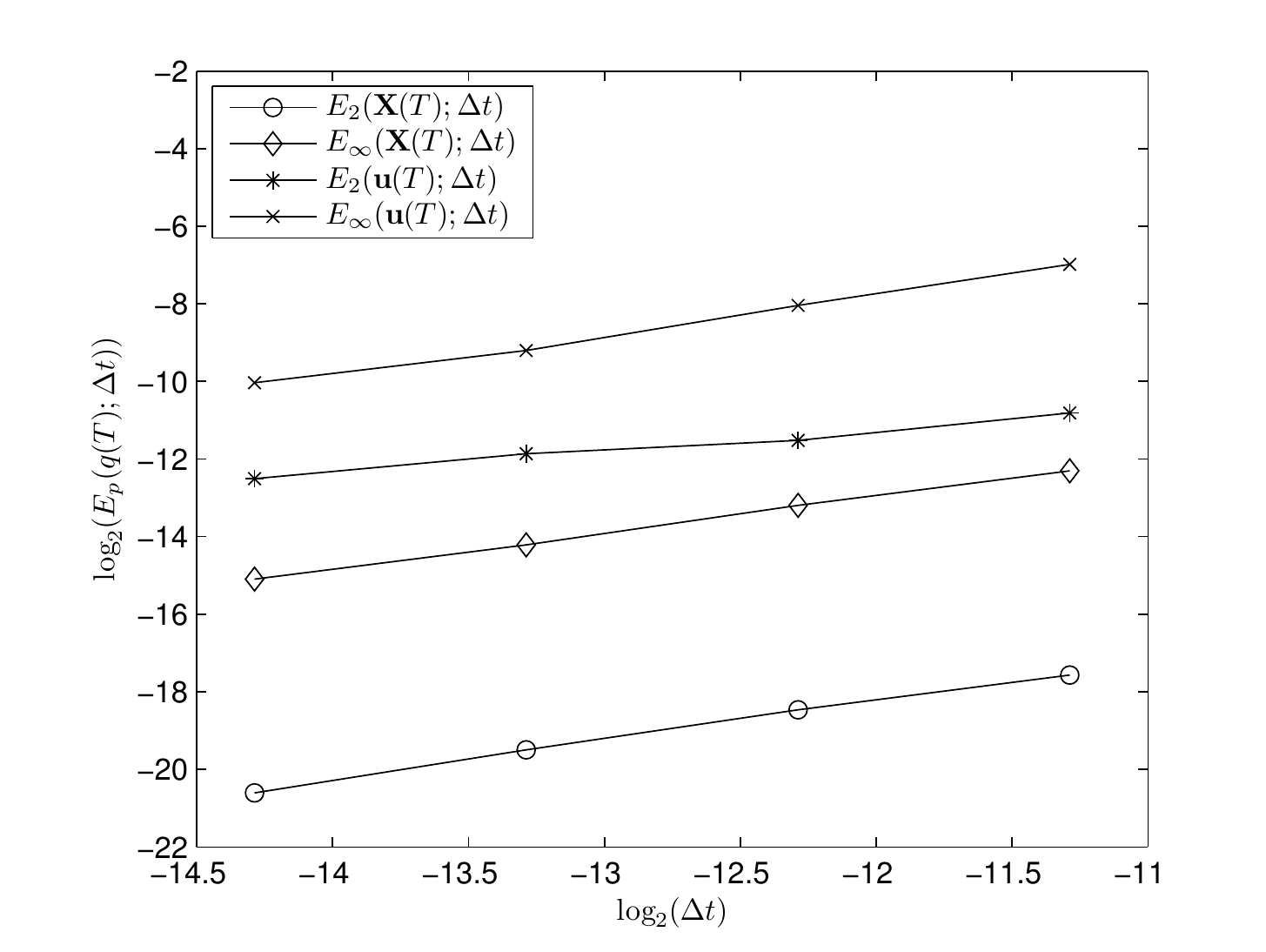}(b)\includegraphics[width=3in]{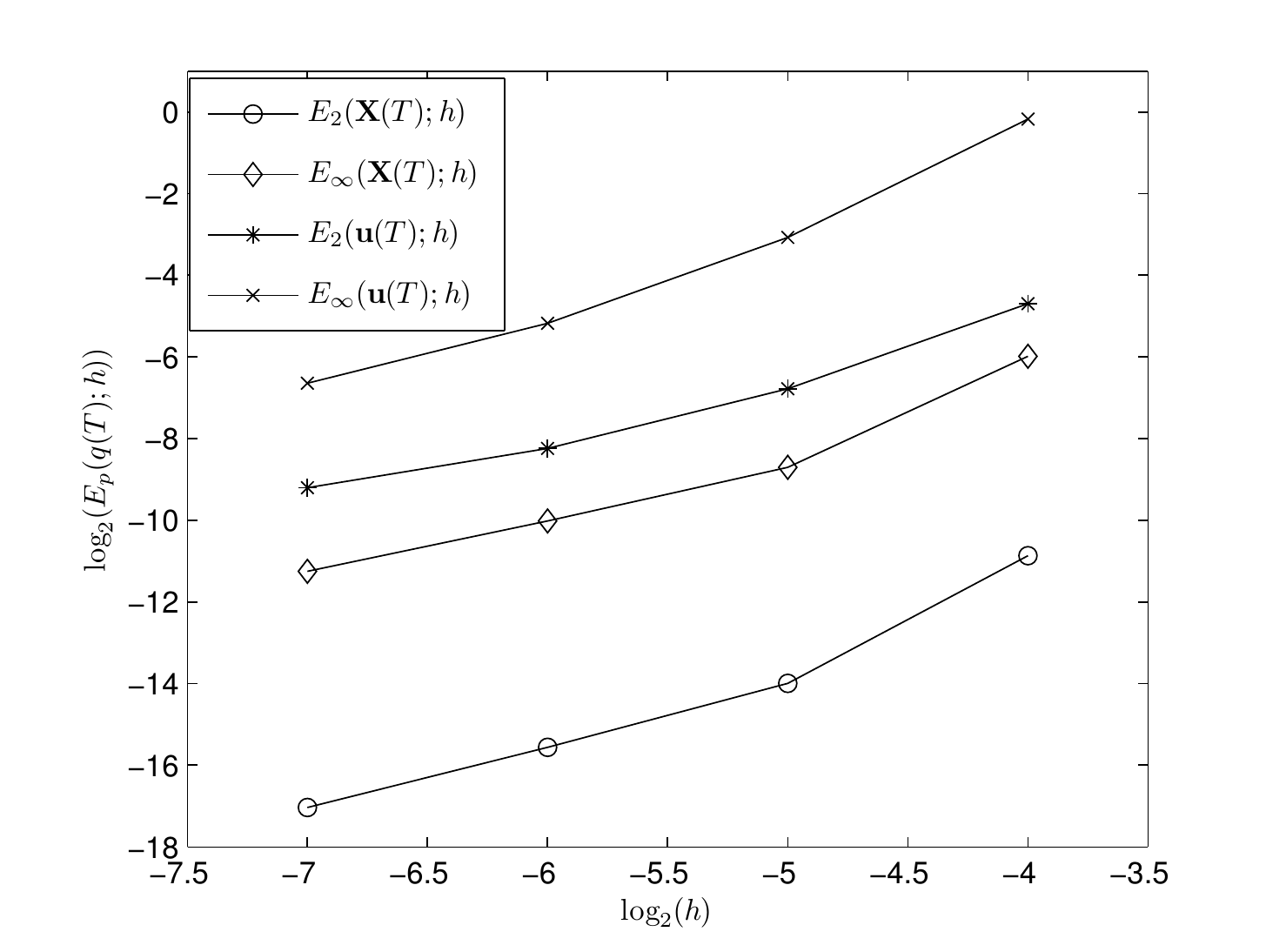}\\
(c)\includegraphics[width=3in]{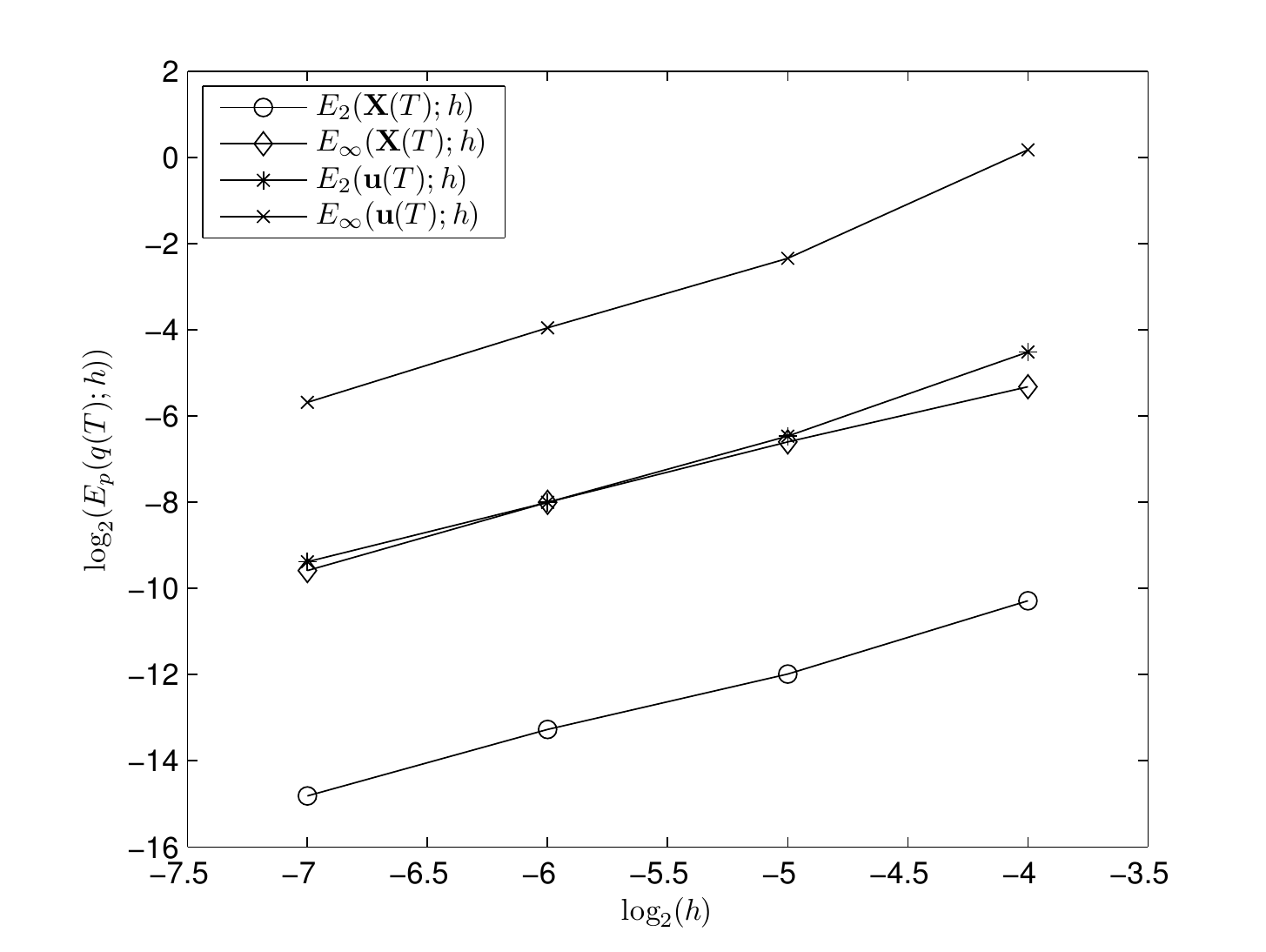}(d)\includegraphics[width=3in]{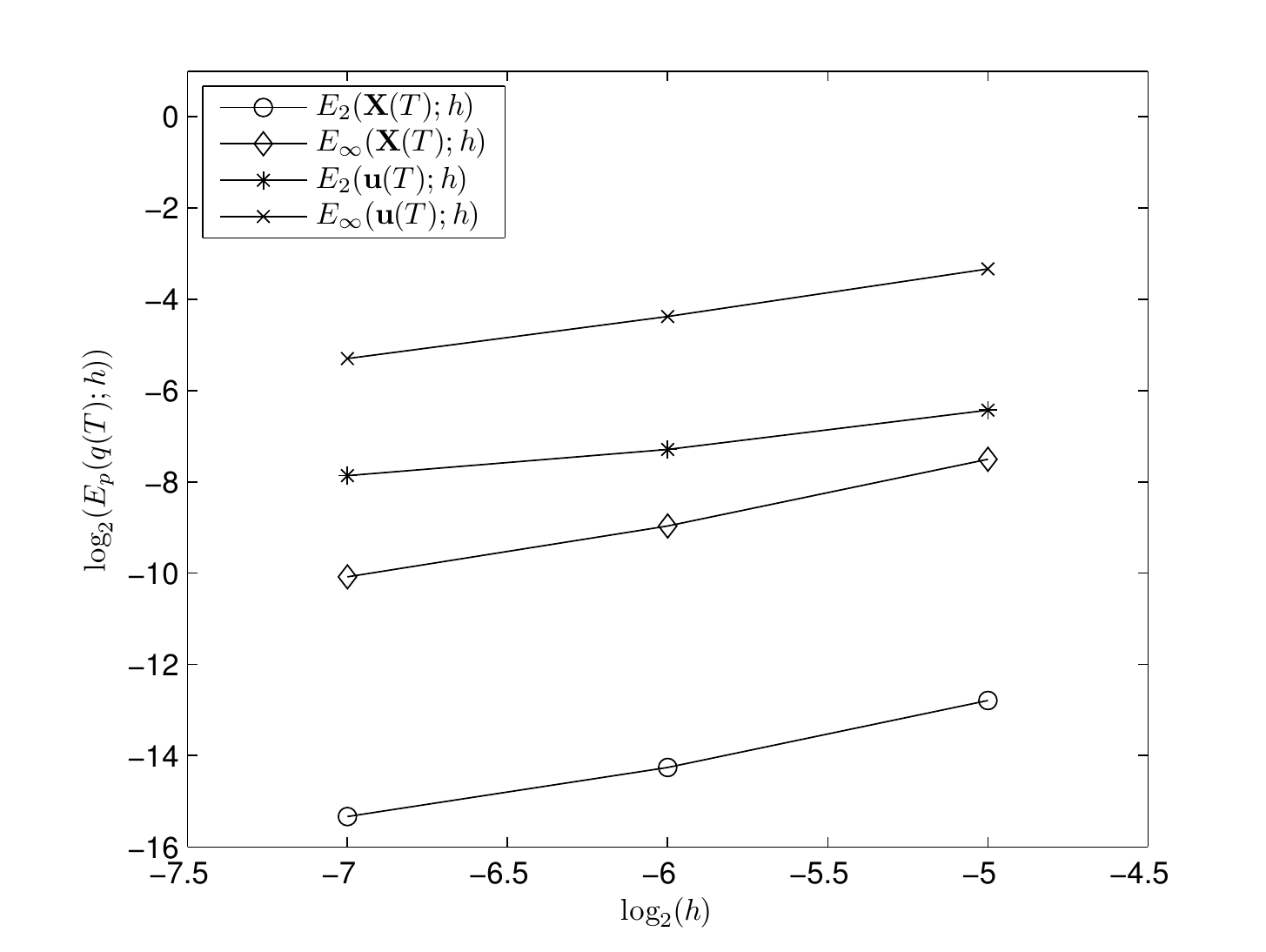}
\par\end{centering}

\caption{\label{fig:EmpConvErr3D}Empirical errors in the 3D simulations with
(a) temporal refinement: $E_{p}(q(T);\Delta t)$ is the $p$-norm
of the error as defined by \prettyref{eq:ErrorTemporal}, (b) spatial
refinement with constant density and viscosity: $E_{p}(q(T);h)$ is
the $p$-norm of the error as defined by \prettyref{eq:ErrorSpatial},
(c) spatial refinement and increased biofilm density, (d) spatial
refinement and increased biofilm viscosity. We show $\log_{2}$ in
the $x$ and $y$ axes so that the empirical convergence rate appears
as the slope of the line segments.}
\end{figure}

Finally, as was done for the 2D case in \prettyref{sub:Eulerian-Grid-Refinement},
we compute the empirical convergence rates for our 3D simulation in
the presence of the biofilm shown in \prettyref{fig:Mushroom-shaped-biofilm2Dsims}(b)
with all of the same fluid parameters used in the 2D analysis. Using
the $p$-norms defined above, we can compute the convergence rates
using \prettyref{eq:temporal convergence rate} and \prettyref{eq:convergenceRate}
(see \prettyref{fig:EmpConvErr3D}(a) and \prettyref{tab:ConvergenceTemporalRefinement3D}).
For the temporal convergence analysis, we used $\omega=\frac{1}{50}$
and $h=\frac{1}{64}$. This analysis resulted in first-order convergence
in all measures except $r_{p}(q(T);\Delta t)$ in which it has an
average convergence rate of about $0.6$. Next, we found empirical
convergence rates for spatial refinement (see \prettyref{tab:ConvergenceEulerianRefinement3D}
and \prettyref{fig:EmpConvErr3D}(b)). As expected, we observe a greater
than first-order convergence rate in both the Eulerian velocity, $\mathbf{u}$,
and the Lagrangian position, $\mathbf{X}$. Next, we conducted a spatial
refinement analysis with a biofilm that has double the density of
the surrounding fluid (see \prettyref{tab:ConvergenceEulerianRefinementDensity3D}
and \prettyref{fig:EmpConvErr3D}(c)). Finally, we did the spatial
refinement study for our simulations with $\mu_{max}=500$, and again
achieved first-order spatial convergence (see \prettyref{tab:ConvergenceEulerianRefinementViscSpr3D}
and \prettyref{fig:EmpConvErr3D}(d)). 

\begin{table}[H]
\begin{centering}
\begin{tabular}{|l|l|l|l|l|l|l|}
\hline 
q & $r_{2}(q(T);.0004)$ & $r_{2}(q(T);.0002)$ & $r_{2}(q(T);10^{-4})$ & $r_{\infty}(q(T);.0004)$ & $r_{\infty}(q(T);.0002)$ & $r_{\infty}(q(T);10^{-4})$\tabularnewline
\hline 
$\mathbf{u}$ & 0.71 & 0.35 & 0.64 & 1.04 & 1.18 & 0.83\tabularnewline
\hline 
$\mathbf{X}$ & 0.90 & 1.03 & 1.11 & 0.90 & 1.02 & 0.88\tabularnewline
\hline 
\end{tabular}
\par\end{centering}

\caption{\label{tab:ConvergenceTemporalRefinement3D}Empirical convergence
rates in the 3D simulations with temporal refinement are shown for
$\mathbf{u}$ and $\mathbf{X}$. $r_{p}(q(T);\Delta t)$ is the convergence
rate in the variable, $q$, at $t=T$ using the $p$-norm and the
three Eulerian step sizes $\Delta t,\,\nicefrac{\Delta t}{2},\,\nicefrac{\Delta t}{4}$.}
\end{table}

\begin{table}[H]
\begin{centering}
\begin{tabular}{|l|l|l|l|l|l|l|}
\hline 
q & $r_{2}(q(T);\frac{1}{16})$ & $r_{2}(q(T);\frac{1}{32})$ & $r_{2}(q(T);\frac{1}{64})$ & $r_{\infty}(q(T);\frac{1}{16})$ & $r_{\infty}(q(T);\frac{1}{32})$ & $r_{\infty}(q(T);\frac{1}{64})$\tabularnewline
\hline 
$\mathbf{u}$ & 2.09 & 1.44 & 1.02 & 2.89 & 2.11 & 1.47\tabularnewline
\hline 
$\mathbf{X}$ & 3.12 & 1.57 & 1.47 & 2.70 & 1.33 & 1.23\tabularnewline
\hline 
\end{tabular}
\par\end{centering}

\caption{\label{tab:ConvergenceEulerianRefinement3D}Empirical convergence
rates in the 3D simulations with spatial refinement are shown for
$\mathbf{u}$ and $\mathbf{X}$. $r_{p}(q(T);h)$ is the convergence
rate in the variable, $q$, at $t=T$ using the $p$-norm and the
three Eulerian step sizes $h,\,\nicefrac{h}{2},\,\nicefrac{h}{4}$.}
\end{table}

\begin{table}[H]
\begin{centering}
\begin{tabular}{|l|l|l|l|l|l|l|}
\hline 
q & $r_{2}(q(T);\frac{1}{16})$ & $r_{2}(q(T);\frac{1}{32})$ & $r_{2}(q(T);\frac{1}{64})$ & $r_{\infty}(q(T);\frac{1}{16})$ & $r_{\infty}(q(T);\frac{1}{32})$ & $r_{\infty}(q(T);\frac{1}{64})$\tabularnewline
\hline 
$\mathbf{u}$ & 1.95 & 1.52 & 1.40 & 2.52 & 1.61 & 1.73\tabularnewline
\hline 
$\mathbf{X}$ & 1.7 & 1.28 & 1.55 & 1.29 & 1.40 & 1.56\tabularnewline
\hline 
\end{tabular}
\par\end{centering}

\caption{\label{tab:ConvergenceEulerianRefinementDensity3D}Empirical convergence
rates with 3D spatial refinement and increased biofilm density. $r_{p}(q(T);h)$
is the convergence rate in the variable, $q$, at $t=T$ using the
$p$-norm and the three Eulerian step sizes $h,\,\nicefrac{h}{2},\,\nicefrac{h}{4}$.
In this experiment, the density of the biofilm is double that of the
surrounding fluid.}
\end{table}

\begin{table}[H]
\begin{centering}
\begin{tabular}{|l|l|l|l|l|}
\hline 
q & $r_{2}(q(T);\frac{1}{64})$ & $r_{2}(q(T);\frac{1}{128})$ & $r_{\infty}(q(T);\frac{1}{64})$ & $r_{\infty}(q(T);\frac{1}{128})$\tabularnewline
\hline 
$\mathbf{u}$ & 0.86 & 0.57 & 1.04 & 0.92\tabularnewline
\hline 
$\mathbf{X}$ & 1.47 & 1.07 & 1.48 & 1.09\tabularnewline
\hline 
\end{tabular}
\par\end{centering}

\caption{\label{tab:ConvergenceEulerianRefinementViscSpr3D}Empirical convergence
rates with 3D spatial refinement and increased biofilm viscosity.
$r_{p}(q(T);h)$ is the convergence rate in the variable, $q$, at
$t=T$ using the $p$-norm and the three Eulerian step sizes $h,\,\nicefrac{h}{2},\,\nicefrac{h}{4}$.
In this experiment, the viscosity of the biofilm is 500 times that
of the surrounding fluid.}
\end{table}

This concludes our validation section, and we now present the results
of our numerical simulations.

\section{\label{sec:Simulations-Results}Simulations Results}

In this section, we present the results of our numerical simulations.
First, we briefly discuss the reality of elastic forces in biofilms.
Then we provide 2D results in \prettyref{sub:Two-Dimensional-Simulations}
and 3D results in \prettyref{sub:Three-Dimensional-Simulations}.

\subsection{\label{sub:Discussion-of-Elastic}Discussion of Elastic Maximum Force,
$F_{max}$}

We now provide a brief discussion of the physical reality of the values
of $F_{max}$ used in our simulations. The \textit{cohesive strength}%
\footnote{The cohesive strength is a measure of the forces that interconnect
the biofilm's cells.%
}\negthinspace{}\negthinspace{} in biofilms has been found experimentally
to be highly heterogeneous, with repeated experimental measurements
on the same biofilm yielding vastly different strength measurements.
For example, 49 cohesive strength measurements taken on only two samples
of \textit{Staphylococcus epidermidis }yielded measurements between
61-5182 Pa \citep{Aggarwal2010}. These biofilms were grown on a $22\,\textrm{mm}$
diameter disc rotating at 75 $\textrm{\ensuremath{\nicefrac{rot}{min}}}$
so the fastest speed, $\sim86\,\nicefrac{\textrm{mm}}{\textrm{s}}$,
was at the perimeter of the disc (i.e. very slow flow growth conditions).
The \textit{adhesive}%
\footnote{The adhesive strength is a measure of the forces that connect a biofilm
to the surface.%
}\negthinspace{}\negthinspace{} and cohesive strengths have also
been shown to vary significantly with changes in growth conditions
such as flow rate and nutrient concentration. Changes in these growth
conditions influence the amount of ECM production in the biofilm as
well as the compactness of the biofilm, which has a direct effect
on its strengths \citep{Chen2005effectsAdhesiveStrength,Ohashi1996,chen1998directBiofilmStrengthTube}.
We note here that the required values we find for $F_{max}$ for the
biofilms to remain attached in our 2D and 3D simulations are consistent
with the cohesive strength measurements provided in \citep{Aggarwal2010}.
Since the diameter \textit{Staphylococcus epidermidis} is about $1\,\mu$m,
in 3D, we multiply the cohesive strengths by $1\,\mu\textrm{m}^{2}$
to get an approximation for the range of forces on the surface area
of one cell. Using the range of 61-5182 Pa yields a range of forces
from $6.1\times10^{-11}$ N to $5.18\times10^{-9}$ N. In 2D, we multiply
the cohesive strengths by the cell diameter to get a rough approximation
for the range of forces on the surface perimeter surrounding one cell.
Using the range of 61-5182 Pa yields a range of forces from $6.1\times10^{-5}$
N to $5.18\times10^{-3}$ N. Our values for $F_{max}$ are at the
low end of these ranges. The actual strength of the biofilm is most
likely larger than our $F_{max}$ values since the positional data
was from a biofilm that was not fragmenting in the flow conditions
in which it was grown, and we used the same flow conditions our simulations.
Thus, in order to see detachment under these flow conditions, we had
to lower the value of $F_{max}$. We could alternatively increase
the flow rate to necessitate a larger $F_{max}$ requirement to avoid
detachment. One eventual goal of this work is that, if the approximate
value of $F_{max}$ is known for a particular type of biofilm, then
our simulations can be used to predict the flow rates required to
break different shaped biofilms.

\subsection{Two-Dimensional Simulations\label{sub:Two-Dimensional-Simulations}}

In this section, we provide results from our 2D simulations, which
represent a cross-section of a biofilm attached to the inside of a
tube and subjected to fluid flow in a computational domain of $150\,\mu\textrm{m}$
by $50\,\mu\textrm{m}$. The parameters for our simulations are given
in \prettyref{tab: 2D sim params}.

\begin{table}
\centering{}%
\begin{tabular}{|l|c|}
\hline 
\multicolumn{2}{|c|}{Parameter Values for the Simulations}\tabularnewline
\hline 
\hline 
Tube Radius & $25\times10^{-6}\,\textrm{m}$\tabularnewline
\hline 
Fluid Dynamic Viscosity & $1.0\times10^{-3}\,\nicefrac{\textrm{kg}}{\textrm{m\ensuremath{\cdot}s}}$\tabularnewline
\hline 
Fluid Density & $998\,\nicefrac{\textrm{kg}}{\textrm{m}^{3}}$\tabularnewline
\hline 
Maximum Fluid Velocity & $10^{-3}\,\nicefrac{\textrm{m}}{\textrm{s}}$\tabularnewline
\hline 
\end{tabular}\\
\caption{\label{tab: 2D sim params}The values of parameters used in the 2D
simulations.}
\end{table}

In all simulations, we implement a breaking condition on the springs
of two times the rest length. The initial configuration for the biofilm
in these simulations is shown in \prettyref{fig:Mushroom-shaped-biofilm2Dsims}.
The spring connections between Lagrangian nodes are put in place at
the beginning of the simulation with every node connected to every
other node less than $d_{c}$ away (the reason for this connection
distance is given above in \prettyref{sub:SimSetup}). The mushroom
shaped biofilm has a height of about $8.5\,\mu\textrm{m}$ and width
of about $8\,\mu\textrm{m}$ (width of about $2\,\mu\textrm{m}$ at
the thinnest part) . We use a non-dimensionalized $\omega=\frac{1}{100}$
to match the radius of \textit{Staphylococcus epidermidis} and choose
$h=\frac{1}{128}$ in all of the simulations shown, so that $\omega>h$. 

\begin{figure}[H]
\begin{centering}
(a)\includegraphics[width=3in]{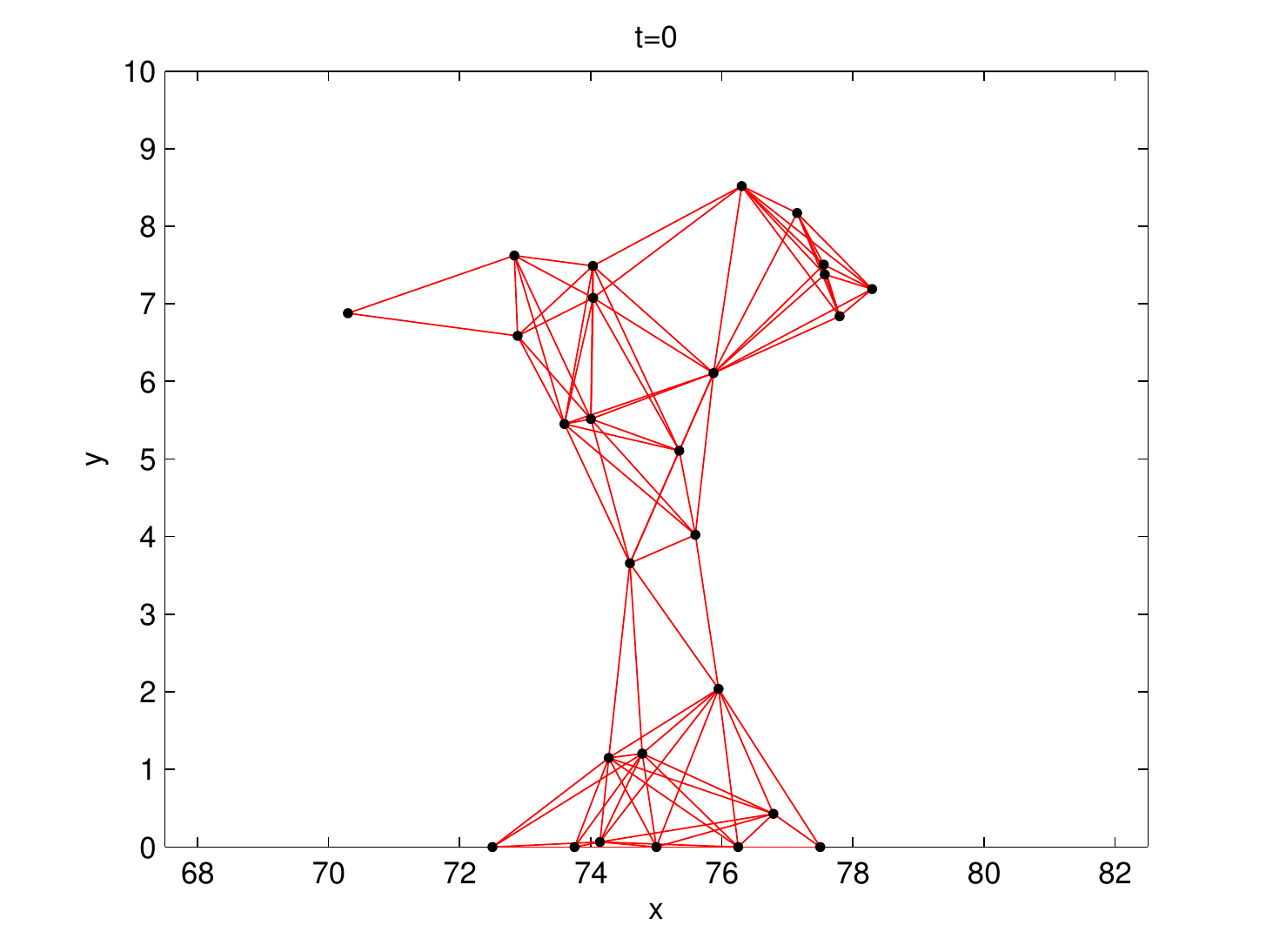}(b)\includegraphics[width=3in]{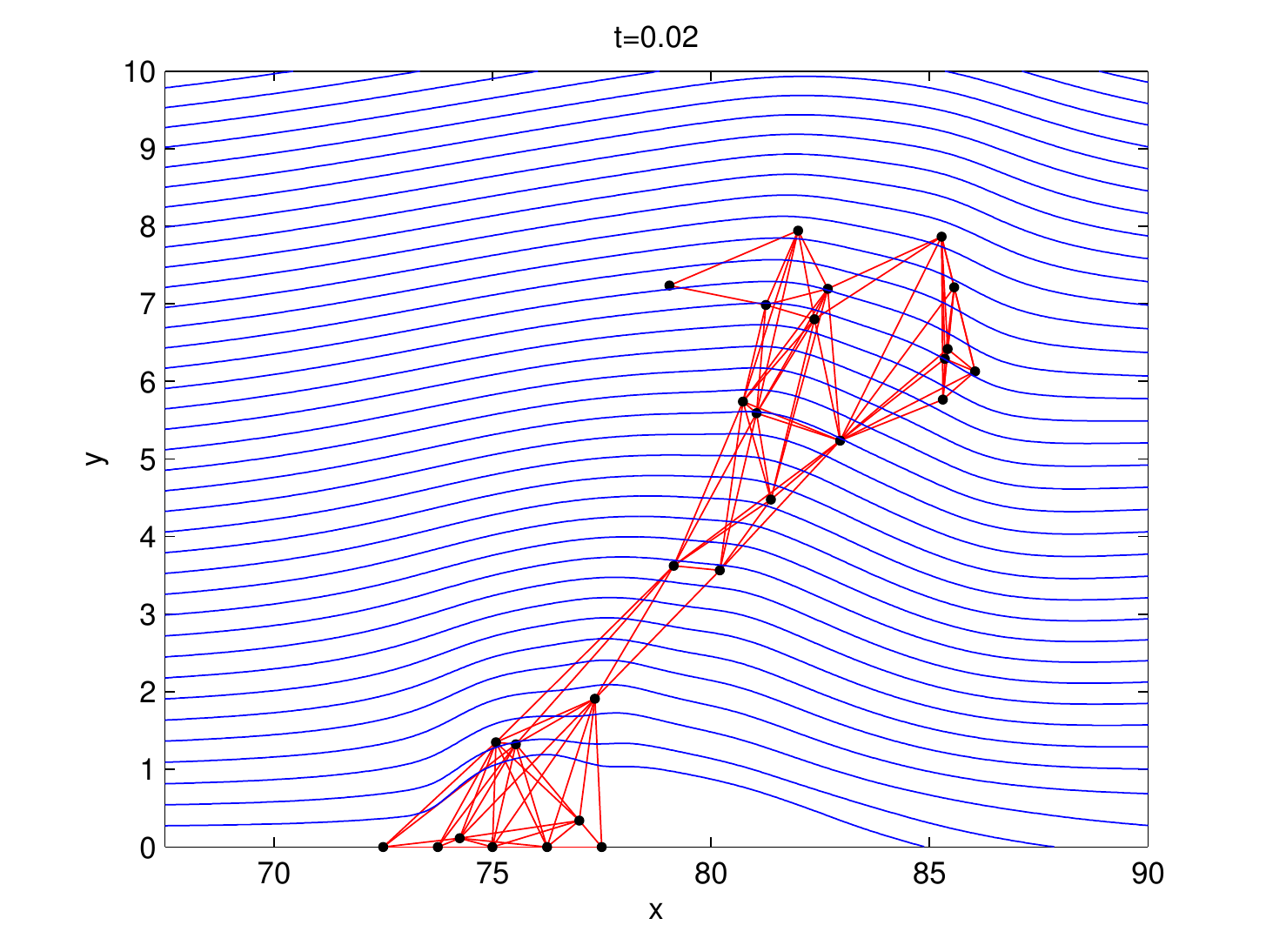}\\
(c)\includegraphics[width=3in]{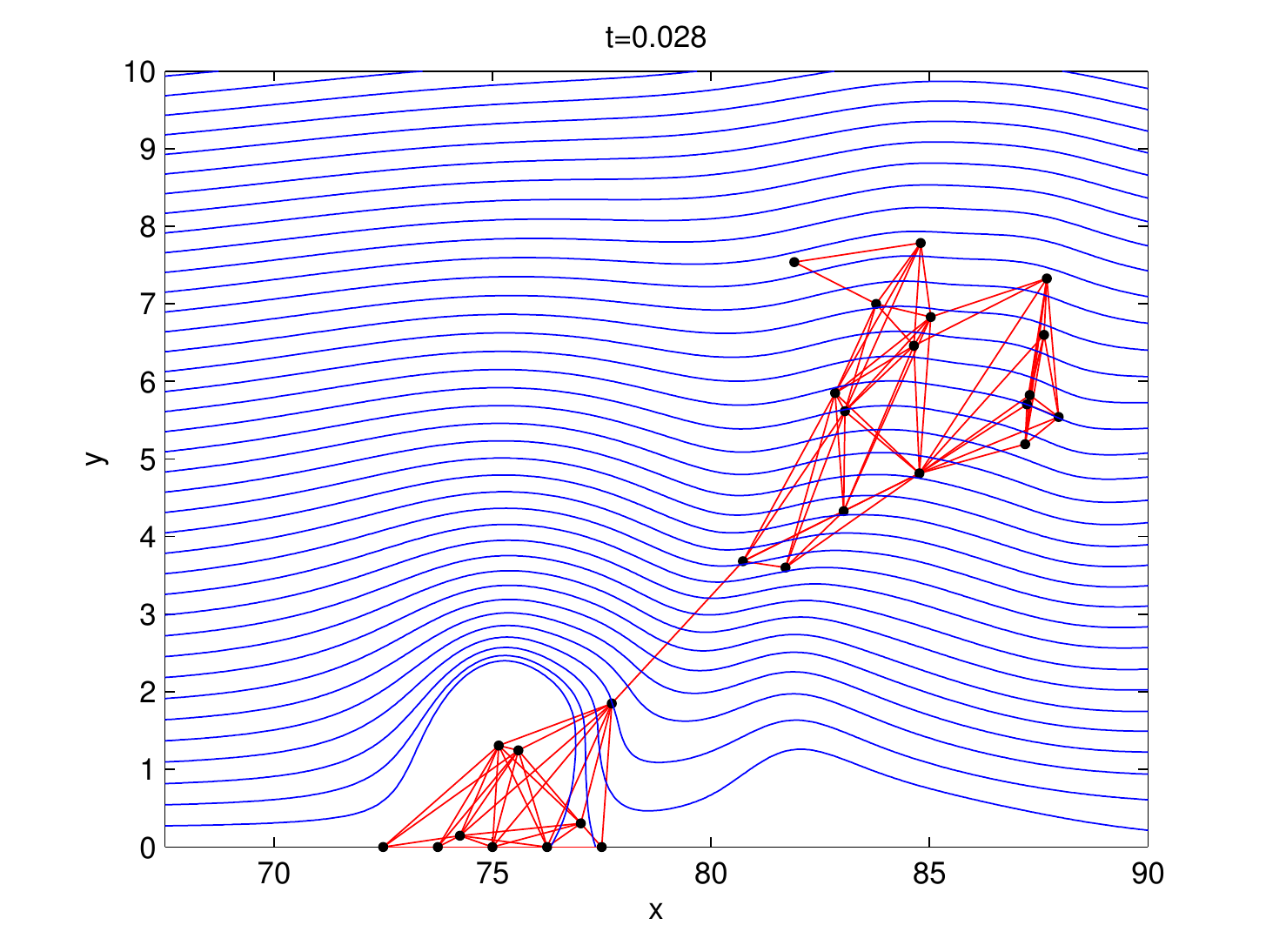}(d)\includegraphics[width=3in]{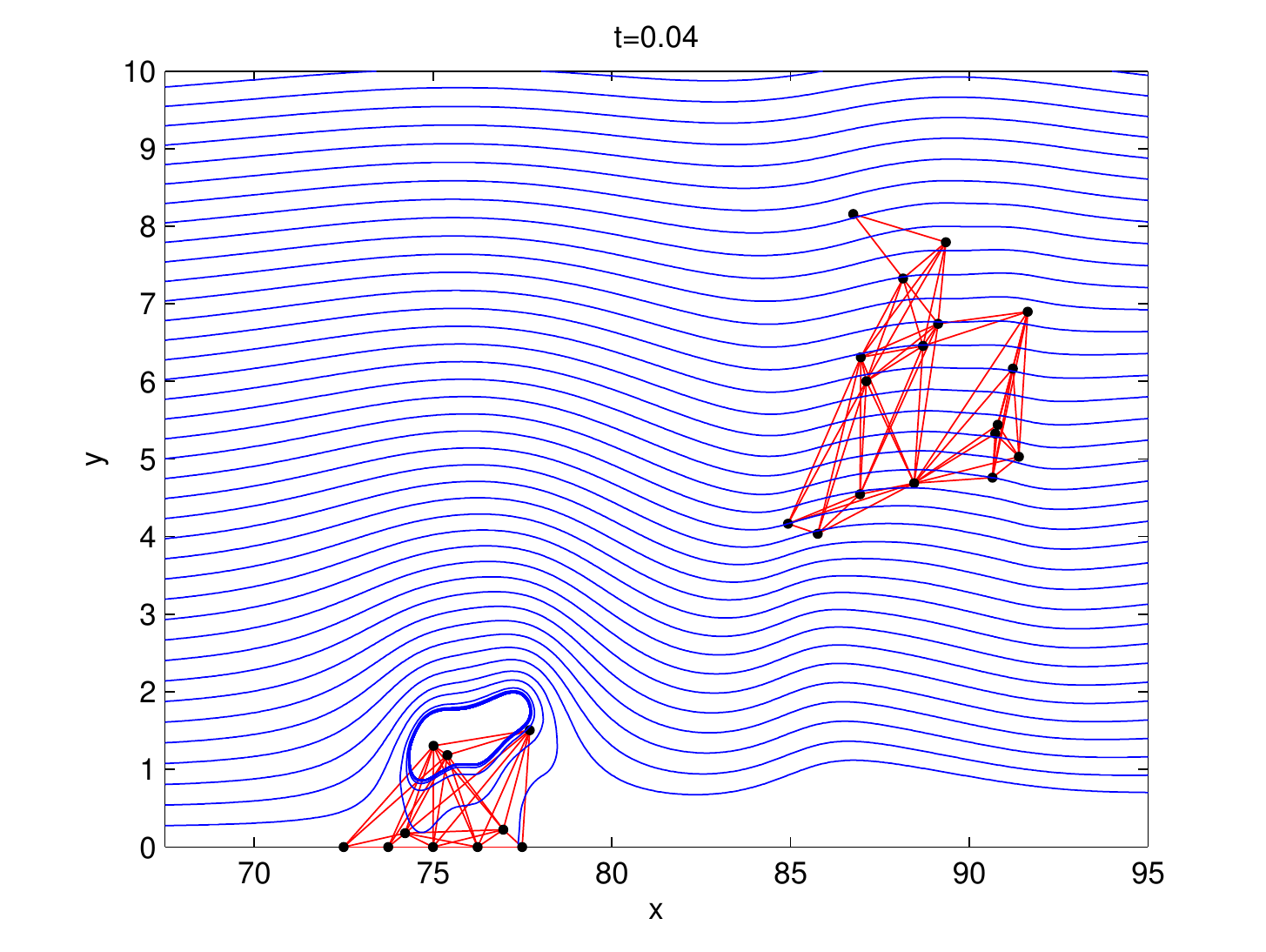}
\par\end{centering}

\caption{\label{fig:ShroomNoaddl_f175}2D Simulation of a mushroom shaped biofilm
with the same density as the surrounding fluid. Time is in seconds
and the distance is in microns. In this simulation, $\rho_{0}=998\,\nicefrac{\textrm{kg}}{\textrm{m}^{3}},\:\rho_{b}=0,\: F_{max}=5.00\times10^{-7}\,\textrm{N}$.
The streamlines follow the velocity field. In this simulation, the
top of the biofilm stretches in the flow, and the top breaks off as
the connections in the the middle separate as they exceed the breaking
criteria of twice the rest length. As expected in a laminar shear
flow, the broken piece then tumbles end over end through the flow.}
\end{figure}

In the first simulation, the maximum spring force, $F_{max}$, is
set to $5.00\times10^{-7}\,\textrm{N}$, and the results are provided
in \prettyref{fig:ShroomNoaddl_f175}. The biofilm bends over in the
flow, and the connections in the thin part of the biofilm break as
they stretch too far. The streamlines in (b), (c), and (d) of \prettyref{fig:ShroomNoaddl_f175}
and in all of the other 2D simulation plots follow the trajectories
given by the velocity field, $\mathbf{u}$. 

We point out that the values of the spring constants are well within
physically realistic values (see the discussion in \prettyref{sub:Discussion-of-Elastic}),
although, in this work, we have chosen these values for the qualities
they give to the simulations rather than experimental evidence of
the elastic strength of biofilms. For example, in these 2D simulations,
we investigated several simulation runs with various spring constants
until we obtained those that exhibited the above described behaviors.

\begin{figure}[H]
\begin{centering}
(a)\includegraphics[width=3in]{shroom2Dnondim_noAddl_t_0_028_f10}(b)\includegraphics[width=3in]{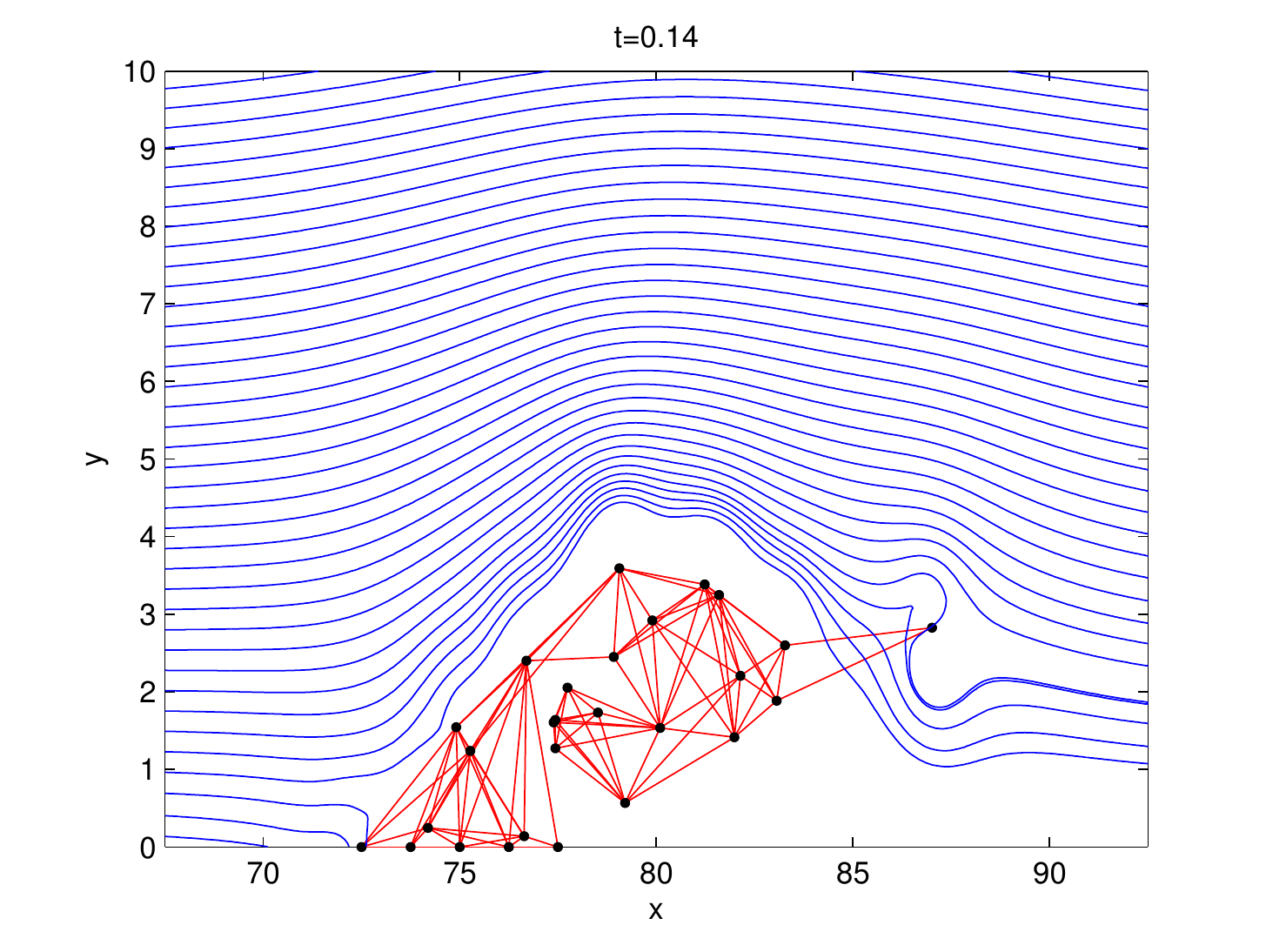}\\
(c)\includegraphics[width=3in]{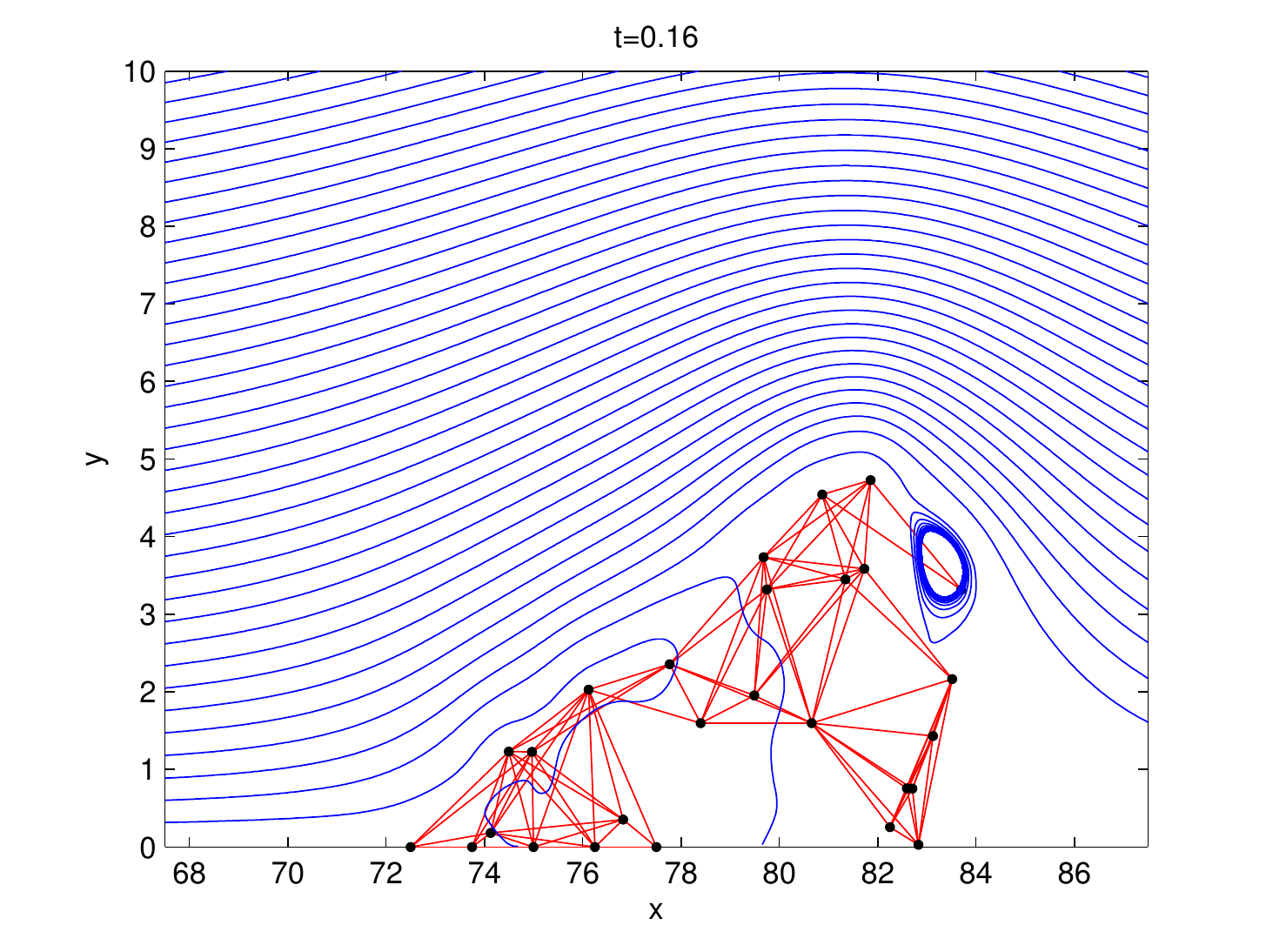}(d)\includegraphics[width=3in]{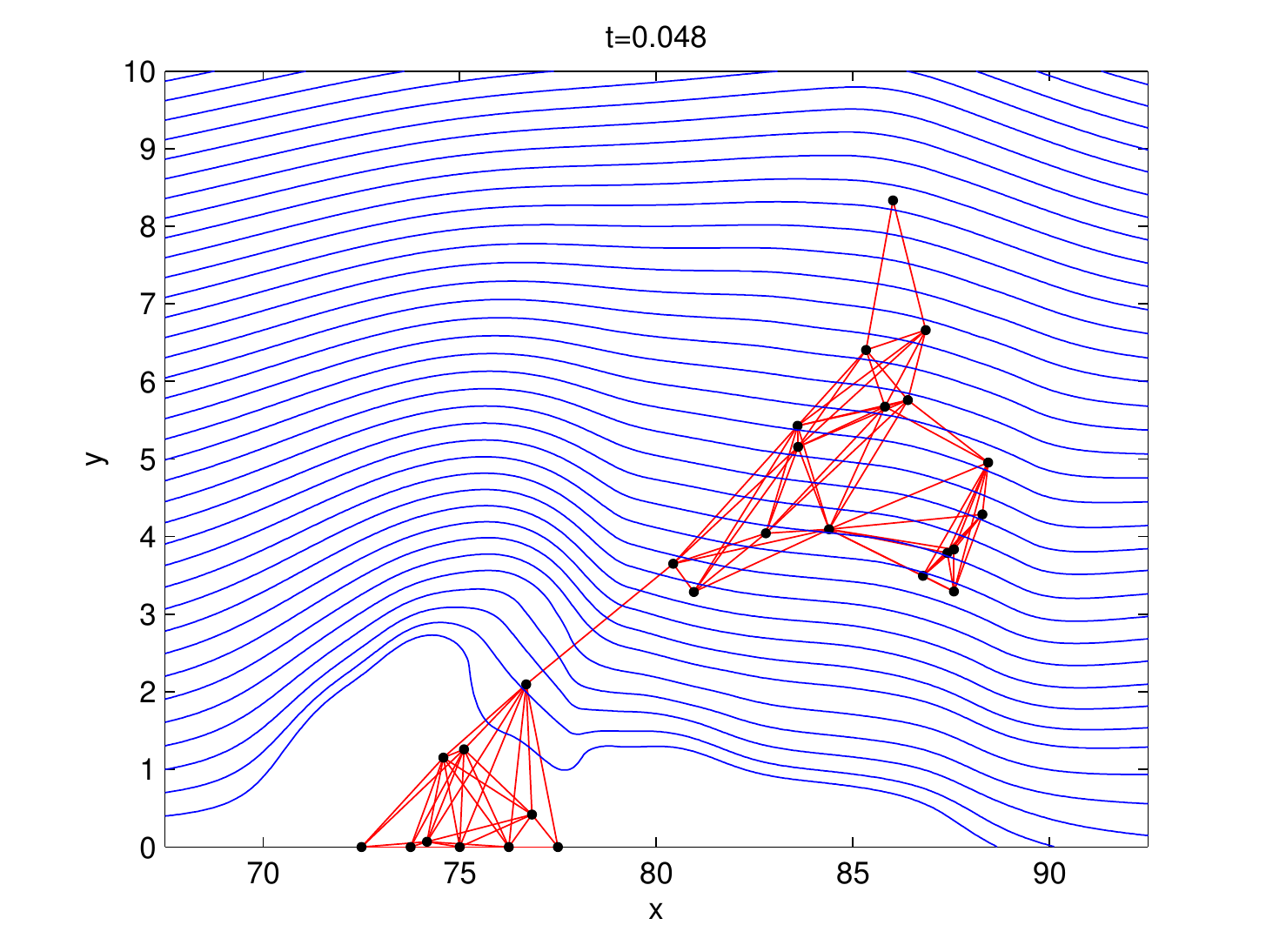}
\par\end{centering}

\caption{\label{fig:ShroomAddl} Snapshots of cell distributions resulting
from 2D simulation of mushroom shaped biofilm with four different
property configurations. In (a), the biofilm has the same density
and viscosity as the surrounding fluid (just before detachment), (b)
increased $F_{max}$, (c) density twice as large as the surrounding
fluid density, and (d) viscosity is $500\times$ that of the surrounding
fluid (just before detachment). The streamlines follow the velocity
field.}
\end{figure}

Next, we conducted a simulation of the same mushroom shaped biofilm
with all of the same parameter values, but we gave the biofilm additional
density of $\rho_{b}=998\,\nicefrac{\textrm{kg}}{\textrm{m}^{3}}$
compared to the ambient fluid. We know this density is larger than
what is seen in actual biofilms (at most 20\% greater density than
water \citep{Masuda1991,Ro1991}), but we chose it to show an exaggerated
example of increasing the biofilm density. These result is provided
in \prettyref{fig:ShroomAddl}(b) and illustrates that the added density
essentially adds momentum to the biofilm. This additional momentum
causes the biofilm to curl over into the slower flow region and thus
prevents detachment. 

Finally, we conducted a simulation of the same mushroom shaped biofilm
with all of the same parameter values as the first simulation, but
we increased $F_{max}$ to $5.00\times10^{-6}\,\textrm{N}$. The effect
of these stronger springs is that the thin part of the biofilm does
not stretch enough to break the connections. The result is depicted
in \prettyref{fig:ShroomAddl}(c). We can see from these simulations
that either increasing the biofilm density or strengthening the springs
causes similar results, but, with the increased density, the biofilm
has more of a curling action. 

In the next simulation, we use all of the same parameters described
in the first simulation, with the addition that the biofilm has a
500$\times$ larger viscosity than the surrounding fluid, so $\mu_{max}=0.5\,\nicefrac{\textrm{kg}}{\textrm{m}\cdot\textrm{s}}$.
Comparing simulation results illustrated in \prettyref{fig:ShroomAddl}(d)
and \prettyref{fig:ShroomAddl}(a), which show the biofilm configurations
just before detachment, we observe a longer time until detachment
in the high viscosity case. This is the expected outcome of increasing
the viscosity in the biofilm. We note here that we used $\omega=\frac{1}{100}$
in the equation for $\mu(\mathbf{x})$ in \prettyref{eq:ViscExp}
because we wanted to spread additional viscosity over the same region
that the elastic forces are spread to. We achieve an even longer detachment
time in the simulation by widening the influence of additional viscosity
by using, for example $\omega=\frac{1}{50}$. 

\subsection{Three-Dimensional Simulations\label{sub:Three-Dimensional-Simulations}}

In this section, we provide results from our 3D simulations, which
use the same parameter values as in the two dimensional simulations
(see \prettyref{tab: 2D sim params}). The difference is that the
simulation in 3D represents flow through a square shaped tube with
a side length of $50\,\mu\textrm{m}$. Note that these 3D simulations
reproduce qualitatively the same results as in the 2D ones.

\begin{figure}[H]
\begin{centering}
(a)\includegraphics[width=3in]{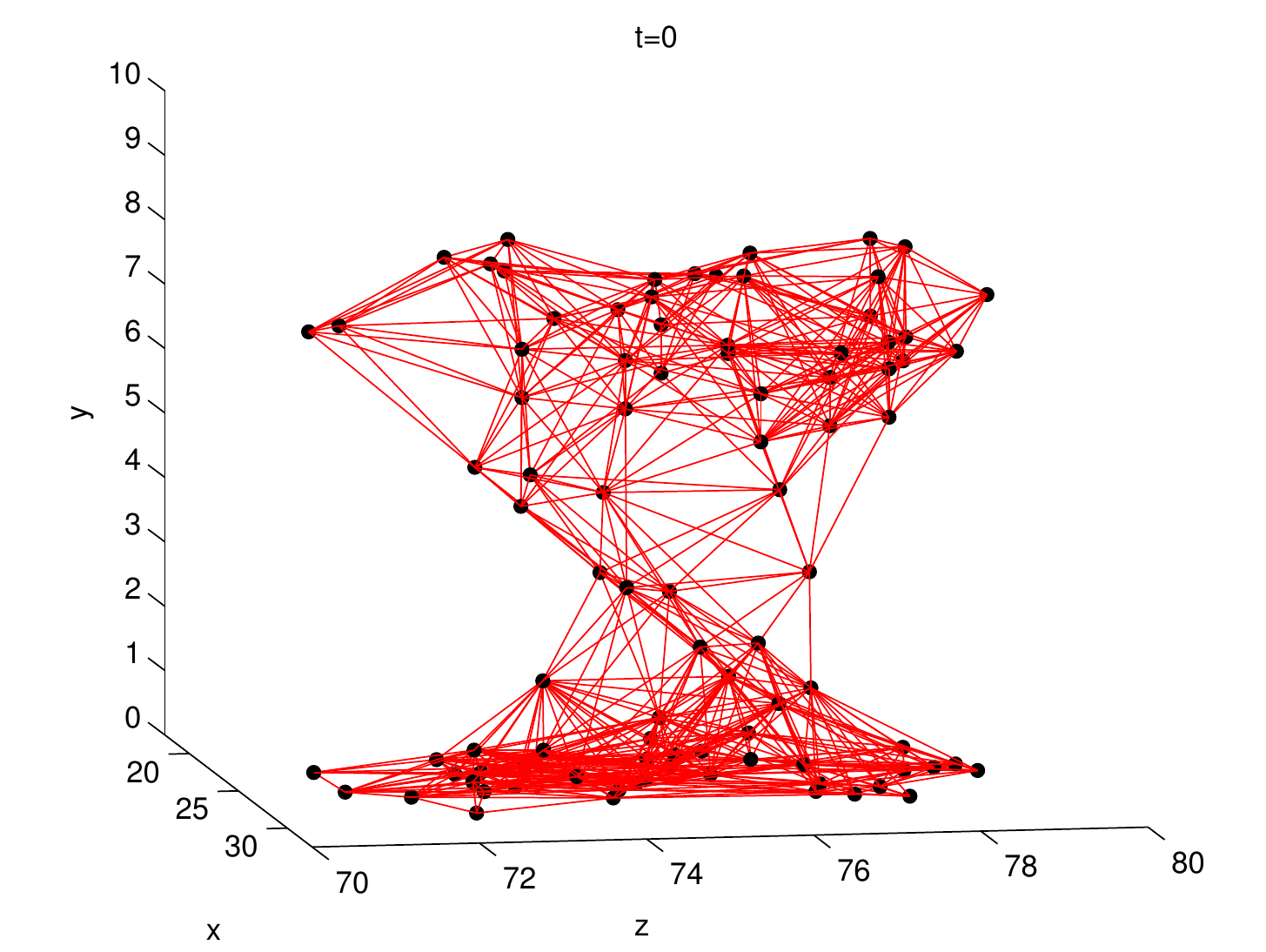}(b)\includegraphics[width=3in]{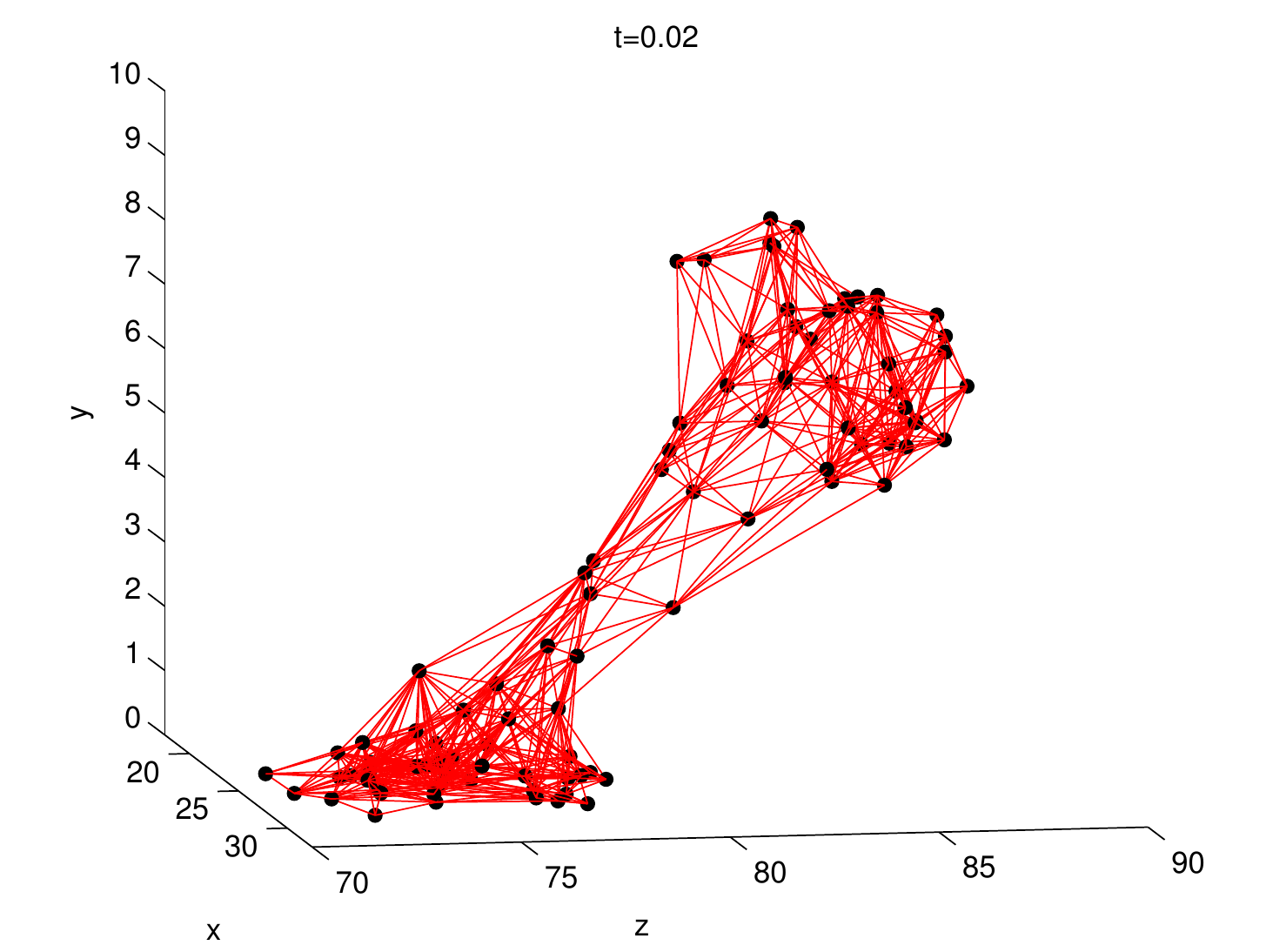}\\
(c)\includegraphics[width=3in]{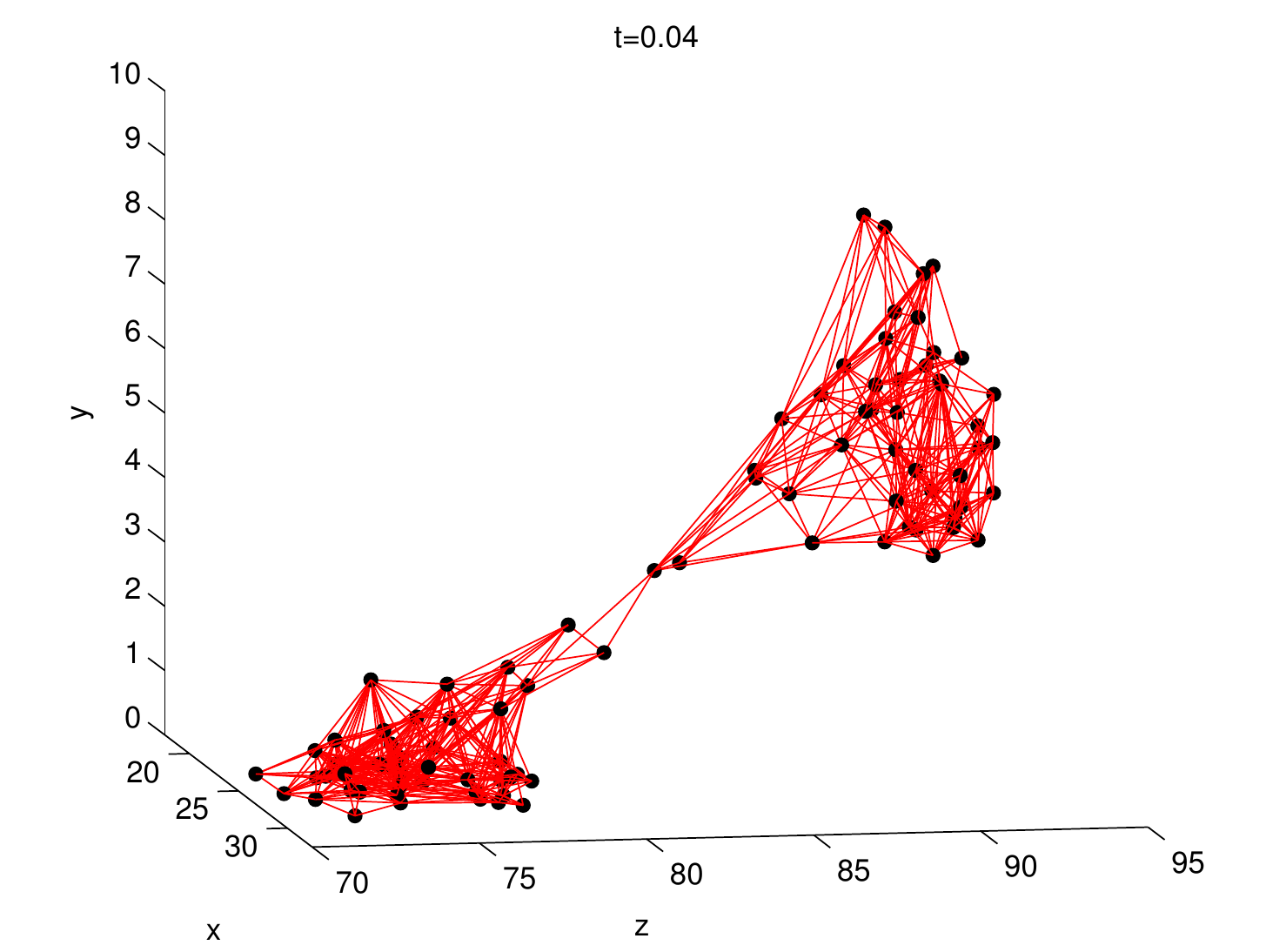}(d)\includegraphics[width=3in]{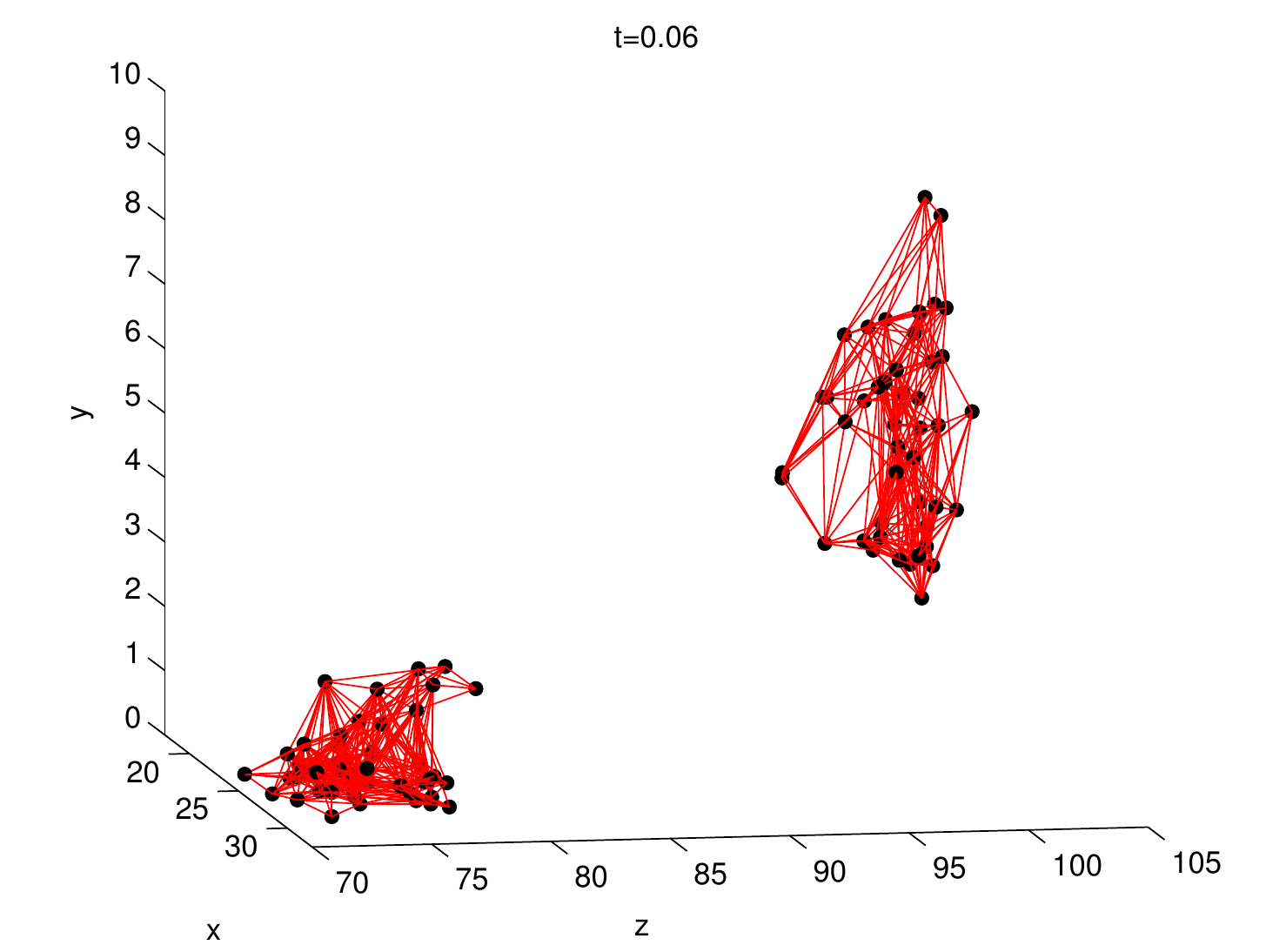}
\par\end{centering}

\caption{\label{fig:Shroom3DNoaddl_0_5}Full 3D simulation of a mushroom shaped
biofilm with the same density as the surrounding fluid. The time is
in seconds and the distance is in microns. As the biofilm stretches
in the flow, the strain in the midsection exceeds the breaking length
of the connections, and the top of the biofilm breaks off into the
flow. Then the broken piece tumbles end over end through the flow,
and the base retracts back. In this simulation, $\rho_{0}=998\,\nicefrac{\textrm{kg}}{\textrm{m}^{3}},\:\rho_{b}=0,\: F_{max}=1.25\times10^{-12}\,\textrm{N}$. }
\end{figure}

\begin{figure}[H]
\centering{}(a)\includegraphics[width=3in]{shroom3Dnondim_noAddl_t_0_04_f0_5}(b)\includegraphics[width=3in]{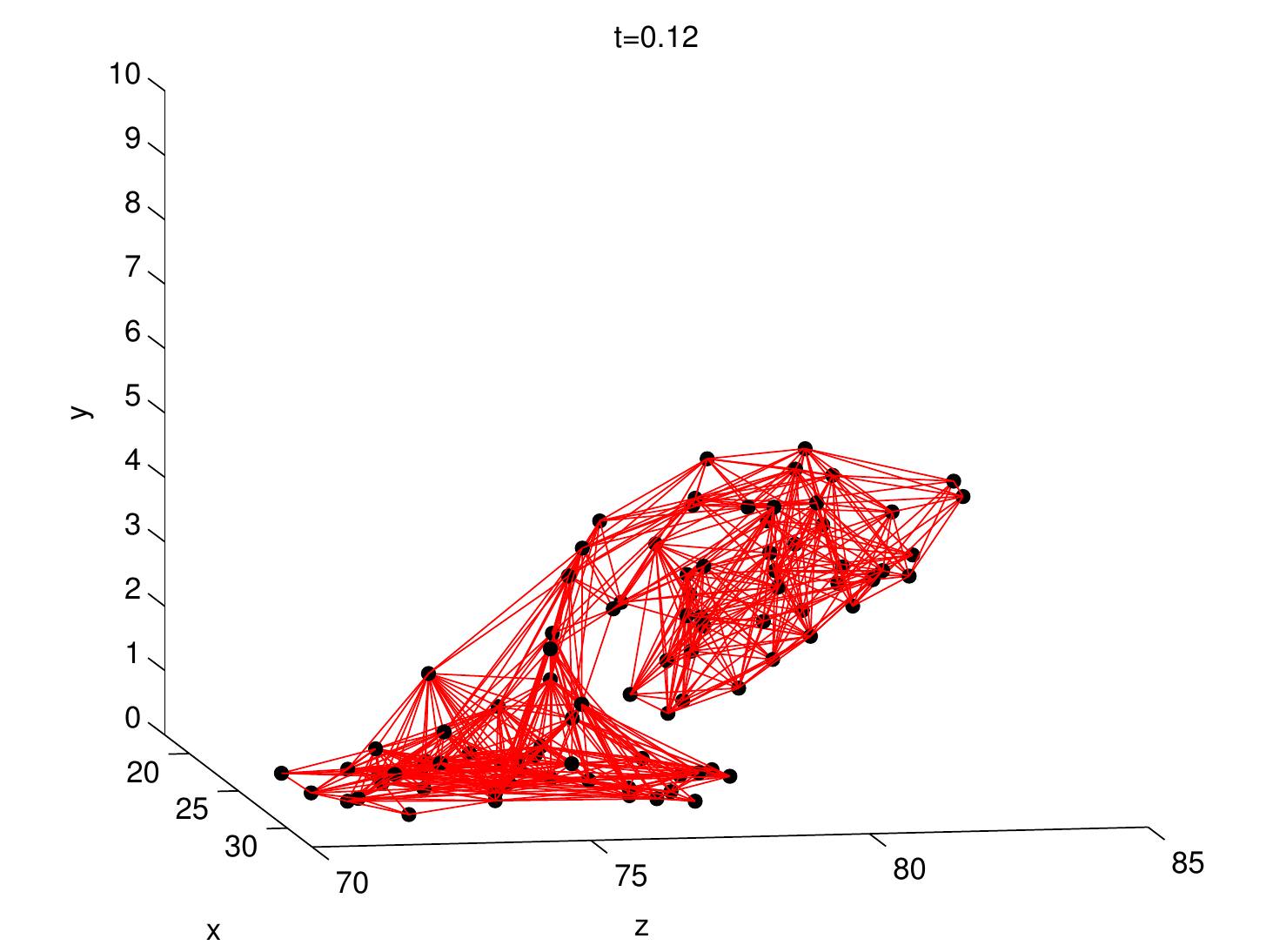}\\
(c)\includegraphics[width=3in]{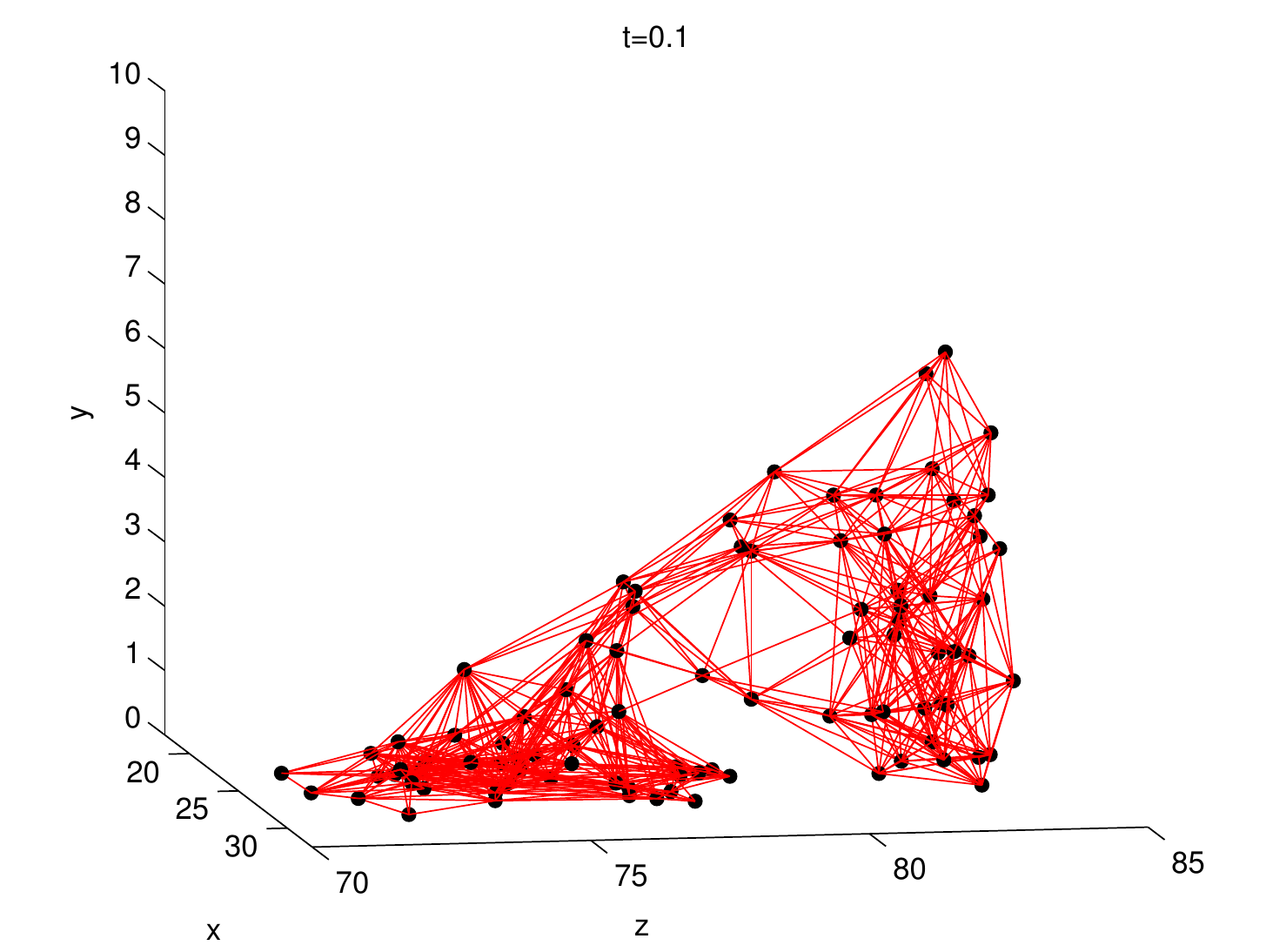}\includegraphics[width=3in]{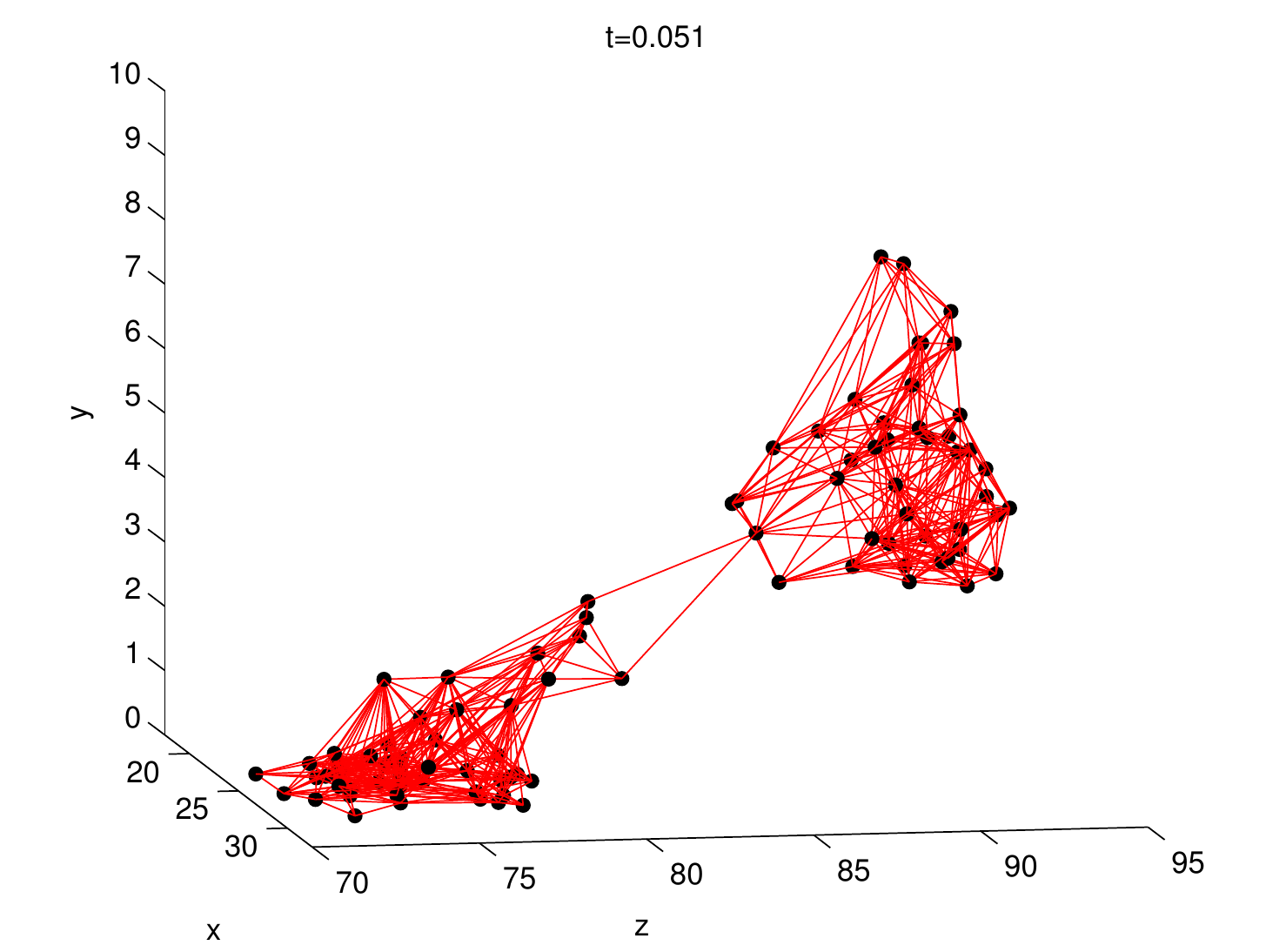}\caption{\label{fig:Shroom3Daddl}Snapshots of cell distributions resulting
from 3D simulation of mushroom shaped biofilm with four different
property configurations. In (a), the biofilm has the same density
and viscosity as the surrounding fluid (just before detachment), (b)
increased $F_{max}$, (c) density twice as large as the surrounding
fluid density, and (d) viscosity is $500\times$ that of the surrounding
fluid (just before detachment). }
\end{figure}

Our 3D simulations were run on a $50\times50\times150\,\mu\textrm{m}$
computational domain. We simulate on a mushroom shaped biofilm with
a height of about $8.5\,\mu\textrm{m}$ and a diameter of about $7.5\,\mu\textrm{m}$.
This shape is carved from the same set of data points described in
\prettyref{sub:SimSetup}. The spring connections between Lagrangian
nodes are put in place at the beginning of the simulation, with every
node connected to every other node less than $d_{c}=3\,\mu\textrm{m}$
away. Note that, for the 3D simulations, we increased $d_{c}$ slightly
to establish enough connections in the biofilm. We again use $\omega=\frac{1}{100}$
to match the radius of \textit{Staphylococcus epidermidis} and choose
$h=\frac{1}{128}$ in all of the simulations shown, so that $\omega>h$.
In the first simulation, the maximum spring force, $F_{max}$, is
set to $1.25\times10^{-12}\,\textrm{N}$. %
{} We again chose the value of these spring constants in order to illustrate
specific behaviors. The results of the first simulation are shown
in \prettyref{fig:Shroom3DNoaddl_0_5}. The mushroom shaped biofilm
bends over and stretches in the flow. The connections in the midsection
of the biofilm exceed their breaking length and the top of the biofilm
breaks off into the flow. Next, we ran a simulation of the same mushroom
shaped biofilm, but we added $\rho_{b}=998\,\nicefrac{\textrm{kg}}{\textrm{m}^{3}}$
additional density to the biofilm compared to the surrounding fluid.
The final result is provided in \prettyref{fig:Shroom3Daddl}(b) and
illustrates that the added density increases the momentum of the biofilm.
This allows for the mushroom to curl over into the flow and increases
the time until detachment. We also ran a simulation of the same mushroom
shaped biofilm, but we increased $F_{max}$ to $1\times10^{-11}\,\textrm{N}$
and kept the biofilm density the same as the surrounding fluid. The
result is provided in \prettyref{fig:Shroom3Daddl}(c). The effect
of these stronger springs is that the thin part of the biofilm does
not stretch enough to break the connections. We can see from these
simulations that either increasing the biofilm density or strengthening
the springs causes similar results, but with the increased density
the biofilm just curls over.

Finally, we provide one 3D simulation to show that, with increased
biofilm viscosity, they produce qualitatively the same behavior as
in the 2D case. In the simulation result shown in \prettyref{fig:Shroom3Daddl}(d),
we use the same parameters as the first 3D simulation, but we use
a viscosity in the biofilm that is a factor of $500$ times that of
the surrounding fluid. Just as in the 2D case, this results in a longer
time until detachment (compare time in \prettyref{fig:Shroom3Daddl}(a)
and \prettyref{fig:Shroom3Daddl}(d)). 

For ease of comparison, we now provide the figures in the order in
which they were discussed in this section.

\section{\label{sec:Realistically-Shaped-Biofilm}Realistically Shaped Biofilm
Simulation}

The biofilm shapes used in \prettyref{sub:Two-Dimensional-Simulations}
and \prettyref{sub:Three-Dimensional-Simulations} were intentionally
carved from the data in a way to provide a weak point at which it
would be most likely to break. This was done in order to illustrate
the effects of varying the different parameters in the simulation.
In this section, we provide results of the simulation on a biofilm
that is a subset of points taken directly from the real biofilm data
set. In reality, this biofilm was surrounded by more cells on all
sides, which would change the behavior of the fluid structure interactions.
However, we use this to show the results of the simulation on a \textit{real}
top heavy biofilm shape that was grown in a lab. The \textit{Staphylococcus
epidermidis} data set discussed in \prettyref{sub:SimSetup} was supplied
as positions in three $30\times30\times15\,\mu\textrm{m}$ sub-domains
of a biofilm.

\begin{figure}[H]
\begin{centering}
(a)\includegraphics[width=3in]{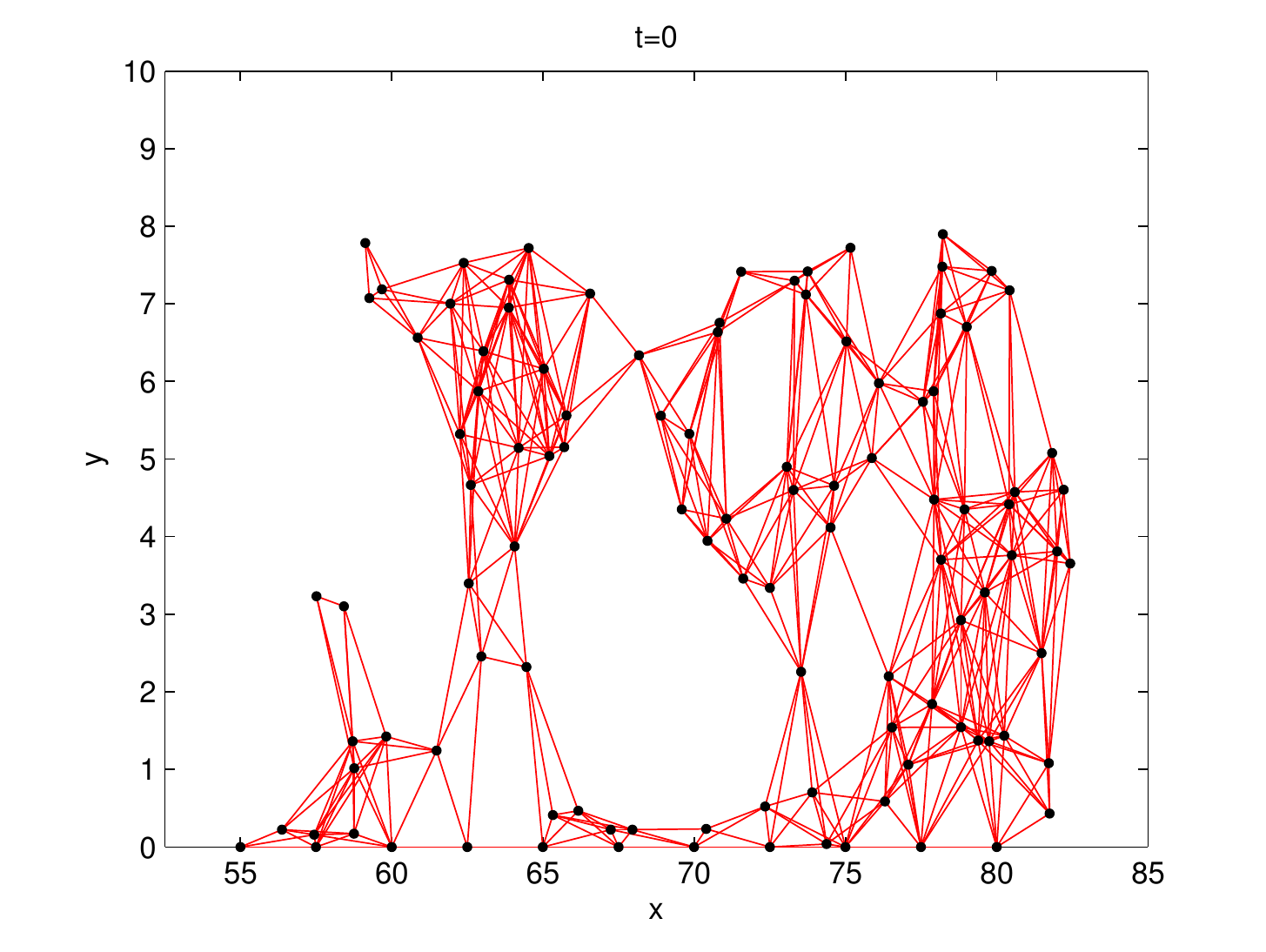}(b)\includegraphics[width=3in]{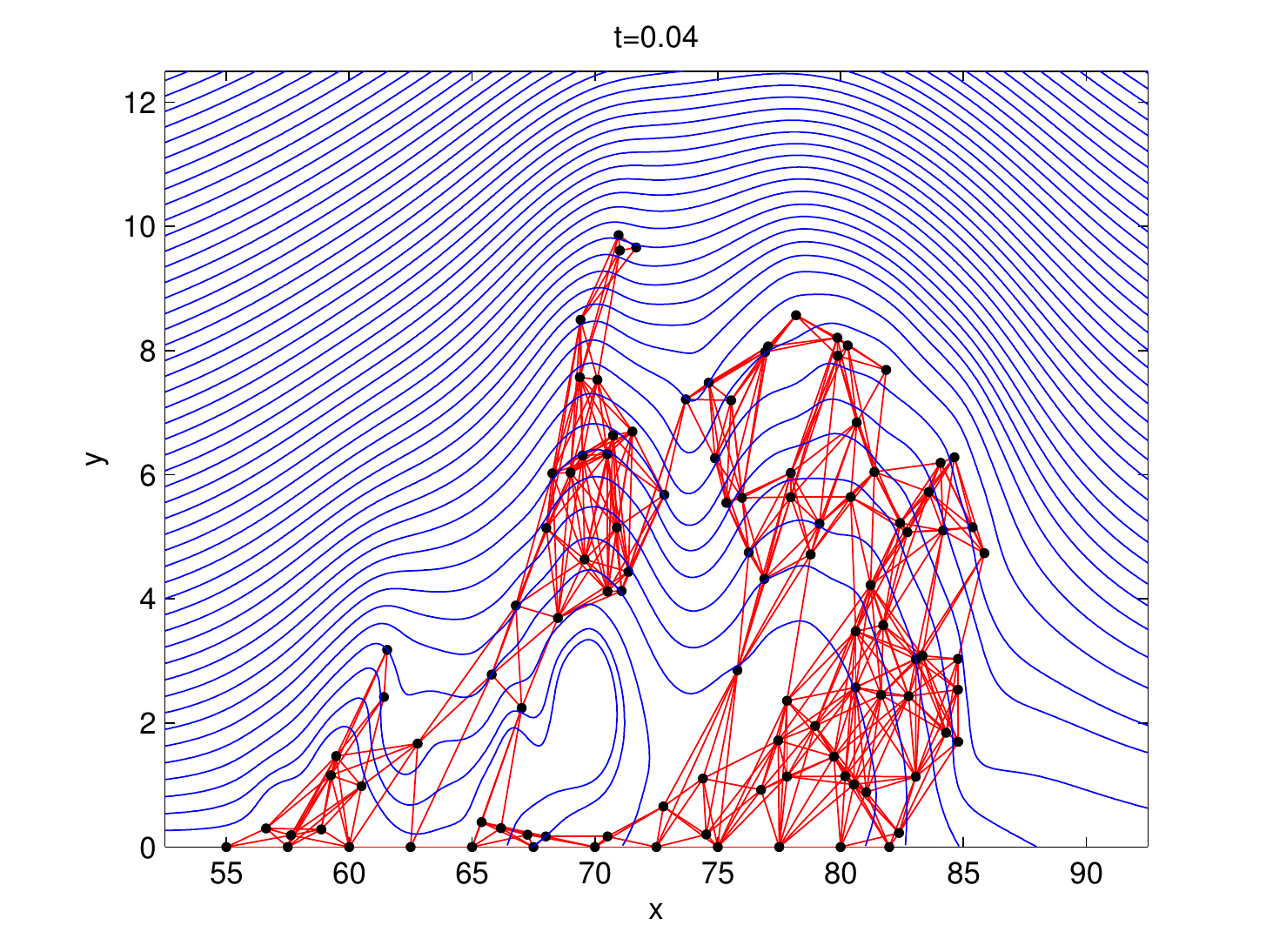}\\
(c)\includegraphics[width=3in]{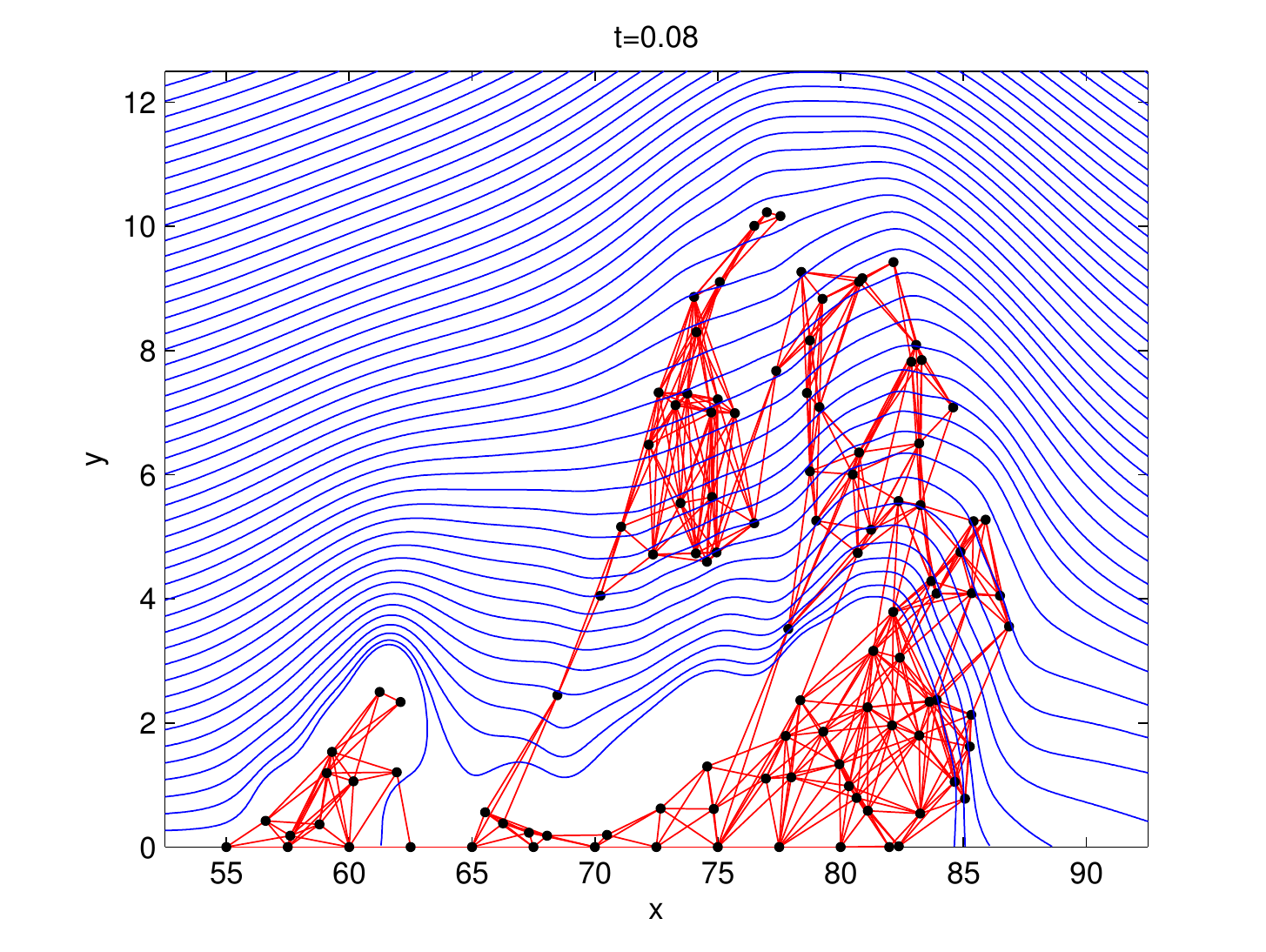}(d)\includegraphics[width=3in]{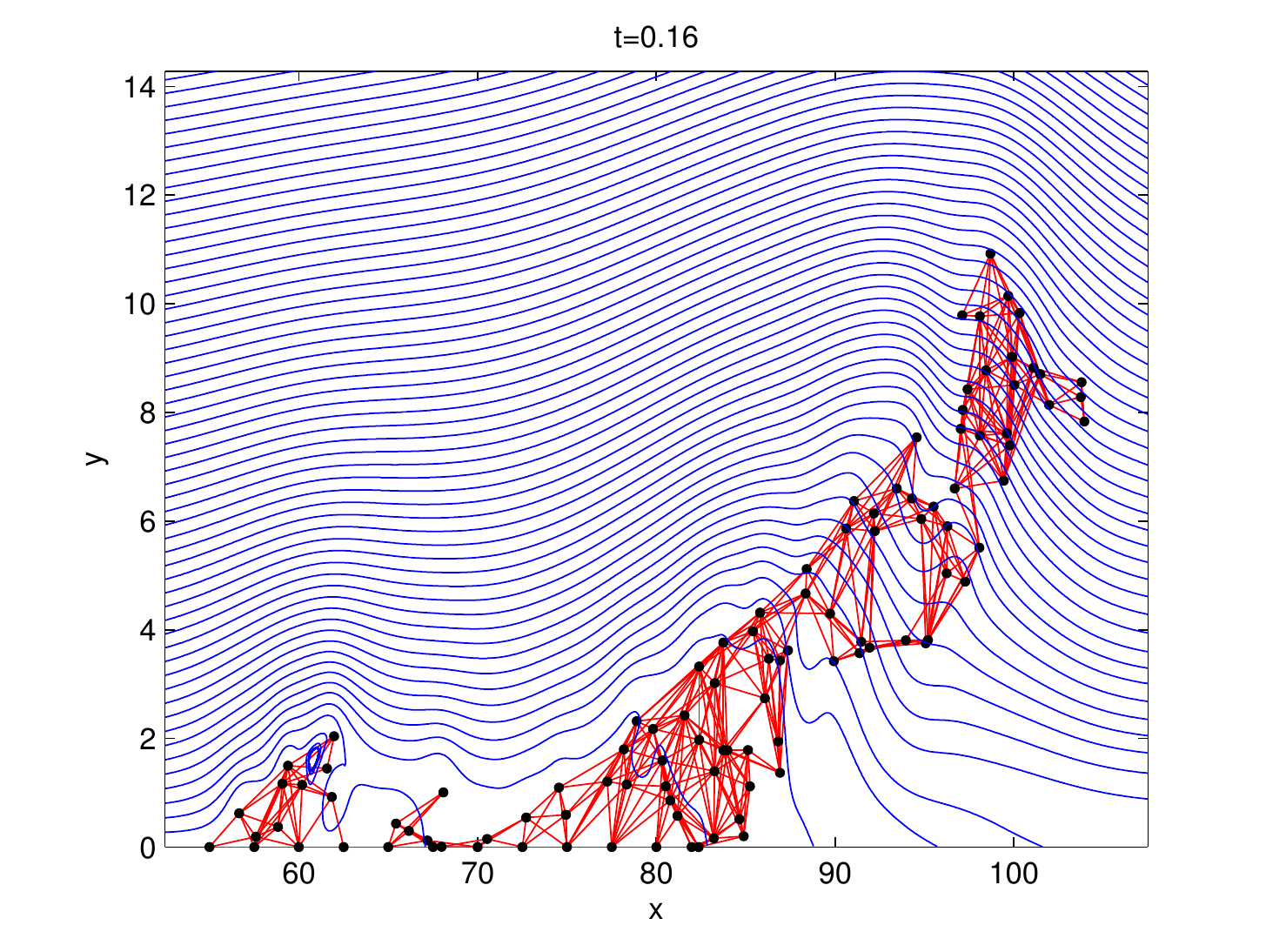}\\
(e)\includegraphics[width=3in]{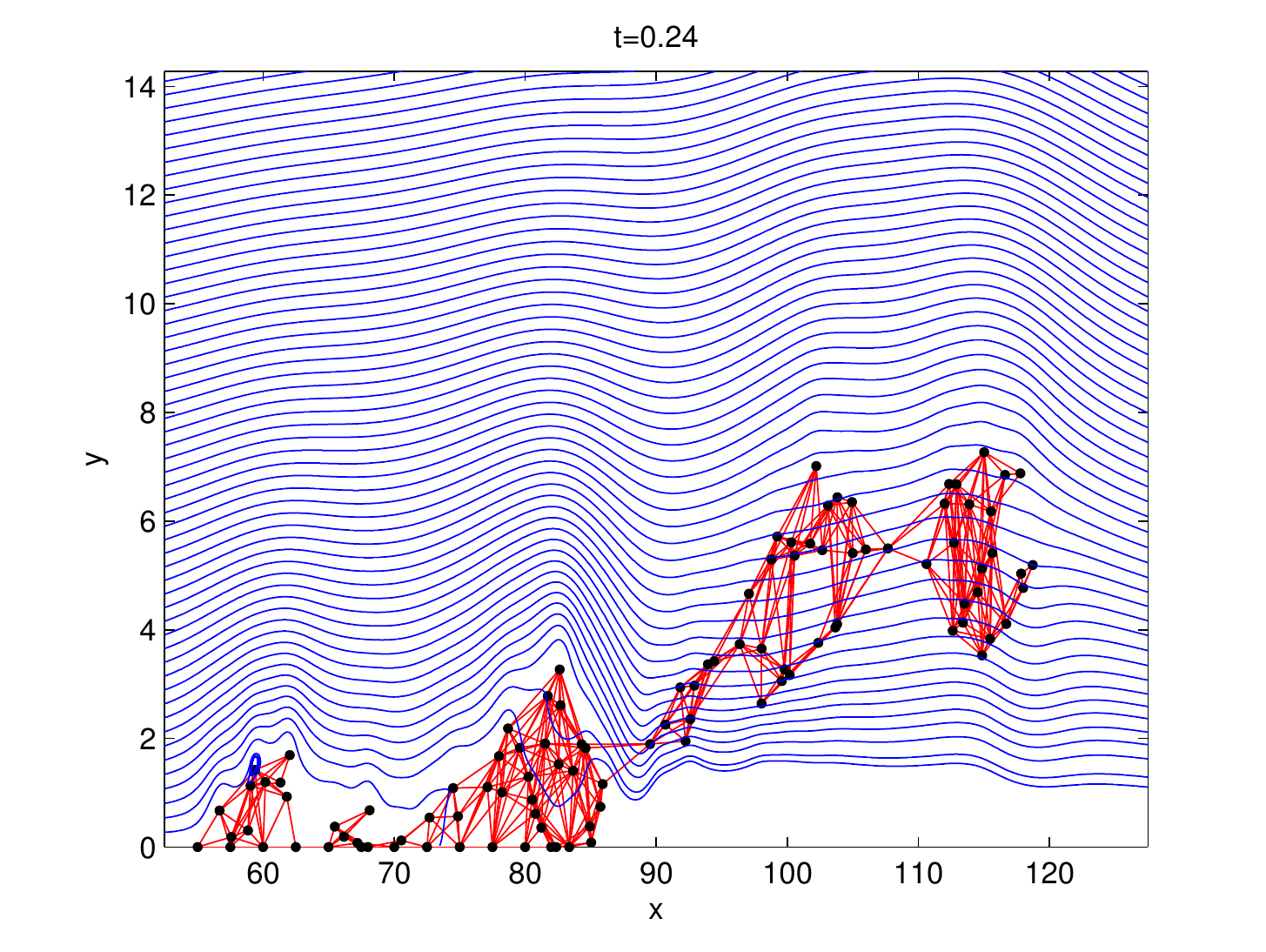}(f)\includegraphics[width=3in]{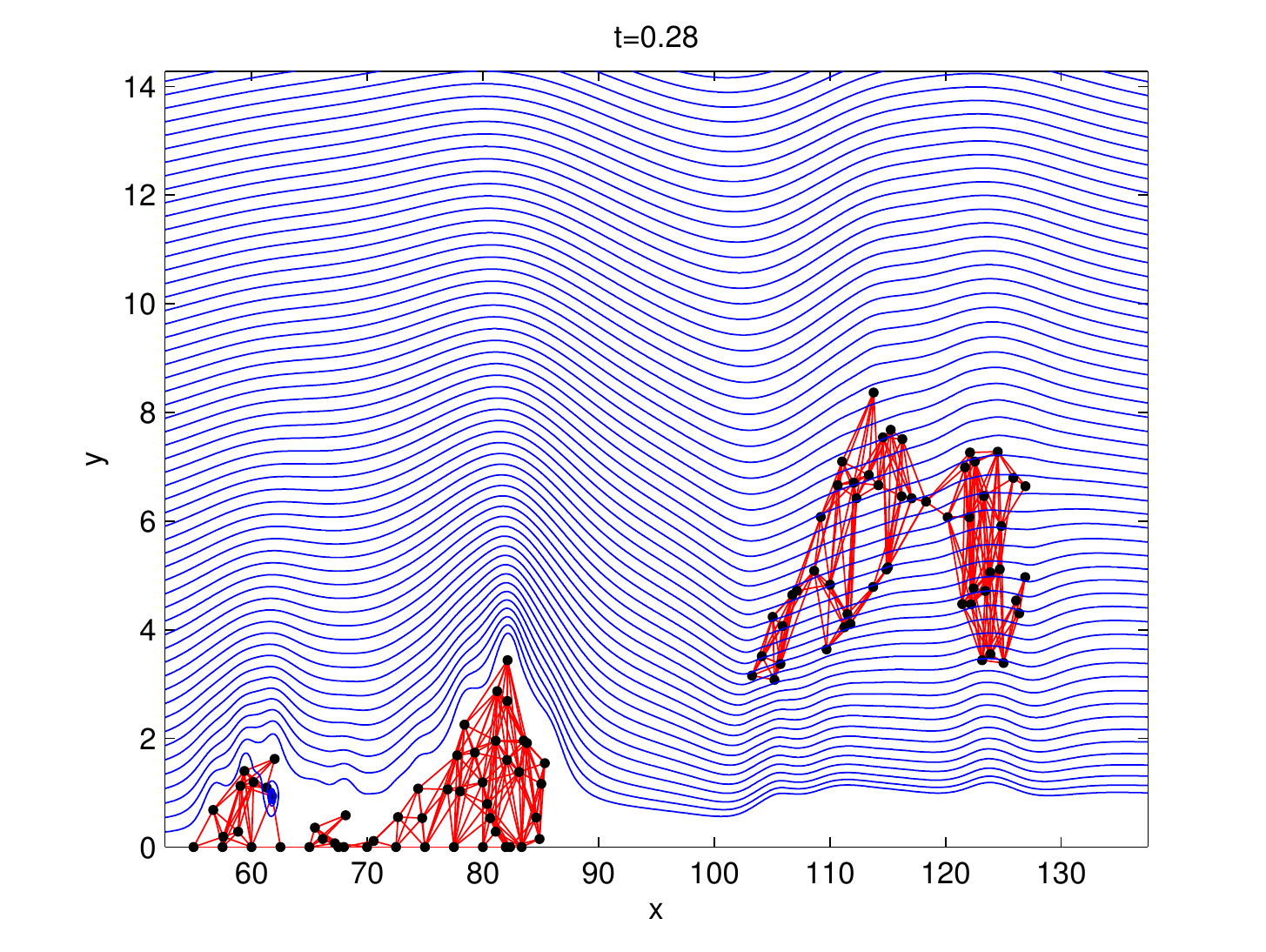}
\par\end{centering}

\caption{\label{fig:RealShroom2Dnoaddl}Simulation on a 2D slice of a real
biofilm with the same density as the surrounding fluid. Time is in
seconds and the distance is in microns. In this simulation, $\rho_{0}=998\,\nicefrac{\textrm{kg}}{\textrm{m}^{3}},\:\rho_{b}=0,\: F_{max}=7.5\times10^{-7}\,\textrm{N}$.
The streamlines follow the velocity field. In this simulation, the
mushroom shaped part pushes against the biofilm behind it (b), then
rolls over the top of it as it breaks from its base (d). Then a large
portion of the biofilm breaks completely off leaving 2 distinct bases
(f).}
\end{figure}

In \prettyref{fig:RealShroom2Dnoaddl}(a), we show the biofilm taken
from a $2\times30\times15\mu\textrm{m}$ subset of one data set that
has been connected with $d_{c}=2.8\,\mu\textrm{m}$. For the 2D representation,
we collapse the $2\,\mu\textrm{m}$ dimension, leaving only the $(x,y)$
coordinates of the data. The most interesting feature of this biofilm
is that in the region from $x=60\,\mu\textrm{m}$ to $x=67\,\mu\textrm{m}$
the biofilm exhibits a mushroom shape similar to the one we used in
\prettyref{sub:Two-Dimensional-Simulations}.

\begin{figure}[H]
\begin{centering}
(a)\includegraphics[width=3in]{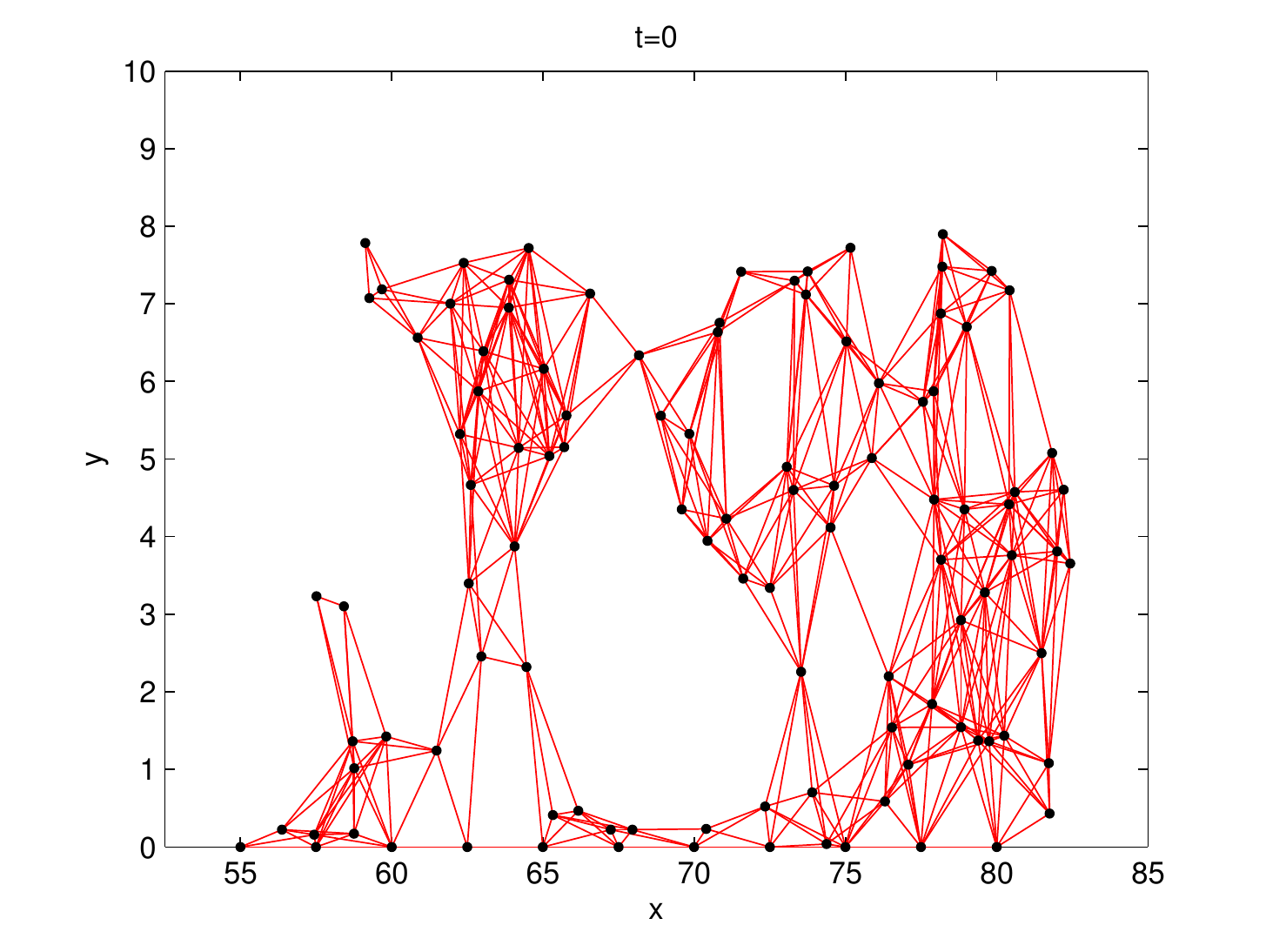}(b)\includegraphics[width=3in]{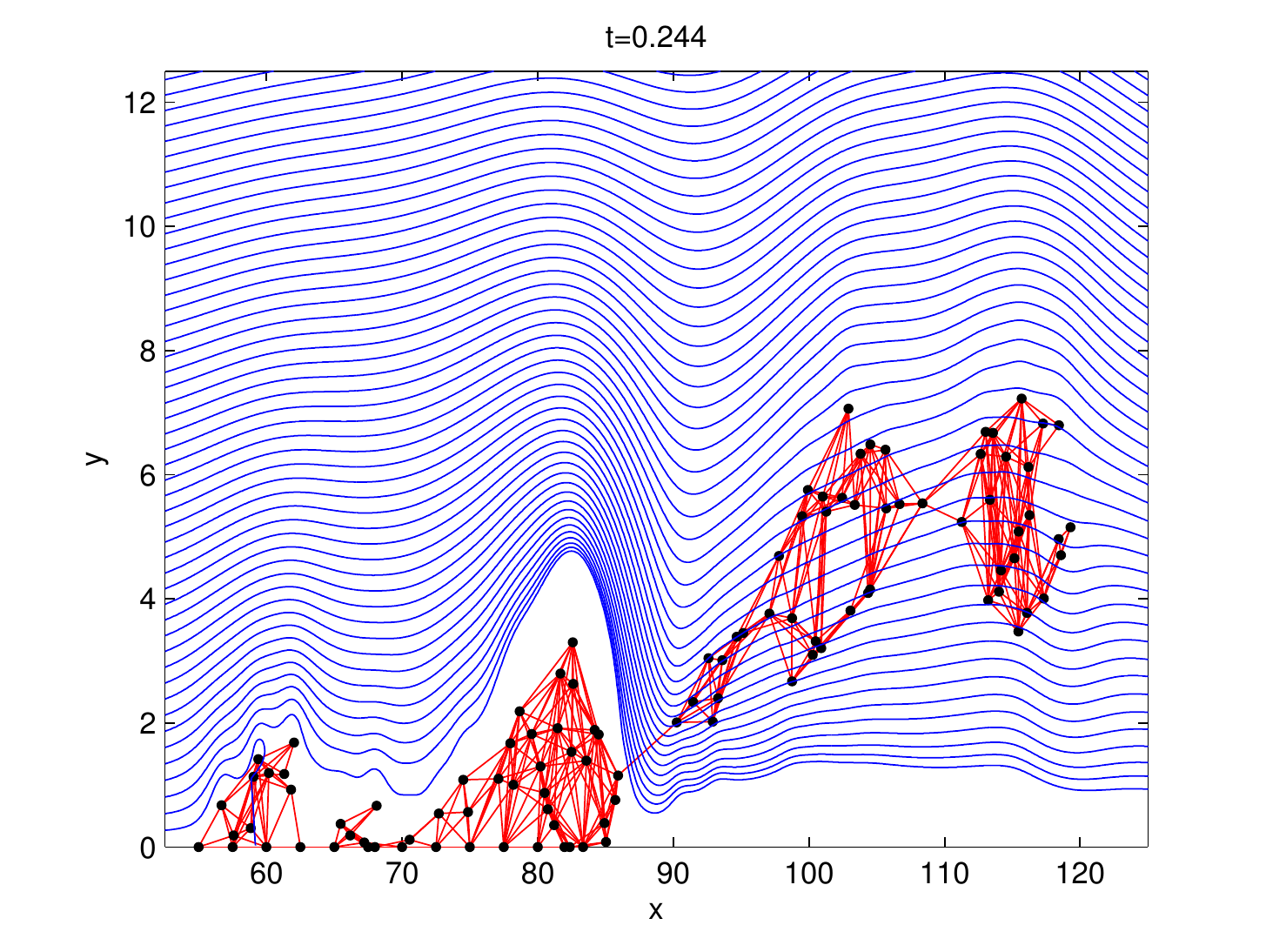}\\
(c)\includegraphics[width=3in]{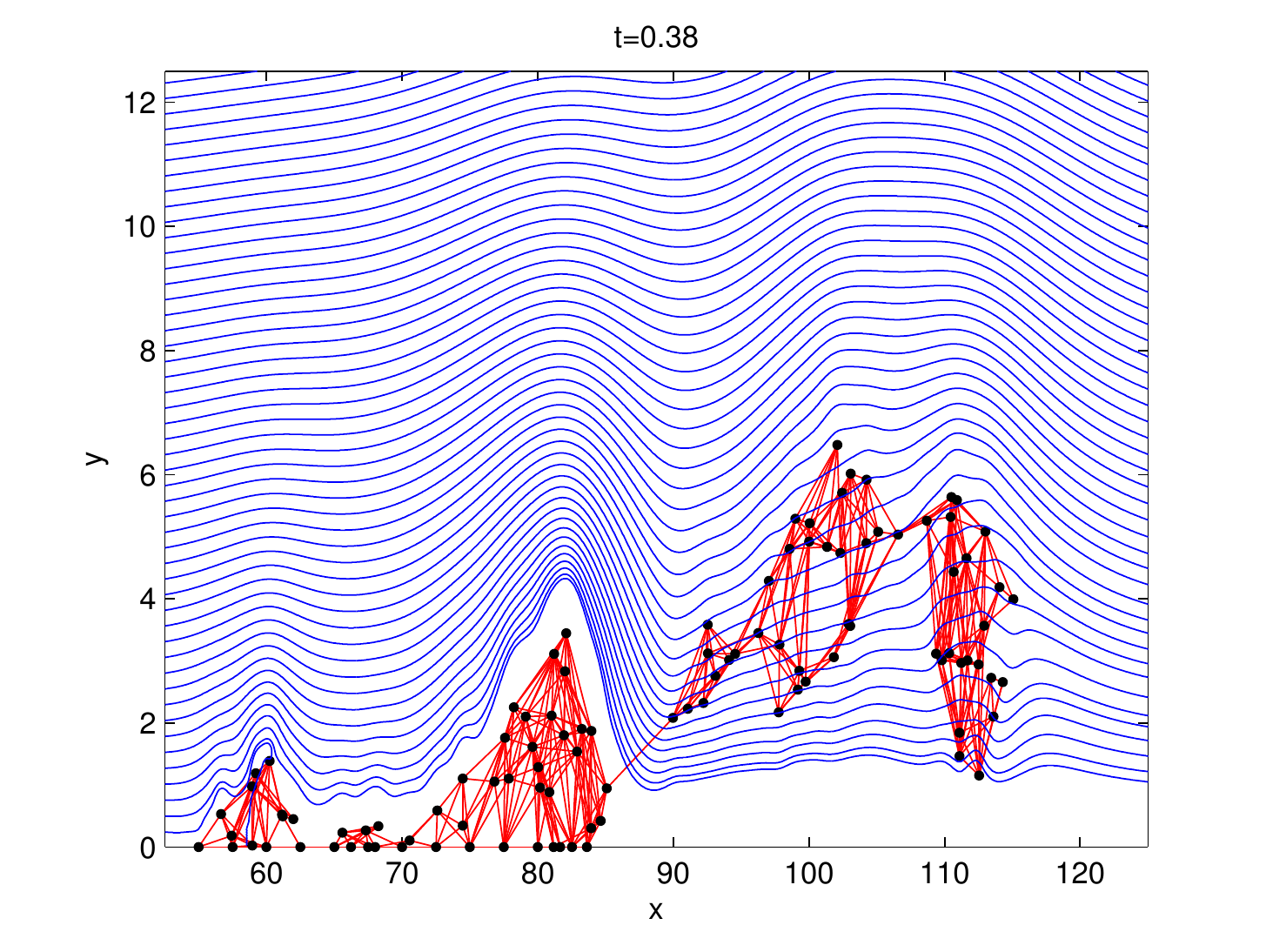}(d)\includegraphics[width=3in]{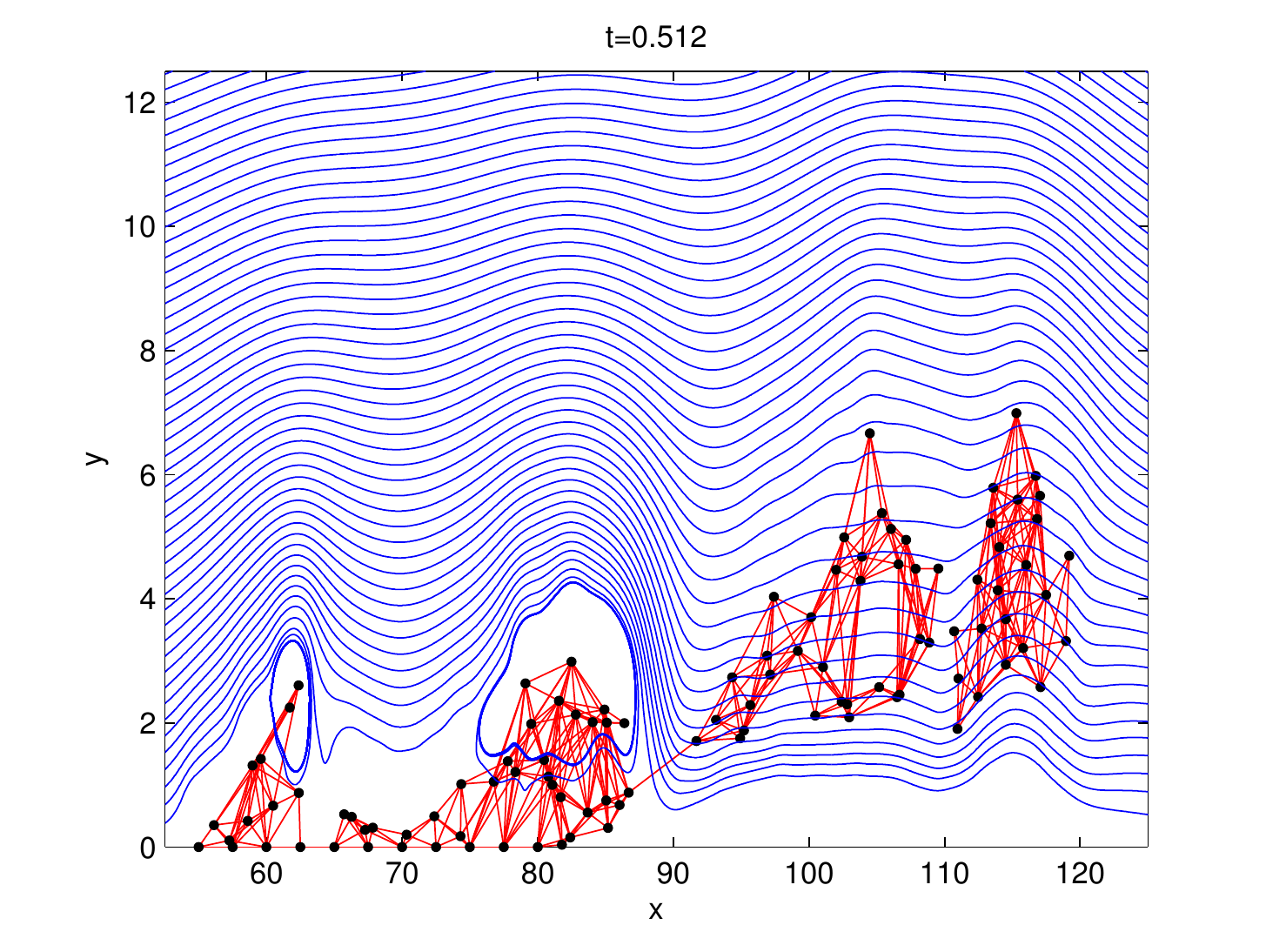}
\par\end{centering}

\caption{\label{fig:RealShroom2Daddl}Snapshots in time showing the detachment
time and configuration at detachment of a 2D slice of a real biofilm
with initial configuration in (a). In (b), biofilm has the same density
and viscosity as the surrounding fluid, (c) 12\% larger density, and
(d) $50\times$ larger viscosity than the surrounding fluid. The streamlines
follow the velocity field. }
\end{figure}

We now provide two simulation results on this realistically shaped
biofilm. The first simulation (see \prettyref{fig:RealShroom2Dnoaddl})
uses a biofilm density equal to the surrounding fluid and uses $F_{max}=7.5\times10^{-7}\,\textrm{N}$.
In this simulation, the mushroom shaped part pushes against the biofilm
behind it, then rolls over the top of it as it breaks from its base,
forming a long streamer-like biofilm.%
\footnote{Streamers are a natural occurrence in biofilms. Examples in \citep{Klapper2002,Rusconi2010}.%
}Then the streamer breaks completely off leaving two distinct attached
structures. In the second simulation, we use a biofilm density of
$\rho_{b}=120\,\nicefrac{\textrm{kg}}{\textrm{m}^{3}}$ and kept everything
else the same. Although the density is only 12\% larger than the surrounding
fluid, it has a large impact on the outcome of the simulation. In
this simulation, the effect of the increased density of the biofilm
is a longer breaking time (compare time in \prettyref{fig:RealShroom2Daddl}(b)
and \prettyref{fig:RealShroom2Daddl}(c)). This occurs since the increased
momentum causes the first detached piece to continue further down,
pulling the whole streamer lower (compare the height of the detaching
pieces). The fluid forces continue to push the streamer until it breaks
into the flow.

In the final simulation, we show that increasing the viscosity in
the ``realistically'' shaped biofilm has larger impact on the results
than in the case of the previous standalone mushroom shaped biofilm.
In this simulation, we increased the biofilm viscosity to 50$\times$
the surrounding fluid with $\mu_{max}=0.05\,\nicefrac{\textrm{kg}}{\textrm{m}\cdot\textrm{s}}$.
Although this is 10$\times$ less than in the previous variable viscosity
simulation, it has a larger impact on this wider biofilm, doubling
the detachment time from the case of constant viscosity (compare \prettyref{fig:RealShroom2Daddl}(b)
and \prettyref{fig:RealShroom2Daddl}(d)). 

\subsection{\label{sub:Numerical-Concerns}Numerical Concerns}

We note that there remains one pressing concern that will be the focus
of our future work. When adapting the multigrid scheme for large values
of $\mu_{max}$, we do not achieve expected speed ups in convergence
rates. We use restriction to transfer the viscosity to the coarse
grids works for small values of $\mu_{max}$, but we found that, for
larger values of $\mu_{max}$, this technique leads to a very slowly
converging solver. Intriguingly, we found that using our restriction
operator to define the coarse grid viscous values and then scaling
the values leads to faster convergence. Specifically, we define the
coarse grid viscosity as
\begin{eqnarray}
\mu_{lh}(\mathbf{x}) & = & \gamma_{lh}I_{\frac{l}{2}h}^{lh}\mu_{\frac{l}{2}h}(\mathbf{x}),\label{eq:ViscCoarseGrid}
\end{eqnarray}
where $lh$ denotes the grid whose mesh width is $l$ times $h$ ($l=2,\,4,\,8,\,16,\ldots$),
$\gamma_{lh}\in(0,1]$ is the scaling which maximizes the convergence
of the solver, and $\mu_{h}(\mathbf{x})$ is defined by \prettyref{eq:ViscExp}.
Through repeated experimentation, we found that using $\gamma_{lh}\in[.7,1]$
resulted in the fastest convergence rates in the 2D and 3D simulations.
This approach is admittedly ad-hoc. However, the consistency with
which we achieved dramatic speed ups strongly suggests the existence
of an underlying mathematical principle to be discovered. 

For our future work, the highest priority is to resolve the problem
of slow convergence for large $\mu_{max}$. There are three approaches
that may lead to resolving this issue. The first would be to mathematically
derive optimal values for the scaling parameters, $\gamma_{lh}$.
Second, we could ensure that our discretization satisfies the Galerkin
condition. Lastly, and probably the best choice, would be to re-implement
the geometric multigrid as an algebraic multigrid method (AMG, see
Ch. 8 of \citep{Briggs2000}).

\section{\label{sec:Biofilm-Simulation-Conclusions}Conclusions}

In this work we developed a simulation to model the flow-induced fragmentation
of biofilms. In this simulation, we have provided a way to adjust
the biofilm density and viscosity, which had not been addressed in
previous IBM biofilm models. We also have control of the fluid flow
rate, density, viscosity, and elastic forces within the biofilm. We
used experimentally measured biofilm bacterial cell locations as initial
positions for our Lagrangian nodes. This is dramatically different
than the traditional IBM, in which methods usually refine the Lagrangian
mesh along with the Eulerian mesh. We adapted the Dirac delta approximation
to scale with the radius of the bacteria rather than with the mesh
width. This implies that the information that transfers from the Lagrangian
grid to the Eulerian grid (i.e.,~density, viscosity, and elastic
force) is spread over a set distance rather than scaling by the mesh
width, $h$. This adapted Dirac Delta approximation improves our numerical
convergence rates as well.

We used a projection method to split the incompressible Navier-Stokes
equations to solve separately for an intermediate velocity and the
pressure, and then used a Gauss-Seidel iterative method with multigrid
to solve the resulting equations. Using an iterative solver, as opposed
to a spectral method, to solve these systems was necessitated by the
fact that biofilms have spatially varying density and viscosity. With
this solver we achieved first order convergence in both space and
time. 

For the numerical simulations, we carved a mushroom shaped biofilm
from the bacterial cell locations and ran simulations with varying
parameters. We first ran the simulation on a simplistic shape in order
to validate the effect of the various parameter changes on the outcome
of the biofilm. By adjusting the maximum elastic force, $F_{max}$
in the biofilm, we controlled the detachment phenomenon. We also showed
that slight changes in the density of the biofilm has a large effect
on the outcome of the simulation. This is an important conclusion
as usually modelers ignore the differences in biofilm density. Finally,
we showed that we can increase the detachment time in the simulations
by increasing the viscosity of the biofilm. Finally, we ran simulations
on more realistically shaped biofilms, which showed how a larger biofilm
with different shapes will react to fluid flow forces. Adjusting these
parameters will be a necessary component when we attempt to match
these simulations to experimental data.

\section{\label{sec:FW} Future Work}

There are several directions in which we plan to take this research
in the future. For example, there are straightforward ways to include
more biologically realistic terms to interpret biofilm internal stress
dynamics, cell volume, and fragmentation dynamics. Additionally, there
are approaches that may greatly improve our numerical scheme convergence
and stability, including alternative multigrid algorithms, implicit
discretizations, and improved immersed boundary implementations.

In its current form, our simulations could be used to make predictions
in detachment times of biofilms, as well as general behavioral responses
of biofilms to various flow conditions. We plan to work closely with
experimentalists to formulate accurate viscoelastic models for the
biofilms and modify our constitutive equations for stress and elasticity
to account for these model choices. We also plan to include the fact
that bacterial cells displace fluid. While we have adapted the Dirac
delta function approximation to transfer the cell parameters ($F$,
$\rho$, $\mu$) to the Eulerian grid, the current simulation does
not actually assign a size to the cells. As a result, the cells are
free to pass through each other. We first plan to alter the model
for the bacterial cell so that it displaces fluid. An important step
in this process will be to identify a collision detection strategy.
We will base ours on potentials for electrostatic, steric, and Van
der Waals forces.

Our current simulation uses a spring-breaking criteria of double the
rest length, and we also assume that the bonds are linearly elastic
until the breaking point. This is not an accurate assumption, as it
is known that biofilms are composed of polymer based ECMs. These structures
are linearly elastic for small strains and then experience plastic
deformation (permanently altering the bonds in the ECM and thus the
rest length) before finally fracturing. In the future, we will use
biofilm yielding data from experiments such as \citep{Aggarwal2010}
to determine accurate approximations for yield points and fracture
points in the biofilm. We plan to include plasticity into the simulations
by changing the equations for stress, \prettyref{eq:springF}, when
the bond has been stretched beyond its yield point.

We have several plans for improving our numerical method. Our current
simulation is limited to first-order accuracy. Guided by the results
in \citep{BrownCortezMinion2001}, we will accurately derive the numerical
boundary conditions for our projection method to ensure second order
accuracy for both velocity and pressure computations. To improve the
accuracy of the immersed boundary method, we could also adapt our
modeling method to either an immersed interface method (\citep{LiIto2006})
in which we adapt the finite difference approximations close to the
interface or a blob projection immersed boundary method as discussed
in \citep{Cortez2000} in order to obtain second-order spatial accuracy.
Another limitation of our current numerical scheme is the time-step
stability restrictions, which limit the size of the elastic forces
between the cells. We plan to eliminate these restrictions altogether
by changing to a semi-implicit or implicit method of transferring
the data between the Eulerian and Lagrangian grids, as is shown by
Newren, Fogelson et al., in \citep{Newren2007}. 

Finally, with large biofilm densities and viscosities, our multigrid
method in its current formulation does not converge as fast as expected.
We plan to fix this by appropriately adapting our implementation of
the geometric multigrid (by satisfying the Galerkin condition) or
by changing to an algebraic multigrid approach.

\section{Acknowledgements}

We thank Dr. Stephen McCormick for the discussions we had about using
multigrid in our simulations.

This work was supported in part by the National Science Foundation
grants PHY-0940991 and PHY-0941227. 

This work utilized the Janus supercomputer, which is supported by
the National Science Foundation (award number CNS-0821794), the University
of Colorado Boulder, the University of Colorado Denver, and the National
Center for Atmospheric Research. The Janus supercomputer is operated
by the University of Colorado Boulder.

\bibliographystyle{abbrv}
\bibliography{mathbioCU}

\begin{thebibliography}{10}

\bibitem{Aggarwal2010}
S.~Aggarwal, E.~H. Poppele, and R.~M. Hozalski.
\newblock {Development and testing of a novel microcantilever technique for
  measuring the cohesive strength of intact biofilms.}
\newblock {\em Biotechnology and bioengineering}, 105(5):924--34, Apr. 2010.

\bibitem{Alpkvist2007}
E.~Alpkvist and I.~Klapper.
\newblock {Description of Mechanical Response Including Detachment Using a
  Novel Particle Method of Biofilm/Flow Interaction}.
\newblock {\em Water Science and Technology}, 55(8-9):265--273, May 2007.

\bibitem{Aravas2008}
N.~Aravas and C.~S. Laspidou.
\newblock {On the calculation of the elastic modulus of a biofilm streamer.}
\newblock {\em Biotechnology and bioengineering}, 101(1):196--200, Sept. 2008.

\bibitem{Atkinson1989}
K.~Atkinson.
\newblock {\em {An Introduction to Numerical Analysis}}.
\newblock Wiley, New York, NY, 2nd edition, 1989.

\bibitem{Balestrino2008}
D.~Balestrino, J.-M. Ghigo, N.~Charbonnel, J.~A.~J. Haagensen, and
  C.~Forestier.
\newblock {The characterization of functions involved in the establishment and
  maturation of Klebsiella pneumoniae in vitro biofilm reveals dual roles for
  surface exopolysaccharides.}
\newblock {\em Environmental microbiology}, 10(3):685--701, Mar. 2008.

\bibitem{Bottino1998}
D.~C. Bottino.
\newblock {Modeling Viscoelastic networks and cell deformation in the context
  of the immersed boundary method}.
\newblock {\em J. Computational Physics}, 147:86--113, 1998.

\bibitem{Briggs2000}
W.~L. Briggs, V.~E. Henson, and S.~F. McCormick.
\newblock {\em {A Multigrid Tutorial}}.
\newblock SIAM, Philadelphia, PA, 2000.

\bibitem{BrownCortezMinion2001}
D.~L. Brown, R.~Cortez, and M.~L. Minion.
\newblock {Accurate Projection Methods for the Incompressible Navier-Stokes
  Equations}.
\newblock {\em Journal of Computational Physics}, 168(2):464--499, Apr. 2001.

\bibitem{chen1998directBiofilmStrengthTube}
M.~J. Chen, Z.~Zhang, and T.~R. Bott.
\newblock {Direct measurement of the adhesive strength of biofilms in pipes by
  micromanipulation}.
\newblock {\em Biotechnology Techniques}, 12(12):875--880, 1998.

\bibitem{Chen2005effectsAdhesiveStrength}
M.~J. Chen, Z.~Zhang, and T.~R. Bott.
\newblock {Effects of operating conditions on the adhesive strength of {{\it
  Pseudomonas fluorescens}} biofilms in tubes.}
\newblock {\em Colloids and surfaces. B, Biointerfaces}, 43(2):61--71, June
  2005.

\bibitem{Cortez2000}
R.~Cortez and M.~Minion.
\newblock {The Blob Projection Method for Immersed Boundary Problems}.
\newblock {\em Journal of Computational Physics}, 161(2):428--453, July 2000.

\bibitem{Dillon1996}
R.~Dillon, L.~Fauci, A.~Fogelson, and D.~{Gaver III}.
\newblock {Modeling Biofilm Processes Using the Immersed Boundary Method}.
\newblock {\em Journal of Computational Physics}, 129(1):57--73, Nov. 1996.

\bibitem{Ferziger2002}
J.~H. Ferziger and M.~Peric.
\newblock {\em {Computational Methods for Fluid Dynamics}}.
\newblock Springer, 3rd edition, 2002.

\bibitem{Huang2009IBM}
W.~X. Huang and H.~J. Sung.
\newblock {An immersed boundary method for fluid-flexible structure
  interaction}.
\newblock {\em Computer Methods in Applied Mechanics and Engineering},
  198(33-36):2650--2661, 2009.

\bibitem{Kissel1984}
J.~C. Kissel, P.~L. McCarty, and R.~L. Street.
\newblock {Numerical Simulation of Mixed-Culture Biofilm}.
\newblock {\em Journal of Environmental Engineering}, 110(2):393--411, Apr.
  1984.

\bibitem{Klapper2010}
I.~Klapper and J.~Dockery.
\newblock {Mathematical Description of Microbial Biofilms}.
\newblock {\em SIAM Review}, 52(2):221, Oct. 2010.

\bibitem{Klapper2002}
I.~Klapper, C.~J. Rupp, R.~Cargo, B.~Purvedorj, and P.~Stoodley.
\newblock {Viscoelastic fluid description of bacterial biofilm material
  properties.}
\newblock {\em Biotechnology and bioengineering}, 80(3):289--96, Nov. 2002.

\bibitem{Kreft1998}
J.-U. Kreft, G.~Booth, and J.~W.~T. Wimpenny.
\newblock {BacSim, a simulator for individual-based modelling of bacterial
  colony growth}.
\newblock {\em Microbiology}, 144:3275--3287, 1998.

\bibitem{Kreft2001}
J.-U. Kreft, C.~Picioreanu, J.~W.~T. Wimpenny, and M.~C.~M. van Loosdrecht.
\newblock {Individual-based modelling of biofilms.}
\newblock {\em Microbiology}, 147(Pt 11):2897--2912, Nov. 2001.

\bibitem{Lau2009}
P.~C.~Y. Lau, J.~R. Dutcher, T.~J. Beveridge, and J.~S. Lam.
\newblock {Absolute quantitation of bacterial biofilm adhesion and
  viscoelasticity by microbead force spectroscopy.}
\newblock {\em Biophysical journal}, 96(7):2935--48, Apr. 2009.

\bibitem{LeVeque2007}
R.~J. LeVeque.
\newblock {\em {Finite Difference Methods for Ordinary and Partial Differential
  Equations: Steady-State and Time-Dependent Problems}}.
\newblock SIAM, Philadelphia, PA, 2007.

\bibitem{LiIto2006}
Z.~Li and K.~Ito.
\newblock {\em {The Immersed Interface Method: Numerical Solutions of PDEs
  Involving Interfaces and Irregular Domains}}, volume~33 of {\em Frontiers in
  Applied Mathematics}.
\newblock SIAM, Philadelphia, PA, 2006.

\bibitem{Luo2008}
H.~Luo, R.~Mittal, X.~Zheng, S.~a. Bielamowicz, R.~J. Walsh, and J.~K. Hahn.
\newblock {An immersed-boundary method for flow-structure interaction in
  biological systems with application to phonation.}
\newblock {\em Journal of computational physics}, 227(22):9303--9332, Nov.
  2008.

\bibitem{Masuda1991}
S.~Masuda, Y.~Watanabe, and M.~Ishiguro.
\newblock {Biofilm properties and simultaneous nitrification and
  denitrification in aerobic rotating biological contactors}.
\newblock {\em Water Science and Technology}, 23:1355--1363, 1991.

\bibitem{Mori2008}
Y.~Mori and C.~S. Peskin.
\newblock {Implicit second-order immersed boundary methods with boundary mass}.
\newblock {\em Computer Methods in Applied Mechanics and Engineering},
  197(25-28):2049--2067, Apr. 2008.

\bibitem{Newren2007}
E.~P. Newren, A.~L. Fogelson, R.~D. Guy, and R.~M. Kirby.
\newblock {Unconditionally stable discretizations of the immersed boundary
  equations}.
\newblock {\em Journal of Computational Physics}, 222(2):702--719, Mar. 2007.

\bibitem{Ohashi1996}
A.~Ohashi and H.~Harada.
\newblock {A novel concept for evaluation of biofilm adhesion strength by
  applying tensile force and shear force}.
\newblock {\em Water Science and Technology}, 34(5-6):201--211, 1996.

\bibitem{Pavlovsky2013}
L.~Pavlovsky, J.~G. Younger, and M.~J. Solomon.
\newblock {In situ rheology of Staphylococcus epidermidis bacterial biofilms}.
\newblock {\em Soft Matter}, 9(1):122, 2013.

\bibitem{Peskin1977}
C.~S. Peskin.
\newblock {Numerical analysis of blood flow in the heart}.
\newblock {\em J. Computational Physics}, 81:372--405, 1977.

\bibitem{Peskin2002}
C.~S. Peskin.
\newblock {The immersed boundary method}.
\newblock {\em Acta Numerica}, 11:479--517, July 2002.

\bibitem{Picioreanu2004}
C.~Picioreanu, J.-U. Kreft, and M.~C.~M. van Loosdrecht.
\newblock {Particle-Based Multidimensional Multispecies Biofilm Model}.
\newblock {\em Applied and environmental microbiology}, 70(5):3024--3040, May
  2004.

\bibitem{Picioreanu2001}
C.~Picioreanu, M.~C. van Loosdrecht, and J.~J. Heijnen.
\newblock {Two-dimensional model of biofilm detachment caused by internal
  stress from liquid flow.}
\newblock {\em Biotechnology and bioengineering}, 72(2):205--18, Jan. 2001.

\bibitem{Picioreanu1999}
C.~Picioreanu, M.~C.~M. van Loosdrecht, and J.~J. Heijnen.
\newblock {Discrete-differential modelling of biofilm structure}.
\newblock {\em Water Science and Technology}, 39(7):115--122, 1999.

\bibitem{PicioreanuLoodsdrecht2000DisContMod}
C.~Picioreanu, M.~C.~M. van Loosdrecht, and J.~J. Heijnen.
\newblock {Effect of diffusive and convective substrate transport on biofilm
  structure formation: A two-dimensional modeling study}.
\newblock {\em Biotechnology and Bioengineering}, 69(5):504--515, Sept. 2000.

\bibitem{Rittmann1982}
B.~E. Rittmann.
\newblock {Comparative performance of biofilm reactor types.}
\newblock {\em Biotechnology and bioengineering}, 24(6):1341--70, June 1982.

\bibitem{Rittmann1980}
B.~E. Rittmann and P.~L. McCarty.
\newblock {Evaluation of steady-state-biofilm kinetics}.
\newblock {\em Biotechnology and Bioengineering}, 22(11):2359--2373, Nov. 1980.

\bibitem{Ro1991}
K.~S. Ro and {J. B. Neethling}.
\newblock {Biofilm density for biological fluidized beds}.
\newblock {\em Research Journal of the Water Pollution Control Federation},
  63(5):815--818, 1991.

\bibitem{Rupp2005}
C.~J. Rupp, C.~A. Fux, and P.~Stoodley.
\newblock {Viscoelasticity of Staphylococcus aureus biofilms in response to
  fluid shear allows resistance to detachment and facilitates rolling
  migration.}
\newblock {\em Applied and environmental microbiology}, 71(4):2175--8, Apr.
  2005.

\bibitem{Rusconi2010}
R.~Rusconi, S.~Lecuyer, L.~Guglielmini, and H.~A. Stone.
\newblock {Laminar flow around corners triggers the formation of biofilm
  streamers.}
\newblock {\em Journal of the Royal Society, Interface / the Royal Society},
  7(50):1293--9, Sept. 2010.

\bibitem{Spiga1994}
M.~Spiga and G.~Morino.
\newblock {A symmetric solution for velocity profile in laminar flow through
  rectangular ducts}.
\newblock {\em International Communications in Heat and Mass Transfer},
  21(4):469--475, July 1994.

\bibitem{Stewart}
E.~J. Stewart, A.~Satorius, J.~G. Younger, and M.~J. Solomon.
\newblock {Impact of osmotic stress and sub-lethal antibiotic concentration on
  intercellular spacing and clustering in {{\it Staphylococcus epidermidis}}
  biofilms}.
\newblock {\em Submitted}.

\bibitem{Strychalski2012}
W.~Strychalski and R.~D. Guy.
\newblock {Viscoelastic Immersed Boundary Methods for Zero Reynolds Number
  Flow}.
\newblock {\em Communications in Computational Physics}, 12(2):462--478, 2012.

\bibitem{Todar2012}
K.~G. Todar.
\newblock {Todar's Online Textbook of Bacteriology}.

\bibitem{WangZhang2010Review}
Q.~Wang and T.~Zhang.
\newblock {Review of mathematical models for biofilms}.
\newblock {\em Communication in Solid State Physics}, 150(21-22):1009--1022,
  2010.

\bibitem{Zamir2000}
M.~Zamir.
\newblock {\em {The Physics of Pulsatile Flow}}.
\newblock Biological Physics Series. Springer-Verlag, New York, NY, 2000.

\bibitem{ZhangCoganWang2008One}
T.~Zhang, N.~G. Cogan, and Q.~Wang.
\newblock {Phase Field Models for Biofilms. I. Theory and One-Dimensional
  Simulations}.
\newblock {\em SIAM Journal on Applied Mathematics}, 69(3):641, 2008.

\bibitem{ZhangCoganWang2008Two}
T.~Zhang, N.~G. Cogan, and Q.~Wang.
\newblock {Phase-field models for biofilms II. 2-D numerical simulations of
  biofilm-flow interaction}.
\newblock {\em Commun. Comput. Phys}, 4(1):72--101, 2008.

\bibitem{Zhu2002}
L.~Zhu and C.~S. Peskin.
\newblock {Simulation of a Flapping Flexible Filament in a Flowing Soap Film by
  the Immersed Boundary Method}.
\newblock {\em Journal of Computational Physics}, 179(2):452--468, July 2002.

\bibitem{Zhuo2011}
J.~Zhuo and R.~Dillon.
\newblock {Using the immersed boundary method to model complex fluids-structure
  interaction in sperm motility}.
\newblock {\em Discrete and Continuous Dynamical Systems - Series B},
  15(2):343--355, Dec. 2011.

\end{thebibliography}

\appendix

\section{\label{sec:Appendix}~}

In this appendix we provide a list of variables and parameters used
in this paper.
\begin{description}
\item [{$b$}] Dashpot Damping Coefficient
\item [{$d_{0}$}] Average Spring Rest Length
\item [{$\delta$}] Dirac Delta Function
\item [{$\delta_{h}$}] Discretized Dirac Delta Function from Peskin
\item [{$\tilde{\delta}$}] Our Modified Discretized Dirac Delta Function
\item [{$D$}] Spatial dimension, $D=2$ for 2D simulations and $D=3$
for 3D simulations
\item [{$\mathbf{e}_{i}$}] Unit Vector in the $i^{th}$ direction
\item [{$\eta$}] Total Number of Lagrangian Points
\item [{\textmd{$\mathbf{f}$}}] Eulerian Force Density
\item [{$\mathbf{F}$}] Lagrangian Force
\item [{$F_{max}$}] Maximum Lagrangian Force
\item [{$h$}] Spatial Discretization of finest grid
\item [{$K$}] Hookean Spring Coefficient
\item [{$\mu$}] Dynamic Fluid Viscosity
\item [{$\mu_{max}$}] Maximum Biofilm Viscosity
\item [{$p$}] Pressure
\item [{\textmd{$\mathbf{q}=(q,\, r,\, s)$}}] Lagrangian Coordinates
\item [{$\rho$}] Density
\item [{$\rho_{0}$}] Uniform Fluid Density
\item [{$\rho_{b}$}] Additional Density in Biofilm
\item [{$s$}] Lagrangian Node Marker
\item [{$t$}] Time
\item [{$T$}] Tension in Spring
\item [{$\mathbf{u}$}] Eulerian Velocity
\item [{$\mathbf{\tilde{u}}$}] Intermediate Velocity
\item [{$\mathbf{U}$}] Lagrangian Velocity
\item [{$\mathbf{x}=(x_{1},x_{2},x_{3})$}] Cartesian Coordinates
\end{description}

\section{\label{sec:Scaling-Parameters}Scaling Parameters}

\noindent 
\begin{table}[H]
\centering{}%
\begin{tabular}{|l|l|c|c|}
\hline 
\multicolumn{4}{|c}{Scaling Parameters}\tabularnewline
\hline 
\hline 
Description & Scaling  & Primary  & Specific values \tabularnewline
 & parameter & dimensions & chosen for \tabularnewline
 &  &  & simulations\tabularnewline
\hline 
Characteristic Length & $L$ & $\left\{ L\right\} $ & $50\: microns$\tabularnewline
\hline 
Characteristic Speed & $u_{0}$ & $\left\{ \nicefrac{L}{t}\right\} $ & $10^{-3}\nicefrac{\ensuremath{m}}{s}$\tabularnewline
\hline 
Characteristic Frequency & $T$ & $\left\{ t\right\} $  & $1\, s$\tabularnewline
\hline 
Reference Pressure Difference & $p_{0}-p_{L_{tube}}$ & $\left\{ mL^{-1}t^{-2}\right\} $ & $.8144\, Pa$\tabularnewline
\hline 
Characteristic Density & $\rho_{0}$ & $\left\{ mL^{-3}\right\} $ & $998\,\nicefrac{kg}{m^{3}}$\tabularnewline
\hline 
Characteristic Viscosity & $\mu$ & $\left\{ mL^{-1}t^{-1}\right\} $ & $10^{-3}\,\nicefrac{kg}{ms}$\tabularnewline
\hline 
Characteristic Force Density & $f_{0}$ & $\left\{ \nicefrac{F}{L^{3}}\right\} $ & $varies$\tabularnewline
\hline 
\end{tabular}\\
\caption{\label{tab: 2D sim params-1}Shows the scaling parameters and their
descriptions.}
\end{table}

\end{document}